%% file: witten41.tex
\documentclass[12pt]{amsart}

\usepackage[usenames,dvipsnames]{color}
\usepackage{amsmath,amssymb,amsthm,graphicx,mathrsfs,url,latexsym,enumerate}
\usepackage{a4wide}
\usepackage{color}
\usepackage{ifthen}
\usepackage{verbatim}
\usepackage{fancyhdr}
\usepackage{url}
\usepackage{bbm}
\usepackage[sans]{dsfont}
\usepackage{MnSymbol}
\usepackage{transparent}
\usepackage{amsbsy}
\usepackage[justification=centering]{caption}
\usepackage[colorlinks=true,linkcolor=Red,citecolor=Green]{hyperref}
\usepackage{wrapfig}

\newcommand{\dist}{\operatorname{dist}}
\newcommand{\Id}{\operatorname{Id}}

\newcommand{\E}{\mathbb E}

\newcommand{\supp}{\operatorname{supp}}

\newcommand{\diag}{\operatorname{diag}}
\newcommand{\C}{\mathbb C}
\newcommand{\CC}{\mathscr{C}}
\newcommand{\N}{\mathbb N}
\newcommand{\R}{\mathbb R}

\newcommand{\Hess}{\operatorname{Hess}}
\newcommand{\Cl}{\operatorname{Cl}}
\newcommand{\Ran}{\operatorname{Ran}}
\newcommand{\Ker}{\operatorname{Ker}}
\newcommand{\bsigma}{\boldsymbol{\sigma}}

\def\<{\langle}
\def\>{\rangle}
\newcommand{\bp}{{\it Proof. }}
\newcommand{\ep}{\hfill $\square$\\}

\def\m{\mathbf{m}}
\def\s{\mathbf{s}}

\def\b{\mathbf{b}}
\def\i{\mathbf{i}}

\newcommand{\be}{\begin{equation}}
\newcommand{\ee}{\end{equation}}
\newcommand{\bes}{\begin{equation*}}
\newcommand{\ees}{\end{equation*}}

\renewcommand{\dh}{d_{\phi,h}}

\newcommand{\dhs}{d_{\phi,h}^{*}}

\numberwithin{equation}{section}
\numberwithin{figure}{section}

\newtheorem{defin}{Definition}[section]
\newtheorem{theorem}[defin]{Theorem}
\newtheorem{corollary}[defin]{Corollary}
\newtheorem{lemma}[defin]{Lemma}

\newtheorem{proposition}[defin]{Proposition}
\newtheorem{remark}[defin]{Remark}
\newtheorem{example}[defin]{Example}

\def\aaa{{\mathcal A}}\def\bbb{{\mathcal B}}\def\ccc{{\mathcal C}}\def\ddd{{\mathcal D}}
\def\eee{{\mathcal E}}\def\fff{{\mathcal F}}\def\ggg{{\mathcal G}} 
\def\iii{{\mathcal I}} \def\jjj{{\mathcal J}}\def\lll{{\mathcal L}}
\def\mmm{{\mathcal M}}\def\nnn{{\mathcal N}} \def\ooo{{\mathcal O}}\def\ppp{{\mathcal P}}
\def\rrr{{\mathcal R}}\def\sss{{\mathcal S}}
\def\uuu{{\mathcal U}}\def\vvv{{\mathcal V}}
  \def\zzz{{\mathcal Z}}

\def\Ar{{\mathscr A}}\def\Cr{{\mathscr C}}
\def\Er{{\mathscr E}}\def\Fr{{\mathscr F}}\def\Gr{{\mathscr G}}
\def\Jr{{\mathscr J}}\def\Lr{{\mathscr L}}
\def\Mr{{\mathscr M}}\def\Pr{{\mathscr P}}
\def\Sr{{\mathscr S}}\def\Tr{{\mathscr T}}
\def\Wr{{\mathscr W}}
\def\Yr{{\mathscr Y}}

\def\ulu{\underline{\mathcal U}}
\def\ala{\underline{\mathcal A}}

\newcommand{\puu}{\wideparen\uuu}

\begin{document}
\title{About small eigenvalues of Witten Laplacian }

\author[L. Michel]{Laurent~Michel}
\address{Universit\'e de Bordeaux, Institut Math\'ematique de Bordeaux, France}
\email{laurent.michel@\allowbreak math.\allowbreak u-bordeaux.fr}

\begin{abstract}
We study the low lying eigenvalues of the semiclassical Witten Laplacian $\Delta_\phi$ associated with a Morse function $\phi$. 
We consider the case where the sequence of Arrhenius numbers $S_1\leq \ldots\leq S_n$ associated with $\phi$ is degenerated, that is the preceding inequalities are not necessarily strict. 
\end{abstract}

\maketitle
\tableofcontents
\section{Introduction}\label{sec:Intro}
\subsection{Motivations}\label{subsec:Motiv}
Introduced in the early eighties by E. Witten to give an analytical proof of Morse inequalities \cite{Wi82}, Witten's Laplacian appeared recently as the cornerstone of many quantitative studies of metastability for kinetic equations (see e.g. \cite{HeNi04,HeKlNi04_01,HeHiSj11_01,DiGLeLLePNec18}).
One  of the simplest examples of metastable dynamics is given by the movement of particle evolving in an energy landscape $\phi$ and submitted to random forces. The position $X_t$ of such a particle at time $t$ satisfies the over-damped Langevin equation
\begin{equation}\label{eq:langevin}
 \dot X_t=-2\nabla \phi(X_t)+\sqrt{2h}\dot B_t
\end{equation}
where  $h$ stands for the temperature of the system and $B_t$ is a Brownian force.\footnote{This equation appears for instance in physics to describe the microscopic 
evolution of a charged gas assuming the mass of the particles is negligible.}
Assuming that the potential $\phi$ has several wells 
a particle starting from a local minimum of the function $\phi$ can use the random force to jump over a saddle point and reach another energy well.
The famous  Eyring-Kramers law describes the 
average time to escape from a well  in the regime of low temperature $h\rightarrow 0$. 
In his pioneering work  \cite{Kr40}, Kramers  predicted in the simplified one dimensional setting  that the average transition 
time $\tau_\phi$ from a local minimum $m$ to the nearest saddle point $s$ is exponentially large with respect to $h^{-1}$, $\tau_\phi\simeq a_\phi e^{ \kappa_\phi/h}$ and he gave additionally formulae for the positive
coefficient $\kappa_\phi$  and the prefactor $a_\phi$ in terms of the second derivatives of $\phi$ in points $m$ and $s$.
Observe that when $h\rightarrow 0$, this average transition time becomes extremely large, which justifies the terminology of metastable state given to the position $m$.  

In practice, Eyring -Kramers law has very important applications in many domains of science where the trajectory \eqref{eq:langevin} is used to implement computational algorithms. In order to compute some thermodynamical quantities 
\begin{equation}\label{eq:thermo}
\E_\mu(f)=\int_{\R^d}f(x)d\mu(x)
\end{equation}
associated with a measure  $\mu$ and an observable $f$, one can introduce any random dynamics $X_t$ ergodic with respect to $\mu$ and  use 
Monte Carlo method to approximate $\E_\mu(f)$ by the  long time average of $f$ along any trajectory (see \cite{LeRoSt10_01} for introduction to these questions).
In many situations one has $d\mu(x)=Z_he^{-\phi(x)/h}$ for some potential $\phi$ and the over-damped Langevin dynamics \eqref{eq:langevin} fulfills the necessary assumptions.
The time needed by the process $X_t$  to explore the whole space $\R^d$ (which insures the validity of the above approximation), is directly linked to the metastable properties discussed previously. Understanding this metastable behavior is then of crucial interest if one needs to evaluate some stopping time or accelerate the convergence for instance. \\

 
 The first mathematical proof of Eyring-Kramers law in a general setting was  obtained recently by a potential theory approach  \cite{BoGaKl05_01} and next by semiclassical methods  \cite{HeKlNi04_01}.
 This later approach and the link with the Witten Laplacian, can be understood easily by looking at Langevin equation \eqref{eq:langevin} at the macroscopic level. Indeed, the evolution 
of any statistical distribution $\rho(t,x)$ of particles is governed by the Kramers-Smoluchowski equation
\be\label{eq:krasmo1}
\partial_t\rho-h\Delta\rho-2\operatorname{div}(\rho\nabla\phi)=0.
\ee
Writing $\rho=e^{-\phi/h}\tilde\rho$, the above equation is equivalent to
$h\partial_t\tilde\rho+\Delta_\phi\tilde\rho=0$, where 
$$
\Delta_\phi=-h^2\Delta+\vert\nabla\phi\vert^2-h\Delta\phi
$$
is the semiclassical Witten Laplacian associated to $\phi$.
This operator is known to be non-negative and under confining assumptions on the function $\phi$ it has  a non-trivial kernel corresponding to the global equilibrium of 
\eqref{eq:krasmo1}. 
As a consequence, the behavior of $\tilde\rho$ when $t\rightarrow\infty$ is driven by the small eigenvalues of 
$\Delta_\phi$. In particular, any state associated to a small eigenvalue looks stable during very long times. 
These are the  metastable states and the inverses of the corresponding eigenvalues yield their lifetimes. 
In \cite{HeKlNi04_01}, Helffer-Klein-Nier obtained a full description of the small eigenvalues of the Witten Laplacian in a quite general setting.
In terms  of Kramers-Smoluchovski equation, their result implies that if the initial probability distribution $\rho_0$ belongs to 
 $L^2(e^{2\phi/h}dx)$, then the solution $\rho$ of \eqref{eq:krasmo1} converges exponentially fast to the equilibrium probability distribution $c_h^{-2}e^{-2\phi/h}$  (where $c_h$ is a normalizing factor) 
\be\label{eq:conv-krasmo}
 \Vert \rho(t)-\frac 1{c_h^2}e^{-2\phi/h}\Vert_{L^2 (e^{2\phi/h}dx)}\leq e^{-\lambda_h t}\Vert \rho_0\Vert_{L^2 (e^{2\phi/h}dx)}.
 \ee
Moreover, the rate of convergence $\lambda_h=hb(h)e^{-2S/h}$ is described by the Eyring-Kramers law, that is:
 \begin{itemize}
 \item[-] $S$ is the biggest height a particle has to pass in order to reach the unique global minimum.
 \item[-] the prefactor $b(h)$ has an asymptotic expansion with respect to the parameter $h$, $b(h)\sim\sum_k b_kh^k$ and its leading term is given by an explicit formula in terms 
 of the Hessian of $\phi$.
 \end{itemize}
 More precisely, the assumptions made in \cite{HeKlNi04_01} imply that there exist a unique minimum $m$ and a unique saddle point $s$ of $\phi$ such that $S=\phi(s)-\phi(m)$. Then, the leading term of $b(h)$ is 
 \be\label{eq:EyKrHeKlNi}
 b_0=\frac{\vert \mu_1(s)\vert}{\pi}\sqrt{\frac{\det\Hess(\phi)(m)}{\det\Hess(\phi)(s)}}
 \ee
 where $\mu_1(s)$ denotes the negative eigenvalue of $\Hess(\phi)(s)$. In the case of a  double well, this formula is exactly the one predicted by Kramers in his paper \cite{Kr40}. 
Later on, the method developed in \cite{HeKlNi04_01} to compute the small eigenvalues of the  Witten Laplacian  was successfully 
used on bounded domains  \cite{HeNi06}, \cite{Lep10} and in the study of semiclassical random walks \cite{BoHeMi15}. 

However, the range of potential $\phi$ covered by these papers  doesn't include many cases  which are very important in practice. 
Roughly speaking, in \cite{HeKlNi04_01}, the author make an assumption on the relative position of minima and saddle points that insures in fine that the small eigenvalues are all of different size. Among the limitations of this assumption, the potential $\phi$ may not have too much saddle points or minima at the same level. 

It turns out that in many physical applications, the energy landscape may not satisfy the above assumption.
In chemical physics, the energy potential of the reaction may have some symmetries in numerous situations.
This is for instance the case when looking at some homogenous systems such as Lennard-Jones clusters (see \cite{Wa06} for example and discussions on this topics). \\
%

The aim of this paper is to study the spectral properties of $\Delta_\phi$ in 
the case where $\phi$ is a general Morse function without restriction on the relative position of its minima.

\subsection{Heuristics on a simple example}\label{subsec:HOASE}
\begin{figure}[!h]
 \center
  \scalebox{0.35}{ 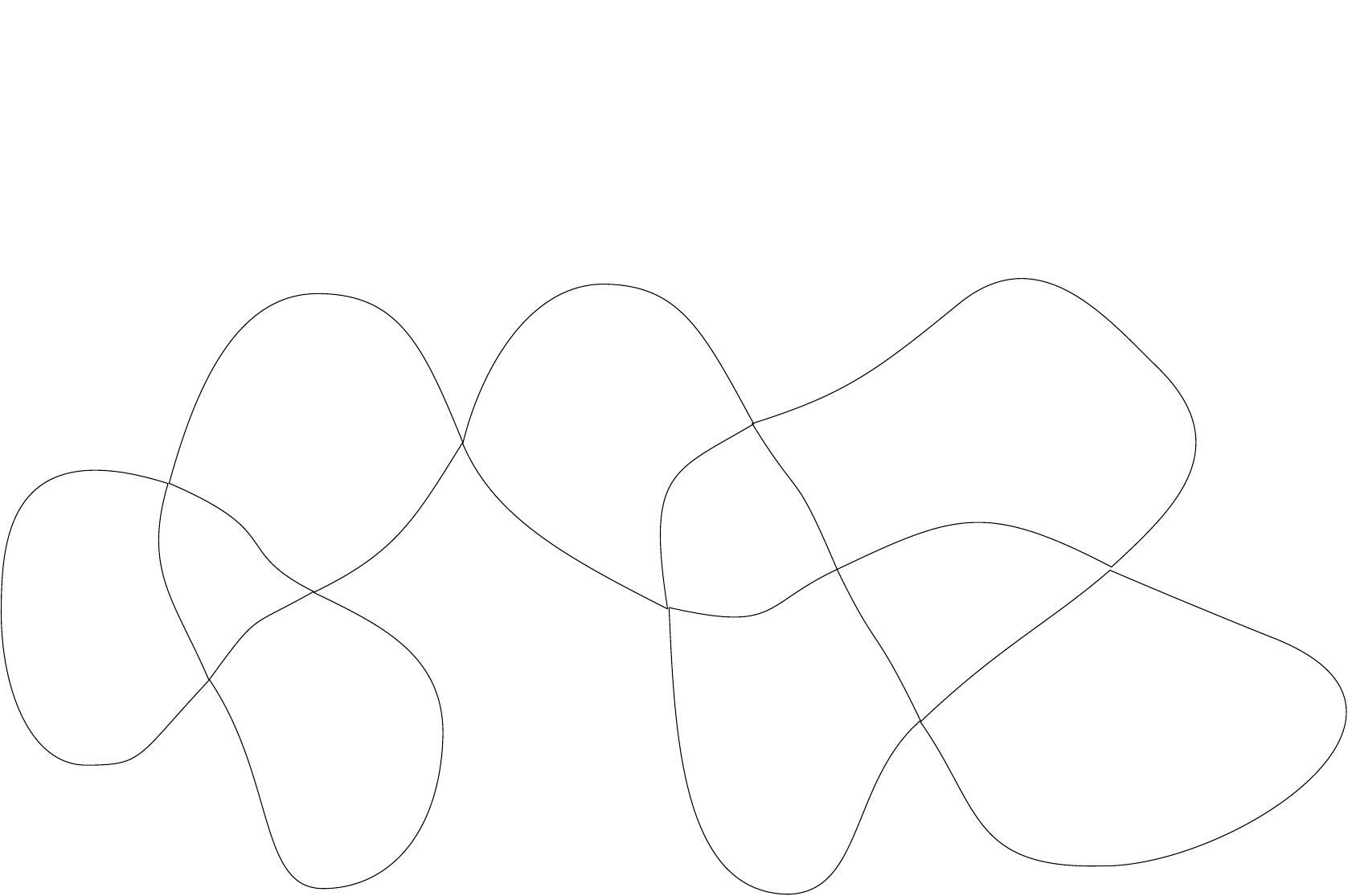 } \hspace{0.1in}
\scalebox{0.6}{ 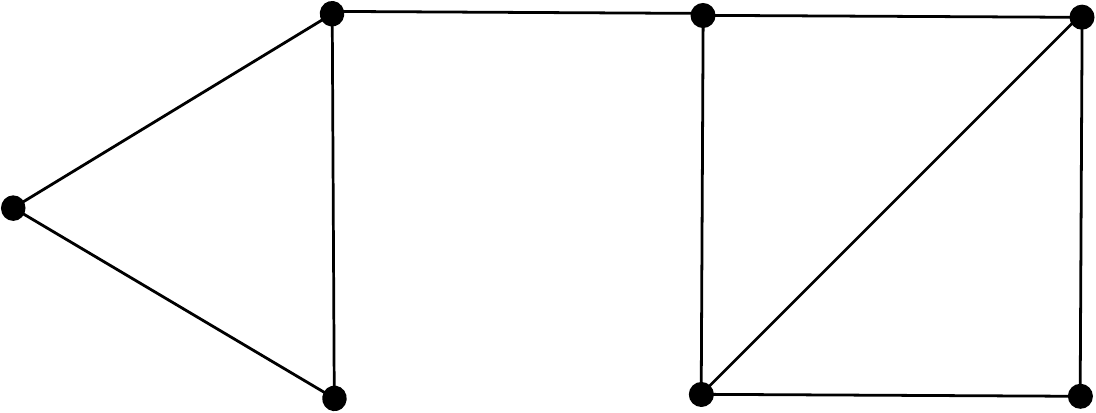}
  \caption{{\em Left:} The sublevel set $\{\phi<\sigma\}$ (dashed region) associated to a potential $\phi$ having a unique saddle value $\sigma$.
   The x's represent local minima, the o's, local maxima. {\em Right:} The graph associated to the
potential on the left.}
    \label{figsublevel}
\end{figure}

A typical example we have in mind is  the following. Suppose that $\phi:\R^d\rightarrow\R$ has $n_0$ minima all at the same level and $n_1$ saddle points all at the same level (see Figure \ref{figsublevel}, where the x represent minima and the o local maxima). Denote by $S=\phi(s)-\phi(m)$ the difference of height between any minimum and any saddle point. In order to simplify the situation assume also that the function 
$ \Hess(\phi)(x)$ has eigenvalues $\pm1$ when $x$ belongs to the set of minima and saddle points.

Such an example doesn't enter in the framework of \cite{HeKlNi04_01}, however it develops some very interesting behavior. 
More precisely, in the very simplified case discussed in this section,  a consequence of Theorem \ref{thm:spectreDtypeII} below is the following
\begin{theorem}\label{th:simplifie}
There exists $\epsilon_0>0$ and $h_0>0$ such that for all $h\in]0,h_0]$, $\Delta_\phi$ has exactly $n_0$ eigenvalues 
$\lambda_k, k=1,\ldots,n_0$ in the interval $[0,\epsilon_0 h]$.
The lowest eigenvalue is $\lambda_1=0$ and the other ones have the form
$$
\lambda_k=hb_k(h)e^{-2S/h}
$$
for all $k=2,\ldots,n_0$. Moreover, the prefactors $b_k(h)$ have an classical expansion $b_k(h)\sim \sum_{j=0}^\infty h^jb_{k,j}$ and the leading terms $b_{k,0}$ are exactly the non zero eigenvalues of the graph $\ggg$ whose vertices are the minima of $\phi$ and whose edges are the saddle points joining two minima (see Figure \ref{figsublevel}).
\end{theorem}

In terms of Kramers-Smoluchovski equation, this theorem  exhibits some metastable states whose lifetimes (given by the inverse of the above eigenvalues) are quantified by the graph $\ggg$.
 At the level of particle, these new rules of computation can be understood as follows.
 Since all the minima are at the same level, the equilibrium state is equidistributed among all the minima. Moreover, since all the saddle points are at the same level,  an ergodic trajectory
 of \eqref{eq:langevin}  will visit all the minima in the same time scale, by traveling  along the edges of the graph $\ggg$. The same graph Laplacian was constructed by Landim et al \cite{LaMiTs15} in a discrete setting.
 
 It could look surprising that the coefficients $b_k$ do not depend on the second derivative of $\phi$ as in the usual case of 
 Eyring-Kramers formula. This is actually due to the normalization assumption made at the beginning of this section. In the case where the Hessian is arbitrary, the above result is still available with a weighted graph instead of $\ggg$. The weights depend explicitly on the Hessian of $\phi$ (see Theorem \ref{thm:spectreDtypeII}).
 \\

\begin{figure}[!h]
 \center
  \scalebox{0.35}{ 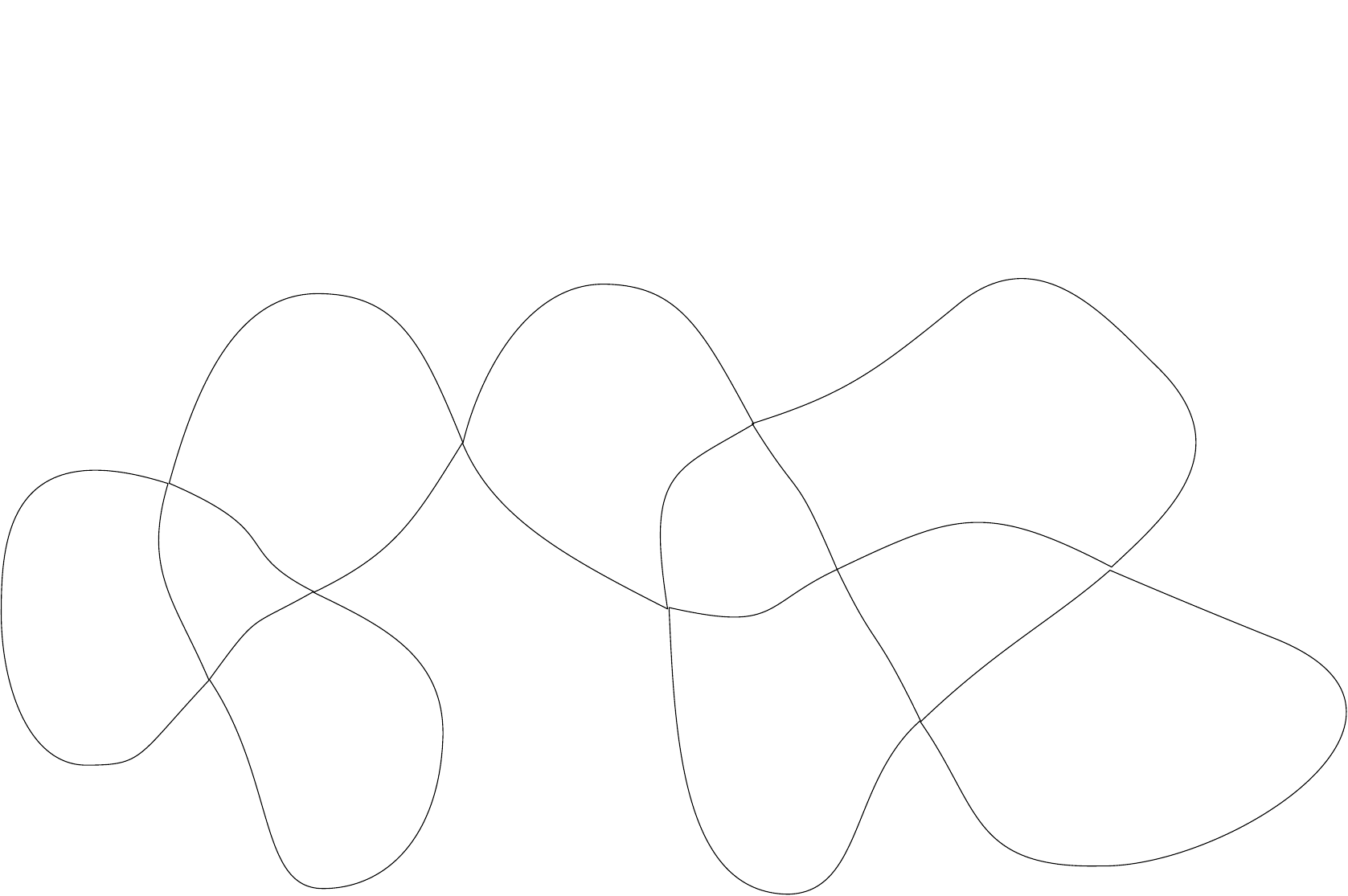 } \hspace{0.1in}
\scalebox{0.6}{ 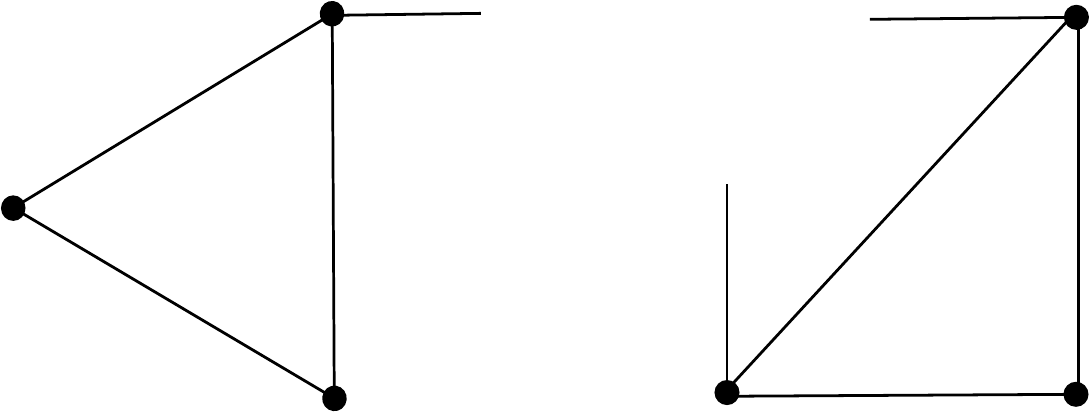}
  \caption{{\em Left:} The sublevel set $\{\phi<\sigma\}$ (dashed region) associated to a potential $\phi$ having a unique saddle value $\sigma$.
   The x's represent local minima, the o's, local maxima. {\em Right:} The two hypergraphs associated to the
potential on the left (the missing vertex corresponds to the minimum $A$)}
    \label{figsublevel-disjoint}
\end{figure}
 To go further and motivate the object introduced in the next section, let us discuss now what happens if we modify slightly the potential $\phi$ in the following way. Suppose that $\phi$ is still like in Figure \ref{figsublevel} excepted that we modify one of the minimal values (higher or lower). In Figure \ref{figsublevel-disjoint},  the modified minimum is denoted by  $A$.
Then, we can associate to this potential the two hypergraphs corresponding to minima at the same level and linked by a saddle value (see  figure \ref{figsublevel-disjoint}). 
If $A$ is an absolute minimum, then equilibrium distribution is concentrated in $A$ and the prefactor $b_k(h)$ will be given by the smallest non zero eigenvalue of the two hypergraphs introduced above (roughly speaking this represents the maximum time needed to reach $A$).
In the contrary, if $A$ is not a global minimum anymore, the equilibrium state is uniformly distributed among all the minima excepted $A$.
In order to visit each site of the equilibrium state, an arbitrary particle will necessarily pass through the point $A$.
This heuristics explains why the computation of the prefactor $b_k(h)$ will involve a more complicated procedure describing the interaction between the two hypergraphs  via the well $A$. 

The main contribution of this paper is to describe quantitatively this phenomena.

\section{Framework and results}\label{sec:FAR}

Let $X$ be either $\R^d$ or a compact manifold of dimension $d$ without boundary and let $\phi:X\rightarrow \R$ be a smooth Morse function.
Consider the semiclassical Witten Laplacian associated to $\phi$:
\begin{equation}\label{eq:definWitten0form}
\Delta_\phi=-h^2\Delta+\vert\nabla\phi\vert^2-h\Delta\phi
\end{equation}
where $h\in]0,1]$ denotes the semiclassical parameter.

If $X$ is a compact manifold, the operator $\Delta_\phi$ is selfadjoint with domain $H^2(X)$ and its resolvent is compact. In the case $X=\R^d$ we make the  additional assumption that 
there exist $C>0$ and a compact $K\subset\R^d$ such that for all $x\in\R^d\setminus K$, one has
\begin{equation}\label{eq:hypgenephi}
 \vert\nabla\phi(x)\vert\geq\frac 1 C,\;\vert \Hess(\phi(x))\vert\leq C\vert\nabla\phi\vert^2,\text{ and }\phi(x)\geq C\vert x\vert.
\end{equation}
Then, $\Delta_\phi$ is essentially selfadjoint on $\ccc^\infty_c(\R^d)$ and thanks to \eqref{eq:hypgenephi}, there exist $h_0>0$ and $c_0>0$ such that
for all $h\in]0,h_0]$, we have 
$$
\sigma_{ess}(\Delta_\phi)\subset [c_0,\infty[.
$$ 
In both situations $X$ compact or $X=\R^d$, it is well-known that $\Delta_\phi$ is non-negative. Hence
$\sigma(\Delta_\phi)\subset[0,\infty[$ and it follows from the above remarks that  $\sigma(\Delta_\phi)\cap[0,c_0[$ is made of  eigenvalues with no accumulation point excepted maybe $c_0$. Moreover 
$e^{-\phi/h}$ is clearly in the kernel of $\Delta_\phi$ and belongs to $L^2(\R^d)$ thanks to \eqref{eq:hypgenephi}, so that  the lowest eigenvalue of
$\Delta_\phi$ is clearly $0$.

Since $\phi$ is a Morse function (and thanks to assumption \eqref{eq:hypgenephi} in the case $X=\R^d$), the set $\uuu$ of  critical points is finite. In the following, for $p=0,\ldots,d$, we will denote by $\uuu^{(p)}$ the set of critical points of $\phi$ of index $p$. Hence, $\uuu^{(0)}$ is the set of minima and 
$\uuu^{(1)}$ the set of saddle points of $\phi$.  Throughout the paper, we will denote $n_j=\sharp\uuu^{(j)}$.

From the pioneer work by Witten \cite{Wi82}, it is well-known that for small $h$, there is a correspondance between the  small eigenvalues of $\Delta_\phi$ and the critical points of $\phi$.
More precisely, by standard localization arguments one can show that there exists  $\epsilon_0>0$ such that for $h>0$ small enough, $\Delta_\phi$ has exactly $n_0$ eigenvalues in the interval
 $[0,\epsilon_0 h]$ that we denote $0=\lambda_1\leq \lambda_2\leq\ldots\leq \lambda_{n_0}$. This result is easily proved in 
 \cite{CyFrKiSi87_01} with $\epsilon_0 h$ replaced by $h^{3/2}$. The proof 
 with $\epsilon_0 h$ can be found in \cite{HeSj85_01}, Prop. 1.7 (see also Prop. 1 of \cite{MiZw18} for a self-contained proof).   Moreover, these eigenvalues are actually exponentially small, that is 
live in an interval $[0,e^{-C/h}]$ for some $C>0$ (see \cite{He88_01} for a proof).
From a  topological point of view, these informations (together with the equivalent estimates for the Witten Laplacian $\Delta_\phi^{(p)}$ acting on $p$-forms) are sufficient to establish a correspondance between the small eigenvalues of $\Delta_\phi^{(p)}$ and the critical points of $\phi$ of index $p$
(this was the key point in the Witten's proof of Morse inequalities).
However, for applications to the description of metastable dynamics, it is important to get some accurate description of the $\lambda_j$'s. Our main theorem will give some asymptotic of these eigenvalues for any Morse function $\phi$, without any assumption on the relative position of minimal and saddle values of $\phi$.\\

%

Before going further, let us introduce few notations that we use in this paper.
For $x_0\in X$ and $r>0$, introduce the geodesic ball $B(x_0,r)=\{x\in X,\,d(x,x_0)<r\}$.

Throughout, we will say that $\s$ is a saddle point if it is a critical point of index $1$. 

Given $a(h),b(h)>0$, two functions of the semiclassical parameter, we say that $a(h)\asymp b(h)$ if there exists some constant
$c_1,c_2>0$ such that for all $h>0$ small one has $c_1b(h)\leq a(h)\leq c_2b(h)$. 
We say that a family of vectors $(a(h))_{h\in]0,1]}$ in a normed vector space $V$ admits a classical expansion if there exists a sequence 
of vectors $(a_n)_{n\in\N}$ independent of $h$ and such that for all $N\in\N$, there exists  some constants $C_N>0$ such that
$$
\Vert a(h)-\sum_{n=0}^Nh^n a_n\Vert_V\leq C_Nh^{N+1},\;\forall h\in]0,1].
$$
We denote $a(h)\sim\sum_{n=0}^\infty h^n a_n$.

As we shall see later, we will have to analyze carefully some finite dimensional matrices which are strongly related to the critical points of
$\phi$. Given any subsets $\bbb_1,\bbb_2$ of $\uuu$, it will be convenient to introduce the finite dimensional  vector space $\Fr(\bbb_j)$
of real valued functions on $\bbb_j$. We shall then denote by $\Mr(\bbb_1,\bbb_2)$ the vector space of linear operators from 
$\Fr(\bbb_1)$ into $\Fr(\bbb_2)$.

\subsection{Labelling of minima}\label{sec:LOM}
Let us now recall the general labelling of minima introduced in \cite{HeKlNi04_01} and generalized in \cite{HeHiSj11_01}.
The main ingredient is the notion of separating saddle point which is defined as follows. 
Given a saddle point $\s$ of $\phi$, and $r>0$ small enough, the set 
$$\{x\in B(\s,r),\;\phi(x)<\phi(\s)\}$$
has exactly two connected components $C_j(\s,r)$, $j=1,2$. The following definition is taken from \cite{HeHiSj11_01}, Definition 4.1.
\begin{defin}\label{def:SSP}
We say that $\s\in X$ is a separating saddle point (ssp) if it is a saddle point and if $C_1(\s,r)$ and $C_2(\s,r)$ 
are contained in two different connected components of 
$\{x\in X,\,\phi(x)<\phi(\s)\}$. We will denote by $\vvv^{(1)}$ the set of ssp.

We say that $\sigma\in\R$ is a separating saddle value (ssv) if it is of the form 
$\sigma=\phi(\s)$ with $\s\in\vvv^{(1)}$. We denote $\underline\Sigma=\phi(\vvv^{(1)})$ the set of ssv.

We say that $E\subset X$ is a critical component if there exists 
$\sigma\in\underline\Sigma$ such that $E$ is a connected component of $\{\phi<\sigma\}$ and if 
$\partial E\cap\vvv^{(1)}\neq \emptyset$. We denote by $\Cr$ the set of critical components.
\end{defin}
\begin{figure}
 \center
  \scalebox{0.6}{ 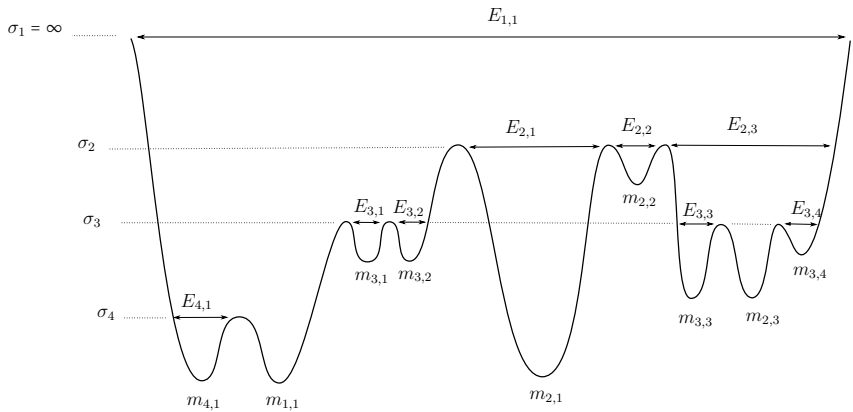}
  \caption{Labelling procedure}
    \label{fig1}
\end{figure}
Let us now describe the labelling procedure of  \cite{HeHiSj11_01}. Since $\phi$ is a Morse function, it has finitely many critical points and so $\underline\Sigma$ is finite. We denote $\sigma_2>\sigma_3>\ldots>\sigma_N$ its elements and for convenience we also introduce a fictive infinite saddle value $\sigma_1=+\infty$ and denote $\Sigma=\underline\Sigma\cup\{\sigma_1\}$.
Starting from $\sigma_1$, we will recursively associate to each $\sigma_i$ a finite family of local minima 
$(\m_{i,j})_j$ and a finite family of critical components $(E_{i,j})_j$ (see Figure \ref{fig1}).

\begin{itemize}
\item[-] Let $X_{\sigma_1}=\{x\in X,\,\phi(x)<\sigma_1=\infty\}=X$.
We let $\m_{1,1}$ be any global minimum of $\phi$ (not necessarily unique) and $E_{1,1}=X$.
\item[-] Next we consider $X_{\sigma_2}=\{x\in X,\,\phi(x)<\sigma_2\}$. This is the union of its finitely many connected components. Exactly one of these components contains $\m_{1,1}$ and the other components are denoted by $E_{2,1},\ldots,E_{2,N_2}$. In each component $E_{2,j}$, we pick up a point $\m_{2,j}$ which is a  global minimum of $\phi_{\vert E_{2,j}}$.
\item[-] Suppose now that the families $(\m_{k,j})_j$ and $(E_{k,j})_j$ have been constructed until rank $k=i-1$. The set 
$X_{\sigma_i}=\{x\in X,\phi(x)<\sigma_i\}$ has again finitely many connected components and we label $E_{i,j},\,j=1,\ldots,N_i$ those of these components that do not contain any $\m_{k,l}$ with $k<i$. In each  $E_{i,j}$ we pick up a point $\m_{i,j}$ which is a  global minimum of $\phi_{\vert E_{i,j}}$. Observe that for all $i\geq 2$, the components $E_{i,j}$ are all critical.
\end{itemize}
We run the procedure  until all the minima have been labelled.
\begin{remark}\label{rem:repart-min}
The above labelling satisfies the following property. For any $\sigma_i\in \Sigma$ and any connected component $A_i$ of 
$\{\phi<\sigma_i\}$, there exists a unique $(k,l)$ such that $k\leq i$ and $\m_{k,l}\in A_i$.
\end{remark}
\bp
Let us start with the existence part of the result. If $A_i$ is one of the $E_{i,j}$ for some $j$, then take $k=i$ and $l=j$. Otherwise, this means that in the labelling procedure, $A_i$ already contained a minimum $\m_{k,l}$ with $k< i$.

Let us prove the uniqueness part. Assume that $\m_{k,l}, \m_{k',l'}\in A_i$ with $k\leq k'\leq i$.
Then $A_i\cap E_{k',l'}\neq\emptyset$ and since $A_i$ is a connected component of $\{\phi<\sigma_i\}$ with $\sigma_i\leq\sigma_{k'}$ it follows that
$A_i\subset E_{k',l'}$. Since $\m_{k,l}\in A_i$, it follows that $\m_{k,l}\in E_{k',l'}$  which is impossible unless $(k,l)=(k',l')$.
\ep

Using the above labelling, H\'erau-Hitrik-Sj\"ostrand made  some significative progress in \cite{HeHiSj11_01} (in the more general situation of Kramers-Fokker-Planck operators, but this applies to Witten Laplacian). 
First, they showed in Theorem 7.1 that the exponentially small eigenvalues $(\lambda_\m(h))_{\m\in\uuu^{(0)}}$ of $\Delta_\phi$ (indexed by the sequence of local minima) satisfy $\lambda_\m(h)\asymp h e^{-2S(\m)/h}$ for the sequence of Arrhenius numbers $(S(\m))_{\m\in\uuu^{(0)}}$ defined by
$S(\m_{i,j})=\sigma_i-f(\m_{i,j})$ with the above notations. However, their method does not permit to prove that $h^{-1}\lambda_\m(h)e^{2S(\m)/h}$ admits a limit when $h\rightarrow 0$.
In order to compute the asymptotic expansion of the eigenvalues $\lambda_\m(h)$, they need to make some additional assumption on the interaction between minima and saddle points (see Assumption 5.1 in \cite{HeHiSj11_01}). This hypothesis, which is a generalization  of  the one made in \cite{HeKlNi04_01}, can be formulated as follows with the  notations of the preceding section:

\vspace{0.2cm}
\emph{Generic Assumption:} For all $i=1,\ldots,N$, $j=1,\ldots,N_i$, the following hold true:
\begin{itemize}
\item[i)]  $\m_{i,j}$ is the unique global minimum of the application $\phi_{\vert E_{i,j}}$.
\item[ii)]  if $E$ is a connected component of $\{\phi<\sigma_i\}$ such  that   $E\cap\vvv^{(1)}\neq\emptyset$, there exists a unique $\s\in\vvv^{(1)}$ such that 
$\phi(\s)=\sup \phi(E\cap\vvv^{(1)})$. In particular, $\phi^{-1}(]-\infty,\phi(\s)[)\cap E$ is the union of exactly two different connected components.
\end{itemize}
Throughout the paper, we denote by (GA) this assumption.
\vspace{0.2cm}

Under this assumption, there exists a bijection between $\uuu^{(0)}$ and $\vvv^{(1)}\cup\{\s_1\}$ where
$\s_1$ is a fictive saddle point associated to $\sigma_1=\infty$ and  for which by convention $\phi(\s_1)=\infty$.
Using this one to one correspondence, the authors exhibit some labelling
$\uuu^{(0)}=\{\m_1,\ldots,\m_{n_0}\}$ and $\vvv^{(1)}\cup\{\s_1\}=\{\s_1,\ldots,\s_{n_0}\}$ such that the small eigenvalues
$\lambda_i(h)$ are of the form
$hb_i(h)e^{-2S_i/h}$ with $S_i=\phi(\s_i)-\phi(\m_i)$.
Moreover, they prove that the  $b_i(h)$ have a classical expansion and compute the leading term of this expansion (see Theorem 5.10 in \cite{HeHiSj11_01}).

As it is stated above, (GA) is not exactly Assumption 5.1 stated in \cite{HeHiSj11_01}. Indeed, it is supposed in \cite{HeHiSj11_01} that  ii) holds true only for $E$ being a critical component. 
However, as indicated by the anonymous referee,  we can easily construct some function $\phi$ satisfying this assumption, for which there is no bijection between $\uuu^{(0)}$ and $\vvv^{(1)}$. To see this, first consider in dimension $1$ a potential $\phi$ with $4$ minima $m_j$, $j=1,\ldots,4$ and $3$ saddle points $s_j$, $j=1,\ldots,3$  such that $m_1<s_1<m_2<s_2<m_3<s_3<m_4$ and such that 
$\phi(m_1)<\phi(m_4)<\phi(m_2)=\phi(m_3)$ and $\phi(s_1)=\phi(s_2)<\phi(s_3)$.
Since the component of $\{\phi<\phi(s_3)\}$ containing $m_1$ is not critical, this function satisfies Assumption 5.1 in \cite{HeHiSj11_01}. It doesn't satisfy  (GA) as stated above. In higher dimension, one can easily generalize this construction to obtain potentials satisfying Assumption 5.1 in \cite{HeHiSj11_01}, with a fixed number of minima and an arbitrary large number of separating saddle points (think for instance to many saddle points between the well containing $m_1$ and the well containing $m_2$). This shows that Assumption 5.1 is not sufficient to insure a  bijection between minima and separating saddle points.

Let us emphasize that the above remark doesn't affect the rest of the work done in \cite{HeHiSj11_01}, where we can easily use the above corrected version of Assumption 5.1 .

Let us observe that the Generic Assumption allows some degeneracy in the sequence $(S_j)$, that is there 
may exists $j$ such that $S_j=S_{j+1}$. 
However, (GA) remains restrictive for the following reasons:
\begin{itemize}
\item[-] It permits only potentials $\phi$ for which  $\uuu^{(0)}$ and $\vvv^{(1)}\cup\{\s_1\}$ have the same cardinal.
\item[-] The eventual degenerate heights are associated to weakly interacting eigenstates in the following sense. Assume for instance that 
$S_j=S_{j+1}$ for some $j=1,\ldots,n_0-1$ and modify slightly the function $\phi$ near the minimum $\m_j$. Then the coefficients $b_j$ is modified whereas  the classical expansion of $b_{j+1}$ remains unchanged.
\end{itemize}
Figures   \ref{fig3}, \ref{fig4} below present some examples of potentials where (GA) is not satisfied. 
These examples as well as an example in higher dimension are discussed in detail in section \ref{subsec:SE}.

In the present  paper, we obtain an asymptotic expansion for the $\lambda_i(h)$ for general Morse functions $\phi$ without any additional assumptions
on the relative position of minima and ssp's.

\subsection{Main result}\label{subsec:mainres}
\begin{figure}
 \center
  \scalebox{0.8}{ 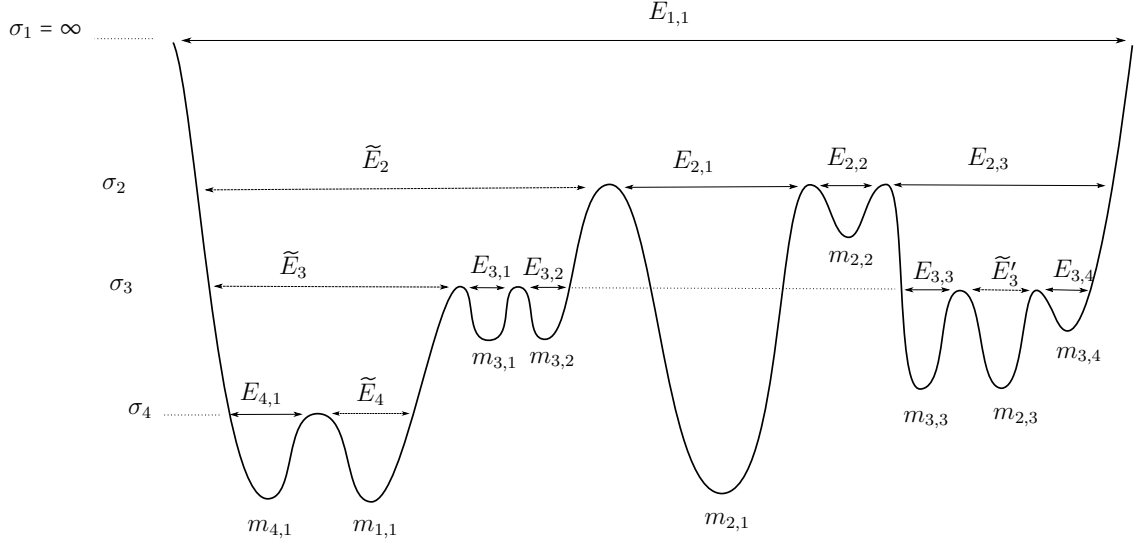}
  \caption{Computations of Example \ref{example1} }
    \label{fig2}
\end{figure}

In order to state our main result, we introduce few notations that will be used throughout the paper.
First, using the above labelling, we define 
 $\bsigma:\uuu^{(0)}\rightarrow \Sigma$ by $\bsigma(\m_{i,j})=\sigma_i$ 
 and  $S:\uuu^{(0)}\rightarrow ]0,+\infty]$ by $S(\m)=\bsigma(\m)-\phi(\m)$.
We let  $\sss=S(\uuu^{(0)})$, then with the notations of the preceding section,  we have
\be\label{eq:defSrond}
\sss=\{\sigma_i-\phi(\m_{i,j}), i=1,\ldots, N,\,j=1,\ldots,N_i\}.
\ee
In all the paper, we denote by $\underline \m=\m_{1,1}$ the (non necessarily unique) absolute minimum of $\phi$ that was chosen at the first step of the labelling procedure, and we let 
\be\label{eq:defulu0}
\ulu^{(0)}=\uuu^{(0)}\setminus\{\underline \m\}.
\ee
Using again the above labelling, we can associate a critical component to any local minimum. More precisely, we define 
 \begin{equation}\label{eq:definE}
 E:\uuu^{(0)}\rightarrow\Cr\cup\{X\}
 \end{equation}
 by $E(\m_{i,j})=E_{i,j}$.
 Observe that by definition, this application is injective.
Using this map, we can associate to each minimum $\m\in\uuu^{(0)}$  a boundary set given by  $\Gamma(\m)=\partial E(\m)$. Thanks to the fact that 
$\phi$ is a smooth Morse function, for any $\m\in \ulu^{(0)}$, the set $\Gamma(\m)$ is a finite union of compact sub-manifold of $X$ of dimension $d-1$ with conic singularities at the saddle points.
For our construction of quasimodes, we also need to introduce the set 
\begin{equation}\label{eq:defHm}
H(\m):=\{\m'\in E(\m)\cap\uuu^{(0)},\;\phi(\m')=\phi(\m)\}.
\end{equation}
Given $\m\in\ulu^{(0)}$, one has $\bsigma(\m)=\sigma_i$ for some $i\geq 2$. Moreover,  since $\sigma_{i-1}>\sigma_i$, there exists a unique connected component of 
$\{\phi<\sigma_{i-1}\}$ that contains $\m$ (observe that this component is not necessarily critical). We denote that component by  $E_-(\m)$, and by
 \begin{equation}\label{eq:definEmoins}
 E_-:\ulu^{(0)}\rightarrow\Omega(X)
 \end{equation}
the corresponding application, where $\Omega(X)$ is  the collection of connected open subsets of $X$.
Thanks to Remark \ref{rem:repart-min}, we know that
for any $\m\in\ulu^{(0)}$, there exists a
unique $\m'\in E_-(\m)\cap\uuu^{(0)}$, denoted by $\hat\m(\m)$, such that $\bsigma(\m')>\bsigma(\m)$.  In particular,
\begin{equation}\label{eq:comparminEmoins}
\forall \m\in\ulu^{(0)},\;\phi(\hat\m(\m))\leq \phi(\m),
\end{equation}
 and we denote  by $\widehat E(\m)$ the connected component of $\{\phi<\sigma(\m)\}$ containing $\hat\m(\m)$. It holds  additionally 
$\widehat E(\m)\subset E_-(\m)$ and we can  easily see that $\widehat E(\m)$ is always a critical component.
Throughout, we denote by
\begin{equation}\label{eq:defEhat}
\widehat E:\ulu^{(0)}\rightarrow\Cr
\end{equation}
and 
\begin{equation}\label{eq:defmhat}
\hat\m:\ulu^{(0)}\rightarrow \uuu^{(0)}
\end{equation}
 the corresponding applications.  
The fact that the  inequality in \eqref{eq:comparminEmoins} is large or strict plays an important role in our analysis.
\begin{defin}\label{def:type}
Let $\m\in\ulu^{(0)}$. We say that $\m$ is of type I if $\phi(\hat\m(\m))<\phi(\m)$. If $\phi(\hat\m(\m))=\phi(\m)$, we say that $\m$ is of type II.
We will denote
$$\uuu^{(0),I}=\{\m\in\ulu^{(0)},\,\m\text{ is of type I}\}$$
$$\uuu^{(0),II}=\{\m\in\ulu^{(0)},\,\m\text{ is of type II}\}$$
We have clearly the following disjoint union $\ulu^{(0)}=\uuu^{(0),I}\cup\uuu^{(0),II}$.
\end{defin}
\begin{example}\label{example1}
 Let us compute the preceding object in the case of the potential $\phi$ represented in Figure \ref{fig1}. The results are presented in Figure \ref{fig2}.
 
 $\bullet$ Let us start with the object associated to $\sigma_2$. By definition, $\widehat E(\m_{2,1})=\widehat E(\m_{2,2})=\widehat E(\m_{2,3})=\widetilde E_2$, where $\widetilde E_2$ is the connected component
 of $\{\phi<\sigma_2\}$ that contains $\m_{1,1}$. Then we have $\hat\m(\m_{2,1})=\hat\m(\m_{2,2})=\hat\m(\m_{2,3})=\m_{1,1}$.

 Since $\phi(\m_{1,1})=\phi(\m_{2,1})<\phi(\m_{2,3})<\phi(\m_{2,2})$, then $\m_{2,1}$ is of type II, whereas $\m_{2,2}$ and $\m_{2,3}$ are of type I.

$\bullet$ Consider now the level $\sigma_3$. One has $E_-(\m_{3,1})=E_-(\m_{3,2})=\widetilde E_2$ and 
$E_-(\m_{3,3})=E_-(\m_{3,4})=E_{2,3}$. Therefore, 
$\widehat E(\m_{3,1})=\widehat E(\m_{3,2})=\widetilde E_3$ where $\widetilde E_3$ is the connected component of $\{\phi<\sigma_3\}$ that contains $m_{1,1}$. Similarly, one has $\widehat E(\m_{3,3})=\widehat E(\m_{3,4})=\widetilde E'_3$ where $\widetilde E'_3$ is the connected component of $\{\phi<\sigma_3\}$ that contains $m_{2,3}$. From this computations, it follows that 
$\hat\m(\m_{3,1})=\hat\m(\m_{3,2})=\m_{1,1}$ and since $\phi(\m_{1,1})<\phi(\m_{3,1})=\phi(\m_{3,2})$ it follows that 
$\m_{3,1}$ and $\m_{3,2}$ are both of type I. On the other hand, $\hat\m(\m_{3,3})=\hat\m(\m_{3,4})=\m_{2,3}$ and since
$\phi(\m_{2,3})=\phi(\m_{3,3})<\phi(\m_{3,4})$ it follows that $\m_{3,3}$ is of type II and $\m_{3,4}$ of type I.
 
 $\bullet$ Finally, $E_-(\m_{4,1})=\widetilde E_3$, $\widehat E(\m_{4,1})=\widetilde E_4$ as represented on Figure \ref{fig2} and 
 $\hat\m(\m_{4,1})=\m_{1,1}$. Since $\phi(\m_{1,1})=\phi(\m_{4,1})$, it follows that $\m_{4,1}$ is of type II.

 \end{example}

The points of type II play an important role in our analysis. 
Given  $\sigma\in\Sigma$, let  $\Omega_\sigma=\Omega^0_\sigma\cup\widehat\Omega_\sigma$ with 
\begin{equation}\label{eq:defomegasigma}
\Omega^0_\sigma=\{E(\m),\m\in\bsigma^{-1}(\sigma)\}
\end{equation}
and $\widehat\Omega_\sigma$ defined by $\widehat\Omega_\sigma=\emptyset$ if $\sigma=\sigma_1$ and 
\begin{equation}\label{eq:defomegahatsigma}
\widehat\Omega_\sigma=\{\widehat E(\m),\,\m\in\bsigma^{-1}(\sigma)\cap\uuu^{(0),II}\} 
\end{equation}
if $\sigma\in \underline\Sigma$.
\begin{defin}\label{defin:releq}
We define an equivalence relation $\rrr$ on $\uuu^{(0)}$ by $\m\rrr\m'$ if and only if
\begin{equation}\label{eq:defreleqU0}
\left\{
\begin{array}{c}
\bsigma(\m)=\bsigma(\m')=\sigma\\
\exists\omega_1,\ldots\omega_K\in\Omega_\sigma\text{ s.t. } \m\in\omega_1,\m'\in\omega_K\text{ and }\forall k=1,\ldots K-1,\,\overline{\omega_k}\cap\overline{\omega_{k+1}}\neq\emptyset
\end{array}
\right.
\end{equation}
\end{defin}
Throughout the paper, we denote by $\Cl(\m)$ the equivalence class of $\m$ for the relation $\rrr$. 
Observe that since $\underline\m$ is the only minimum such that $\bsigma(\m)=\infty$, then 
$\Cl(\underline\m)=\{\underline\m\}$. 

 Let us denote by $(\uuu^{(0)}_\alpha)_{\alpha\in\aaa}$ the equivalence classes of $\rrr$
with $\aaa$ a finite set. We have evidently
\begin{equation}\label{decompU0-1}
\uuu^{(0)}=\bigsqcup_{\alpha\in\aaa}\uuu^{(0)}_\alpha.
\end{equation}
We need also to consider the set $\underline\aaa$ defined by $\underline\aaa=\aaa\setminus\{\underline\alpha\}$ where $\uuu^{(0)}_{\underline\alpha}=\{\underline \m\}$ is the equivalence class of the absolute minimum chosen for $\phi$. 
Throughout, we will denote $q_\alpha=\sharp \uuu^{(0)}_\alpha$. We will also  use the following partition of $ \uuu^{(0)}_\alpha$  for any $\alpha\in\ala$:
 \be\label{eq:defU0alphaI-II}
 \uuu^{(0),I}_\alpha:=\uuu^{(0)}_\alpha\cap\uuu^{(0),I},\;\;\uuu^{(0),II}_\alpha:=\uuu^{(0)}_\alpha\cap\uuu^{(0),II}.
 \ee

\begin{proposition}\label{prop:objetcteclass}
Let $\alpha\in\underline \aaa$. The applications $\bsigma, E_-,\widehat E$ and $\hat\m$ are constant on $\uuu^{(0)}_\alpha$.
\end{proposition}
\bp
For $\bsigma$, it is a direct consequence of the definition.
Suppose now that $\m,\m'\in\ulu^{(0)}$ satisfy $\m\rrr\m'$ and $\m\neq\m'$. 
Then, $\m$ and $\m'$ belong to the same connected component of $\{\phi\leq\sigma(\m)\}$. Hence, the uniqueness part in the definition of $E_-$ shows that $E_-(\m)=E_-(\m')$.
Since $E_-(\m)=E_-(\m')$, then the identity $\hat \m(\m)=\hat \m(\m')$ follows directly from the definition of $\hat\m$. This implies automatically
$\widehat E(\m)=\widehat E(\m')$.
\ep

 Thanks to the above proposition, given $\alpha\in\ala$,  we will denote respectively $\bsigma(\alpha),E_-(\alpha),\widehat E(\alpha)$ and $\hat\m(\alpha)$ instead of
 $\bsigma(\m),E_-(\m),\widehat E(\m)$, $ \hat\m(\m)$ for some $\m\in \uuu^{(0)}_{\alpha}$. 
 \begin{defin}\label{def:deftypeclasse}
  We  say that 
 \begin{itemize}
 \item[-] $\alpha$ is of type I, if $\phi(\hat\m(\alpha))<\phi(\m)$ for all $\m\in\uuu^{(0)}_\alpha$
 \item[-] $\alpha$ is of type II, if there exists  $\m\in\uuu^{(0)}_\alpha$ such that $\phi(\hat\m(\alpha))=\phi(\m)$.
 \end{itemize}
 \end{defin}

Recall that the height function $S:\uuu^{(0)}\rightarrow \R$ and the set of heights $\sss=S(\uuu^{(0)})$ were defined by \eqref{eq:defSrond} and above.
For any $\alpha\in \aaa$, we let 
\be\label{eq:defSalpha}
\sss_\alpha=S(\uuu^{(0)}_\alpha)\;\;\text{ and }\;\;p(\alpha)=\sharp \sss_\alpha
\ee
There exists some integers $\nu_1^\alpha<\nu_2^\alpha<\ldots<\nu_{p(\alpha)}^\alpha$ such that
$$
\sss_\alpha=\{S_{\nu_1^\alpha},\ldots,S_{\nu_{p(\alpha)}^\alpha}\}
$$

In the theorem below we sum up in a rather vague way the  description of these eigenvalues that we obtained  in Sections \ref{sec:COTASV} and \ref{sec:POMT}.

\begin{theorem}\label{th:main} 
There exist $c>0$ and  some symmetric positive definite matrices $\mmm^\alpha$, $\alpha\in\aaa$ such that counted with multiplicity, on has
$
\sigma(\Delta_\phi)=\bigcup_{\alpha\in\aaa}\sigma(\mmm^\alpha)(1+\ooo(e^{-c/h}))
$
with
$$
\sigma(\mmm^\alpha)=\bigcup_{j=1}^{p(\alpha)}he^{-2h^{-1}S_{\nu_j^\alpha}}\sigma(M^{\alpha,j})
$$
for some  symmetric positive definite matrices  $M^{\alpha,j}$ having a classical expansion with invertible leading term given in Theorem \ref{thm:SVLalpha}.
\end{theorem}

Let us make a few comments on this theorem.

First, observe that  since $M^{\alpha,j}$ has a classical expansion with invertible leading term $M^{\alpha,j}_0$, then its eigenvalues 
$\zeta^{\alpha,j}_r$, $r=1,\ldots,r^{\alpha,j}$ have a classical expansion 
$$
\zeta^{\alpha,j}_r(h)\sim\sum_k h^k\zeta^{\alpha,j}_{r,k}
$$
with $\zeta^{\alpha,j}_{r,0}$ eigenvalue of the matrix $M^{\alpha,j}_0$. 

Compared to previous results obtained under the Generic Assumption, the main difference is that the prefactor $\zeta_{r,k}^{\alpha,j}$ are more difficult to compute since they are obtained as the eigenvalues of the matrices 
$M^{\alpha,j}$. When (GA) is satisfied, the $M^{\alpha,j}$ are $1\times 1$ matrices whose spectrum is direct to obtain. In the general case, this is not true anymore and the construction of the matrices $M^{\alpha,j}$ is more involved.
In particular, it  depends dramatically on the number $p(\alpha)=\sharp S(\uuu^{(0)}_\alpha)$. Observe that this number is also equal to the number of different values taken by $\phi$ on the equivalence class $\uuu^{(0)}_\alpha$.

If $p(\alpha)=1$, the coefficients of $M^{\alpha,j}$ depend only on the couples $(\m,\s)$ for which 
$\phi(\s)-\phi(\m)=S_{\nu_j^\alpha}$. Excepted the fact that the different eigenvalues $\zeta^{\alpha,j}_r, r=1,\ldots,r^{\alpha,j}$ are linked together, the situation is similar to that encountered in the generic case. Actually, we prove in appendix that if (GA) is satisfied then $\Cl(\m)$ is reduced to one point for any $\m$, and in particular $p(\alpha)=1$ for all $\alpha$.

In the case where $p(\alpha)\geq 2$, the matrix is more difficult to compute. It comes from an application of  
Schur complement's method and it depends on some couples $(\m,\s)$ for which the height $\phi(\s)-\phi(\m)$ is smaller than $S_{\nu^\alpha_j}$.
In other words, the lifetime of the metastable state $\m$ is not entirely described by the height that is needed to jump over in order to reach
the nearest lower energy position. It depends also on some interactions with some higher energy states that are not present in the classical 
Eyring-Kramers formula.
To our knowledge, this is is the first time that such a phenomena is exhibited.\\

\begin{figure}
 \center
  \scalebox{0.5}{ 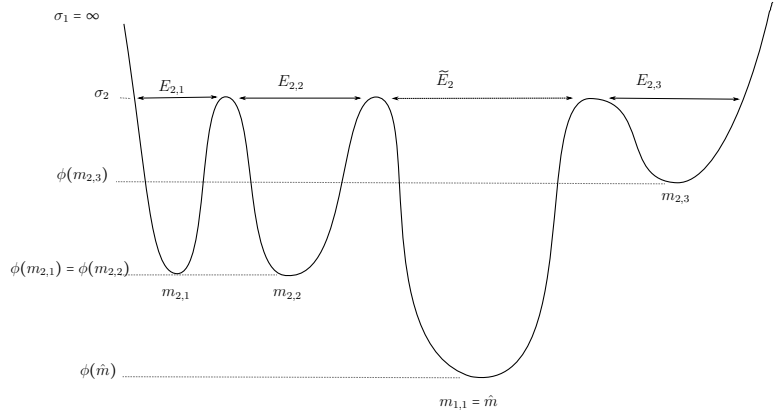}
  \caption{A potential with $p(\alpha)=1$ }
    \label{fig3}
\end{figure}

\begin{figure}
 \center
  \scalebox{0.5}{ 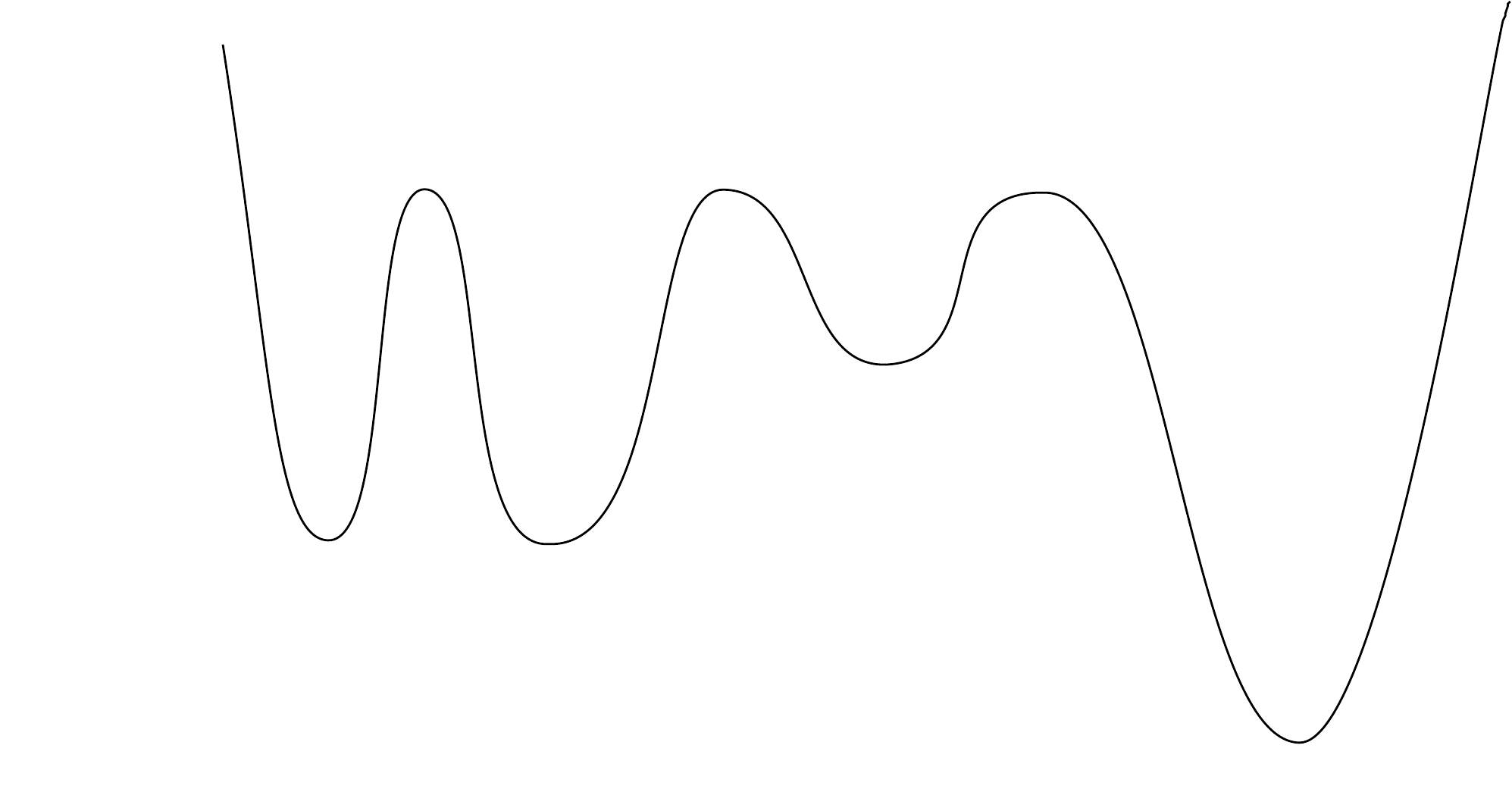}
  \caption{A potential with $p(\alpha)=2$ }
    \label{fig4}
\end{figure}

Let us now compute $p(\alpha)$ on explicit examples. Let us fix $n=2$ and consider the potentials $\phi$  
given respectively by Figures \ref{fig3}, \ref{fig4}.
In both cases $\hat\m(\m_{2,1})=\hat\m(\m_{2,2})=\hat\m(\m_{2,3})=\m_{1,1}$ that we denote by $\hat\m$ for short.
Since $\phi(\hat\m)<\phi(\m_{2,j})$ for all $j$, then there is no point of type II, $\uuu^{(0),II}=\emptyset$ and hence
$\Omega_{\sigma_2}=\{E_{2,1},E_{2,2}, E_{2,3}\}$.
Therefore, one can compute easily the equivalence classes of $\rrr$ in both cases:
\begin{itemize}
\item[-] in the case of Figure \ref{fig3}, we have 3 equivalence classes: $c_1=\{\m_{1,1}\}$, $c_2=\{\m_{2,1},\m_{2,2}\}$ and $c_3=\{\m_{2,3}\}$.
The potential $\phi$ is constant on each equivalence class, and hence $p(c_1)=p(c_2)=p(c_3)=1$.
\item[-] in the case of Figure \ref{fig4}, we have 2 equivalence classes: $c_1=\{\m_{1,1}\}$, $c_2=\{\m_{2,1},\m_{2,2},\\ \m_{2,3}\}$.
The potential $\phi$ takes two different values on $c_2$: $p(c_2)=2$.
\end{itemize}
We will come back to these examples at the end of the paper and compute explicitly the spectrum of $\Delta_\phi$ in both cases.


\begin{figure}
 \center
  \scalebox{0.5}{ 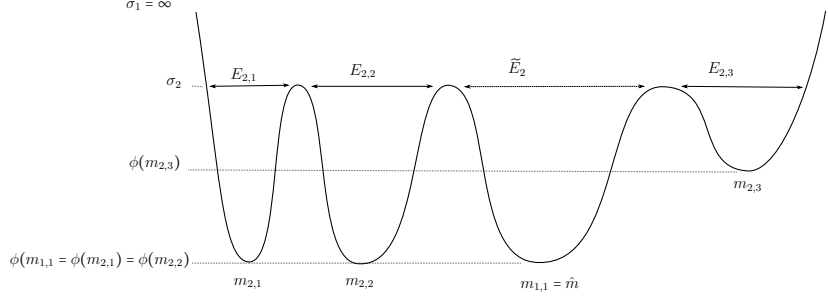}
  \caption{An example with points of type II }
    \label{fig5}
\end{figure}

Let us finish this discussion by an example where $\uuu^{(0),II}\neq\emptyset$. Consider the potential given by Figure \ref{fig5}.
In that case $\hat\m(\m_{2,1})=\hat\m(\m_{2,2})=\hat\m(\m_{2,2})=\m_{1,1}$ that we denote by $\hat\m$ for short.
Since $\phi(\hat\m)=\phi(\m_{2,1})=\phi(\m_{2,2})<\phi(\m_{2,3})$, then $\m_{2,1}$ and $\m_{2,2}$ are of type II and $\m_{2,3}$ is of type I.
We still have $\Omega^0_{\sigma_2}=\{E_{2,1},E_{2,2}, E_{2,3}\}$ but contrary to the previous case 
 $\widehat\Omega_{\sigma_2}=\{\widetilde E_2\}$ is non empty. It follows that 
 $\Omega_{\sigma_2}=\{E_{2,1},E_{2,2}, E_{2,3},\widetilde E_2\}$ and  $\rrr$ admits two equivalence classes:
  $c_1=\{\m_{1,1}\}$, $c_2=\{\m_{2,1},\m_{2,2},\m_{2,3}\}$. 
  The potential $\phi$ takes two different values on $c_2$ and hence $p(c_2)=2$.

%

%

\subsection{General strategy of the proof}\label{subsec:GSOP}

Let us recall the general strategy followed in \cite{HeKlNi04_01}. The starting point is to use the supersymmetric structure of the Witten Laplacian.
For $0\leq k\leq n$, let $\Omega^k(X)=\ccc^\infty(X,\Lambda^kT^*X)$ be the space of $k$-differential forms and denote $d:\Omega^k(X)\rightarrow \Omega^{k+1}(X)$  the exterior derivative and
$d^*: \Omega^k(X)\rightarrow \Omega^{k-1}(X)$  its adjoint for the natural pairing.
The  Witten complex associated to the function $\phi$ is defined by the semiclassical weighted de Rham differentiation
$$\dh=e^{-\phi/h}\circ h d\circ e^{\phi/h}=hd+d\phi^\wedge$$
and its adjoint  
$$\dhs=e^{\phi/h}\circ h d^{*}\circ e^{-\phi/h}=hd^*+d\phi^\righthalfcup.$$
Then the semiclassical Witten Laplacian is defined on the forms of any degree by
\begin{equation}\label{eq:defWittenSusy}
\Delta_\phi=\dhs\circ\dh+\dh\circ\dhs.
\end{equation}
When restricted to the space of $p$-forms we denote this operator by $\Delta^{(p)}_\phi$ (observe that in the case $p=0$, the above formula yields easily \eqref{eq:definWitten0form}). 
Then, we have the following intertwining relation
\begin{equation}\label{eq:intertwin}
\dh \Delta_\phi^{(p)}=\Delta_\phi^{(p+1)}\dh
\end{equation}
and its analogous for the coderivative
\begin{equation}\label{eq:intertwinstar}
\dhs \Delta_\phi^{(p+1)}=\Delta_\phi^{(p)}\dhs.
\end{equation}

For any $p=0,\ldots,d$, it follows from  \eqref{eq:hypgenephi} that $\Delta_\phi^{(p)}$ (as an unbounded operator on $L^2$) is essentially self-adjoint on the space of compactly supported smooth forms. We still denote by 
$\Delta_\phi^{(p)}$ its unique self-adjoint extension.
Then $\Delta_\phi^{(p)}$ is non-negative and thanks to  \eqref{eq:hypgenephi}, there exists $c_0>0$ such that 
$\sigma_{ess}( \Delta_\phi^{(p)})\subset[c_0,+\infty[$  for any 
$h>0$ small enough (in the case where $X$ is a compact manifold,  $\Delta_\phi^{(p)}$ has actually compact resolvent).
Moreover, there exists $\epsilon_p>0$ such that for $h>0$ small enough, it has exactly
$n_p$ eigenvalues in the interval $[0,\epsilon_p h]$ where $n_p$ denotes the number of critical points of index $p$ of $\phi$. We shall denote by $E^{(p)}$ the spectral subspace 
associated to these small eigenvalues of $\Delta_\phi^{(p)}$. Then $\dim E^{(p)}=n_p$ and relations \eqref{eq:intertwin}, \eqref{eq:intertwinstar} show that 
\begin{equation}\label{eq:actionderivEp}
\dh(E^{(p)})\subset E^{(p+1)}\text{ and }\dhs(E^{(p+1)})\subset E^{(p)}.
\end{equation}
This shows in particular that $\dh$ acts from $E^{(0)}$ into $E^{(1)}$ and we shall denote by $\lll$ this operator.
Similarly $\Delta_\phi^{(0)}$ acts on $E^{(0)}$ and we denote by $\mmm$ this operator. By \eqref{eq:defWittenSusy}, we get 
$$
\mmm=\lll^*\lll.
$$ 
The general strategy used in \cite{HeKlNi04_01} (that we will follow in the present work), is to construct  appropriate bases
of $E^{(0)}$ and $E^{(1)}$ in which one can compute handily the singular values of $\lll$. The main idea to construct such bases is to build accurate quasimodes for $\Delta_\phi$ and to project them on the spaces $E^{(j)}$. The construction of the quasimodes is performed in Section \ref{sec:COAQ}. The quasimodes for $1$-forms are the ones constructed in \cite{HeSj85_01}. The main properties of these quasimodes will be recalled in Section \ref{subsec:QMF1F}. Concerning the quasimodes on $0$-forms, one can not use  the ones constructed in \cite{HeKlNi04_01} since many important properties that are required for our analysis fail to be true in the present situation (for instance, the quasi-orthogonality). In  Section \ref{subsec:QMF0F}, we use the partition of $\uuu^{(0)}$ into equivalence classes of $\rrr$ to construct a family of quasimodes on $0$-forms adapted to our setting. 
Each quasimode will be associated to a minimum $\m\in\uuu^{(0)}$.

In Section \ref{sec:PFSVA}, we compute  the matrix $\lll$ in the above basis. One  arrives to a block diagonal matrix 
$\diag(\lll^\alpha,\alpha\in\aaa)$ whose singular values are the singular values of each block.

Section  \ref{sec:COTASV} is devoted to the computation of singular values of the above blocks. The main difficulty is that  given two minima $\m,\m'$ in the same equivalence class, one has not necessarily $S(\m)=S(\m')$.
 For equivalence classes satisfying this property (that is $p(\alpha)=1$), 
each block $\lll^\alpha$ of the matrix $\lll$ has a typical size $e^{-S(\alpha)/h}$ and the situation could be handled quite easily. 
But more complicated cases may arise where quasimodes yielding different heights $S(\m)$ are interacting. In order to treat the full general case, we 
use Schur complement method combined with an induction on $p(\alpha)$.  Running the induction step requires to exhibit a specific structure of the matrices under consideration (see Sections \ref{subsec:IOTML} and  \ref{subsec:GSOTM}). In Section \ref{subsec:TSOGM}, we prove a general result for such matrices that we use to conclude in Section \ref{subsec:TSVOL}.

In Section \ref{sec:POMT}, we prove Theorem \ref{th:main}.

In a separate appendix, we collect several result in linear algebra. We also provide a list of notations used in the paper.\\

{\bf Acknowledgement:}  The author would like to thank D. Le Peutrec for numerous discussions; in particular for pointing out a fondamental argument in the induction step of the proof of the main result.  This paper was written while the author was in leave at Stanford University and U.C. Berkeley.
He was supported by the European Research Council, 
ERC-2012-ADG, project number 320845:  Semi Classical Analysis of Partial Differential
Equations, by the Stanford Mathematic Department and by the France Berkeley Fund, FBF 2016-0071, Semiclassical study of randomness and dynamics.

\section{ Construction of adapted quasimodes}\label{sec:COAQ}

\subsection{Gathering minima by equivalence class}\label{subsec:GMBEC}
Let us start this section with  a proposition 
 collecting some elementary facts about $E,E_-$ and $\widehat E$.
 \begin{proposition}\label{prop:propelemE}
 Let $\m,\m'\in\uuu^{(0)}$ such that $\m\neq \m'$. Then, we have the following
 \begin{enumerate}
 \item[i)] If $\sigma(\m)=\sigma(\m')$, then:
 \begin{itemize}
 \item[i.a)] $E(\m)\cap E(\m')=\emptyset$
 \item[i.b)] either $E_-(\m)=E_-(\m')$ or $E_-(\m)\cap E_-(\m')=\emptyset$ 
 \item[i.c)] if  $E_-(\m)=E_-(\m')$  then $\widehat E (\m)=\widehat E(\m')$ otherwise $\widehat E(\m)\cap \widehat E(\m')=\emptyset$ 
 \end{itemize}
 \item[ii)] If $\sigma(\m)>\sigma(\m')$, then
   \begin{itemize}
 \item[ii.a)]either $E(\m)\cap E(\m')=\emptyset$ or $E_-(\m')\subset E(\m)$
 \item[ii.b)]either $E_-(\m)\cap E_-(\m')=\emptyset$ or $E_-(\m')\subset E_-(\m)$
\end{itemize}
 \end{enumerate}
 \end{proposition}
 \bp
Let $\m\ne\m'$ be two minima. Assume first that $\sigma(\m)=\sigma(\m')=\sigma$. 
Since $\m\neq\m'$ and $\bsigma^{-1}(\infty)=\{\underline\m\}$, then one has necessarily $\m,\m'\in\ulu^{(0)}$. In particular, 
$E_-(\nu),\widehat E(\nu)$, $\nu=\m,\m'$ are well-defined.
Moreover,  $E(\m)$ and $E(\m')$ being  two connected components of 
$\{\phi<\sigma\}$, one has either $E(\m)=E(\m')$ or $E(\m)\cap E(\m')=\emptyset$. Since $\m\neq \m'$ and $E$ is injective, then $E(\m)\cap E(\m')=\emptyset$, which proves i.a). 

Since $E_-(\m)$ and $E_-(\m')$ are two connected component of the same set $\{\phi<\tau\}$ for some $\tau>\sigma(\m)$, then i.b) is obvious.

Suppose now that $E_-(\m)=E_-(\m')$. Since $\sigma(\m)=\sigma(\m')$ then $\hat\m(\m)=\hat\m(\m')$. Moreover, $\widehat E(\m)$ being the unique connected component of $\{\phi<\sigma(\m)\}$ containing $\hat\m(\m)$, we get  $\widehat E(\m)=\widehat E(\m')$. If $E_-(\m)$ and $E_-(\m')$ are disjoint, then $\widehat E(\m)$ and $\widehat E(\m')$ are also disjoint since $\widehat E(\m)\subset E_-(\m)$ and $\widehat E(\m')\subset E_-(\m')$. This completes the proof of i.c).

Let us now prove ii) and assume that $\sigma(\m)>\sigma(\m')$. Once again, since $\bsigma^{-1}(\infty)=\{\underline\m\}$, then
$\m'\in\ulu^{(0)}$. If $E(\m')\cap E(\m)\neq \emptyset$, then $E_-(\m')\cap E(\m)\neq \emptyset$.
Moreover, $E_-(\m')$ is a connected component of $\{\phi<\tau\}$ for some $\tau\leq \sigma(\m)$. Since $E(\m)$ is a connected component of 
$\{\phi<\sigma(\m)\}\supset \{\phi<\tau\}$, then $E_-(\m')\subset E(\m)$ which proves ii.a).

The point  ii.b) is proved by similar arguments.
 \ep

Let us now  decompose the set of separating saddle points according to the equivalence classes. Given $\alpha\in\ala$, introduce the closed set 
\begin{equation}\label{eq:definF}
G(\alpha)=\bigcup_{\m\in\uuu^{(0)}_\alpha}\overline{E(\m)}
\end{equation}
and for any $\alpha\in\ala$ let
\begin{equation}\label{eq:definV1alpha}
\vvv^{(1)}_{\alpha}=\{\s\in \vvv^{(1)},\,\phi(\s)=\bsigma(\alpha)\}\cap G(\alpha).
\end{equation}
For any $\alpha\in\ala$, let
\be\label{eq:defparenU}
\puu^{(0)}_\alpha=\uuu^{(0)}_\alpha\cup\{\hat\m(\alpha)\}
\ee
and define an application $\Gamma_\alpha$ from $\puu^{(0)}_\alpha$ into the 
closed subsets of $X$ by
\be\label{eq:defGamalpha}
\left\{
\begin{array}{c}
\Gamma_\alpha(\m)=
\Gamma(\m) \text{ if }\m\in\uuu^{(0)}_\alpha\\
\Gamma_\alpha(\hat\m(\alpha))=\partial \widehat E(\alpha).\phantom{****}
\end{array}\right.
\ee
where $\Gamma$ is defined below \eqref{eq:definE}.
\begin{remark}\label{rem:defm1m2} 
Since $\widehat E(\m)\subsetneq E(\hat\m)$, then the application $\Gamma_\alpha$ is slightly different from the application $\Gamma$ defined in below \eqref{eq:definE}.
Observe also that for all $\m\in \puu^{(0)}_\alpha$, $\Gamma_\alpha(\m)$ is the boundary of the connected component of 
$\{\phi<\phi(\s)\}$ that contains $\m$. 

\end{remark}

\begin{lemma}\label{lem:defm1m2}
The collection $(\vvv^{(1)}_\alpha)_{\alpha\in\ala}$ is a partition of $\vvv^{(1)}$.
Moreover, for all $\alpha\in\ala$ and $\s\in\vvv^{(1)}_\alpha$, there exists $\m_1(\s)\in \uuu^{(0)}_\alpha$ and 
$\m_2(\s)\in \puu^{(0)}_\alpha$ such that 
\be\label{eq:defm1m2}
\s\in\Gamma_\alpha(\m_1)\cap\Gamma_\alpha(\m_2)
\ee
One can chose $\m_1,\m_2$ in order that 
$S(\m_1)\leq S(\m_2)$ (that is $\phi(\m_1)\geq \phi(\m_2)$). Up to permutation, the couple $(\m_1(\s),\m_2(\s))$ is unique.
\end{lemma}

\bp
Let $\s\in\vvv^{(1)}$, then $\phi(\s)\in\underline\Sigma$ and there exists $k\geq 2$ such that $\phi(\s)=\sigma_k$. By definition,  there exists two different connected components $E_1,E_2$ of $\{\phi<\sigma_k\}$ such that $\s\in\overline E_1\cap\overline E_2$. From the existence part of Remark \ref{rem:repart-min} there exist $\m_{l,i}\in E_1$ and $\m_{l',i'}\in E_2$ with $l'\leq l\leq k$. Moreover, one has necessarily $l=k$. 
Otherwise $\bsigma(\m_{l,i})>\sigma_k$ and since $\overline E_1\cap\overline E_2\neq\emptyset$, this would imply that  $\m_{l',i'}\in E(\m_{l,i})$ which is impossible since $l'\leq l$. 
Hence we have $l=k$. Therefore   $E_1$ is equal to  $E(\m_{l,i})$ with  $\m_{l,i}\in\uuu^{(0)}_\alpha$, which proves that $\s\in\vvv^{(1)}_\alpha$. Moreover, $E_2$  is either of the form $E_2=E(\m_{l',i'})$  with $\m_{l',i'}\in\uuu^{(0)}_\alpha$ (if $l'=k$) or $E_2=\widehat E(\m_{l,i})$ (if $l'<k$). 
Setting $\m_1(\s)=\m_{l,i}\in \uuu^{(0)}_\alpha$ and $\m_2(\s)=\m_{l',i'}\in\puu^{(0)}_\alpha$, one has
$\s\in\Gamma_\alpha(\m_1)\cap\Gamma_\alpha(\m_2)$ and  since $l\geq l'$ one has also $\phi(\m_1)\geq\phi(\m_2)$.

Let us now prove that the union of the $\vvv^{(1)}_\alpha$ for  $\alpha\in\ala$ is disjoint.
Suppose that $\s\in\vvv^{(1)}_\alpha\cap \vvv^{(1)}_\beta$. Then $\sigma(\alpha)=\phi(\s)=\sigma(\beta)$. Moreover, there exists 
$\m\in\uuu^{(0)}_\alpha$ and $\m'\in \uuu^{(0)}_\beta$ such that $\s\in \overline{E(\m)}\cap \overline{E(\m')}$. This proves that $\m\rrr\m'$
and hence $\alpha=\beta$.

The uniqueness of $(\m_1,\m_2)$ up to permutation is obvious.
\ep


Let us now introduce an extra partition that will be useful in the sequel.
\begin{lemma}\label{lem:partitionValpha1} For all
$\alpha\in \aaa$ there exists  a partition
$\vvv^{(1)}_\alpha=\vvv^{(1),\b}_\alpha\sqcup\vvv^{(1),\i}_\alpha$
such that the following hold true:
\begin{enumerate}
\item[i)] for any  $\s\in \vvv^{(1),\i}_{\alpha}$, $\m_1(\s)$ and $\m_2(\s)$ belong to $\uuu^{(0)}_{\alpha}$.
\item[ii)] the set $\vvv^{(1),b}_\alpha$ is non-empty and 
for all $\s\in \vvv^{(1),b}_\alpha$ one has $\m_1(\s)\in \uuu^{(0)}_{\alpha}$, $\m_2(\s)=\hat\m(\alpha)$ and 
$$\s\in\Gamma_\alpha(\m_1(\s))\cap\Gamma_\alpha(\hat\m(\alpha)).$$
\end{enumerate}
\end{lemma}
\bp
Define 
$\vvv^{(1),\i}_{\alpha}=\{\s\in \vvv^{(1)}_\alpha,\,\m_1(\s),\m_2(\s)\in\uuu^{(0)}_\alpha\}$.
 Then $i)$ is true by definition. Moreover, defining
 $\vvv^{(1),\b}_{\alpha}=\vvv^{(1)}_{\alpha}\setminus\vvv^{(1),\i}_{\alpha}$, one has automatically the partition property and it remains to prove ii). 
 
Since $\alpha\in\ala$, the set  
 $
 \overline {\widehat E(\alpha)}\cap\left(\cup_{\m\in \uuu^{(0)}_\alpha}\overline {E(\m)}\right)
 $
 is non-empty and contained in $\vvv^{(1),\b}_{\alpha}$. This proves  that   $\vvv^{(1),\b}_{\alpha}$ is not empty.
Suppose now that $\s\in\vvv^{(1),b}_\alpha$. It follows from 
Lemma \ref{lem:defm1m2} that  $\m_1(\s)\in\uuu^{(0)}_\alpha$ and $\m_2(\s)\in \widehat \uuu^{(0)}_\alpha$. 
But by definition of $\vvv^{(1),\b}_\alpha$, $\m_2(\s)$ can not belong to 
$\uuu^{(0)}_\alpha$, which implies by definition that  $\m_2(\s)=\hat\m(\alpha)$. This completes the proof of ii).
\ep

\subsection{Quasimodes for $0$-forms}\label{subsec:QMF0F}
In this section we construct a family of quasimodes of $\Delta_\phi^{(0)}$ associated to the minima of $\phi$. Each of these quasimodes will be of the form 
$x\mapsto h^{-\frac d 4}\chi_\m(x)e^{-(\phi(x)-\phi(\m))/h}$ with some suitable cut-off functions $\chi_\m$ associated to a minimum $\m\in\uuu^{(0)}$.

Following \cite{HeKlNi04_01}, we can associate to each minimum $\m\in\uuu^{(0)}$ a cut-off function $\chi_\m$  in the following way. 
For $\m=\underline\m$, we simply take $\chi_{\underline\m}=1$. For $\m\in\ulu^{(0)}$ we introduce some small parameters 
$\epsilon,\tilde\epsilon, \delta>0$ with $\tilde\epsilon<\epsilon$ and 
we define 
\begin{equation}\label{eq:defin-bassinlarge}
E_{\epsilon,\tilde\epsilon,\delta}(\m)=\left(\big(E(\m)\setminus\bigcup_{\s\in \vvv^{(1)}\cap\Gamma(\m)}B(\s,\epsilon)\big)+B(0,\tilde \epsilon)\right)
\bigcup\left(\bigcup_{\s\in(\uuu^{(1)}\setminus\vvv^{(1)})\cap\Gamma(\m)}B(\s,\delta)\right)
\end{equation}

\begin{proposition}\label{prop:propcutoff}
Let $\chi_\m$ be any function in $\CC^{\infty}_c(E_{\epsilon,2\tilde\epsilon,2\delta}(\m))$ such that $\chi_\m=1$ on $E_{\epsilon,\tilde\epsilon,\delta}(\m)$.
There exists $\epsilon_0>0$, $\delta_0>0$ and $C>0$ such that for all $0<\delta<\delta_0$, all  $0<\epsilon<\epsilon_0$ and all $0<\tilde\epsilon<\epsilon/4$, the following hold 
true:
\begin{enumerate}
\item[a)] if $x\in \supp(\chi_\m)$ and $\phi(x)<\sigma(\m)$, then $x\in E(\m)$
\item[b)] there exists $c_\epsilon>0$ such that for all 
$x\in \supp(\nabla\chi_\m)$, we have
\begin{itemize}
\item[-] either $x\notin \cup_{\s\in \vvv^{(1)}\cap\Gamma(\m)}B(\s,\epsilon)$ and 
$$\sigma(\m)+c_\epsilon^{-1}<\phi(x)<\sigma(\m)+c_\epsilon$$
\item[-] or $x\in \cup_{\s\in \vvv^{(1)}\cap\Gamma(\m)}B(\s,\epsilon)$ and 
$$\vert \phi(x)-\sigma(\m)\vert\leq C\epsilon.$$
\end{itemize}
\item[c)] for all $\s\in\uuu^{(1)}\setminus(\vvv^{(1)}\cap\Gamma(\m))$, one has 
$dist(\s,\supp\nabla\chi_\m)\geq \delta$. 
If moreover $\s\in\supp(\chi_\m)$ then $\s\in E(\m)$ (in particular $\chi_\m(\s)=1$).
\item[d)] suppose that $\m\in\uuu^{(0)}_{\alpha},\,\alpha\in\aaa$ and let $\s\in\vvv^{(1)}\cap\supp(\chi_\m)$. Then, there exists $\beta\in\ala$  such that   $\bsigma(\beta)<\bsigma(\alpha)$, 
$\s\in\vvv^{(1)}_{\beta}$ and $\cup_{\m'\in \uuu^{(0)}_{\beta}}E(\m')\subset \{x\in X,\,\chi_{\m}(x)=1\}$.
\end{enumerate}

\end{proposition}
\bp Observe that the construction of the cut-off functions $\chi_\m$ is slightly different to that of the $\chi_{k,\epsilon}$ in Proposition 4.2 in \cite{HeKlNi04_01} (in particular because there can exist more than one separating saddle point  on $\partial E(\m)$). 

Let  $\delta_1=\min\{ \vert \s-\s'\vert,\;\s,\s'\in\uuu^{(1)},\;\s\neq\s'\}$ and $\delta_2=\min\{\dist(\s,\Gamma(\m)),\; \s\in E(\m)\cap\uuu^{(1)}\}$.
Let $0<\delta<\frac 14 \min(\delta_1,\delta_2)$ and 
$\epsilon_0>0$ such that there exists $C>0$ such that for all $0<\epsilon<\epsilon_0$ and all $\s\in\vvv^{(1)}$, one has 
$$
\vert \phi(x)-\phi(\s)\vert < C\epsilon,\;\forall x\in B(\s,\epsilon).
$$
This is possible since $\phi$ is a smooth function.
Then a) and b) above can be proved in a similar way as in Proposition 4.2 in \cite{HeKlNi04_01} and c) is a direct consequence of our choice of 
$\delta$.

Let us now prove d). By definition, if $\s\in\vvv^{(1)}\cap\supp(\chi_\m)$, then $\s\in E(\m)$ (here we use the condition $0<\tilde\epsilon<\epsilon/4$). Hence, there exists $\beta\neq \alpha$ such that 
$\s\in \vvv^{(1)}_\beta$ and one has additionally $\bsigma(\beta)<\bsigma(\alpha)$. By definition of the sets $E(\m)$, this implies that 
$\cup_{\m'\in \uuu^{(0)}_{\beta}}E(\m')\subset E(\m)\setminus\cup_{\s'\in \vvv^{(1)}\cap\Gamma(\m)}B(\s',\epsilon)$ for any $\epsilon\in]0,\epsilon_0[$ 
with $\epsilon_0>0$ small enough independent of $\delta$. This implies the results.
\ep

We are now in position to define the quasimodes in a recursive way on the values of $\bsigma(\alpha)$. 
\begin{itemize}
\item[-] We start with the quasimode associated to $\underline \m$. We set 
\begin{equation}\label{definGS}
f_{\underline \m}^{(0)}(x)=c(\underline \m,h)h^{-d/4}e^{(\phi(\underline\m)-\phi(x))/h}
\end{equation}
where $c(\underline\m,h)$ is a normalizing constant such that $\Vert f_{\underline \m}\Vert_{L^2}=1$.
Due to the fact that $\phi$ may have several global minima, the function $f_{\underline \m}^{(0)}$ does not concentrate 
only on $\underline \m$ but on the reunion of all global minima. Hence the normalizing factor $c(\underline \m,h)$ is computed by adding the contributions coming from each of these minima via quadratic approximation. 
More precisely, it follows from the Laplace method that 
$c(\underline \m,h)\sim\sum_{k=0}^\infty h^k\gamma_k(\underline \m)$ with the function $\gamma_0$ given by 
\begin{equation}\label{eq:defnormcte}
\gamma_0(\m)^{-2}=\pi^{\frac d 2} \sum_{\m'\in H(\m)}\vert \det\Hess\phi(\m')\vert^{-\frac 12}
\end{equation}
where by  definition \eqref{eq:defHm} one has
\begin{equation*}
H(\m):=\{\m'\in E(\m)\cap\uuu^{(0)},\;\phi(\m')=\phi(\m)\}.
\end{equation*}

Eventually, observe that $f_{\underline\m}^{(0)}$ is an exact quasimode: $\Delta_\phi f_{\underline \m}^{(0)}=0$.

\item[-] Suppose now that $k\in\{2,\ldots,K\}$ and that the quasimodes $f_\m^{(0)}$ have been constructed for $\m\in \bigcup_{\alpha'\in\ala,\bsigma(\alpha')\leq\sigma_{k-1}}\uuu^{(0)}_\alpha$, and 
let us define $f_\m^{(0)}$ for $\m\in \uuu^{(0)}_{\alpha}$ with $\bsigma(\alpha)=\sigma_k$. 
The  form of the quasimode associated to $\m$ depends on the  type of $\m$ as introduced in Definition \ref{def:type}.
\begin{itemize}
\item[$\star$] If $\m$ is of type I, then we define the $f_\m^{(0)}$ as in \cite{HeKlNi04_01} by 
\begin{equation}\label{eq:definQMusual}
f_{\m}^{(0)}(x)=c(\m,h)h^{-d/4}\chi_\m(x)e^{(\phi(\m)-\phi(x))/h}
\end{equation}
where $\chi_\m$ is the cut-off function associated to $\m$ defined in Proposition \ref{prop:propcutoff} and $c(\m,h)$ is again a normalizing constant such that $\Vert f_{\m}^{(0)}\Vert_{L^2}=1$. As before, we have to add all the contributions of minima in $E(\m)$ at the same height as $\m$. We get  
$c(\m,h)\sim\sum_{k=0}^\infty h^k\gamma_k(\m)$ with $\gamma_0(\m)$ given by \eqref{eq:defnormcte}.

\item[$\star$] Let us now construct quasimodes associated to minima $\m$ of type II. We assume here that $\uuu^{(0),II}\neq\emptyset$ and we define
\be\label{eq:defhatU0II}
\widehat \uuu^{(0),II}_{\alpha}=\uuu^{(0),II}_{\alpha}\cup\{\hat \m\}.
\ee
where for short, 
we denote $\hat\m=\hat\m(\alpha)$ and $q_\alpha^{II}=\sharp \uuu^{(0),II}_{\alpha}$.

Let us introduce an additional  cut-off function around $\hat \m$ that we define as follows.
Recall that $\widehat E( \alpha)$ denotes the connected component of  $\{x\in E_-(\m),\;\phi(x)<\sigma(\m)\}$ that contains $\hat \m$. As before, we introduce some parameters 
$\epsilon,\tilde\epsilon,\delta>0$ with $\tilde\epsilon<\epsilon$ and we define
$$
\widehat E_{\epsilon,\tilde\epsilon,\delta}(\alpha)=\left(\big(\widehat E( \alpha)\setminus\bigcup_{\s\in \vvv^{(1)}\cap\partial \widehat E(\alpha)}B(\s,\epsilon)\big)+B(0,\tilde \epsilon)\right)
\bigcup\left(\bigcup_{\s\in(\uuu^{(1)}\setminus\vvv^{(1)})\cap\partial \widehat E(\alpha)}B(\s,\delta)\right)
$$
Then, we let $\hat\chi_{\hat \m}$ be any function in $\CC^{\infty}_c(\widehat E_{\epsilon,2\tilde\epsilon, 2\delta}(\alpha))$ such that $\hat\chi_{\hat \m}=1$ on $\widehat E_{\epsilon,\tilde\epsilon,\delta}(\alpha)$.
For $\m\in \uuu^{(0),II}_{\alpha}$, we let $\hat\chi_{\m}=\chi_{\m}$, with $\chi_{\m}$ defined in Proposition \ref{prop:propcutoff}. 
We want to construct the quasimode as linear combination of the $\hat\chi_{\m} e^{-\phi/h}$, $\m\in\widehat\uuu^{(0),II}_\alpha$. In order to chose the coefficients, let us introduce
 $\fff_{\alpha}=\Fr(\widehat\uuu^{(0),II}_\alpha)$ the finite vector space of functions from $\widehat \uuu^{(0),II}_{\alpha}$ into $\R$.
This space has dimension $q_\alpha^{II}+1$ and is  endowed with the usual euclidean structure
$$
\<\theta,\theta'\>_{\fff_\alpha}=\sum_{\m\in \widehat \uuu^{(0),II}_{\alpha}}\theta(\m)\theta'(\m).
$$
 We denote by $N$ the associated norm.
Eventually, we  define  
$\theta^\alpha_0\in\fff_{\alpha}$
by 
\begin{equation}\label{eq:deftheta0}
\theta_0^\alpha(\m)=\frac {c^\alpha_0(h)}{c(\m,h)}
\end{equation}
where $ c(\m,h)$ is the unique positive constant such that 
the function 
$$\tilde f_\m:= c(\m,h) h^{-\frac d 4} \hat \chi_\m e^{(\phi(\m)-\phi(x))/h}$$
 satisfies $\Vert \tilde f_\m\Vert_{L^2}=1$ and 
$c^\alpha_0(h)$ is defined by  $N(\theta^\alpha_0)=1$. 
Let us now extend the definition of the set 
$H(\m)$ in the following way. Given 
$\alpha\in\aaa$ and $\m\in\widehat \uuu^{(0),II}_\alpha$ we define
\be\label{eq:defhatHm}
\widehat H_\alpha(\m)=
\left\{
\begin{array}{c}
H(\m)\text{ if } \m\in\uuu^{(0)}_\alpha\phantom{*********************}\\
\{\m'\in\widehat E(\alpha)\cap\uuu^{(0)},\phi(\m')=\phi(\hat\m)\}\text{ if } \m\in \widehat\uuu^{(0),II}_\alpha\setminus \uuu^{(0)}_\alpha.
\end{array}
\right.
\ee
Observe that if $\alpha$ is of type II, since $E(\hat\m(\alpha))$ is larger than $\widehat E(\alpha)$, then 
$H(\hat\m(\alpha))$ and $\widehat H_\alpha((\hat\m(\alpha))$ may be different.
From the preceding definition, it follows that for all $\m\in \widehat \uuu^{(0),II}_\alpha$, $c(\m,h)$ admits a classical expansion
$c(\m,h)=\sum_kh^k \gamma_k(\m)$ with  
\be\label{eq:gamma0gen}
\gamma_0(\m)^{-2}=\pi^{\frac d 2} \sum_{\m'\in \widehat H_\alpha(\m)}\vert \det\Hess\phi(\m')\vert^{-\frac 12}.
\ee
Therefore, we can compute the constant $c^\alpha_0(h)$, and we get
$$c^\alpha_0(h)=\pi^{-\frac d 4}\big(\sum_{\nu\in\widehat \uuu^{(0),II}_\alpha}\sum_{\m'\in \widehat H_\alpha(\nu)}\vert \det\Hess\phi(\m')\vert^{-\frac 12} \big)^{-\frac 12}+\ooo(h).$$
Here the index $\alpha$ is used to indicate that the function is associated to $\uuu^{(0),II}_\alpha$.
\end{itemize}
\end{itemize}
\begin{lemma}\label{lem:graamschmidt}
There exist some functions $\theta^\alpha_1,\ldots,\theta^\alpha_{q_\alpha^{II}}\in\fff_{\alpha}$ such that the following hold true:
\begin{itemize}
\item[i)]  $\{\theta^\alpha_j,j=0,\ldots,q_\alpha^{II}\}$ is an orthonormal basis of $\fff_{\alpha}$
\item[ii)] the functions $\theta^\alpha_j$ admit a classical expansion 
$$\theta^\alpha_j=\sum_{k\geq 0}h^k\theta^{\alpha,k}_j$$
and for all $j\geq 1$, the leading terms $\theta^{\alpha,0}_j$ are orthogonal to the function 
$\theta_0^{\alpha,0}(\m)=\frac {c^\alpha_0(0)}{\gamma_0(\m)}$.
\end{itemize}
\end{lemma}
\bp
First observe that $\theta^\alpha_0$ admits a classical expansion
$\theta^\alpha_0\sim \sum_{j\geq 0}h^j\theta^{\alpha,j}_0$ with $\theta_0^{\alpha,0}(\m)=\frac {c^\alpha_0(0)}{\gamma_0(\m)}$.
 Since $(\theta^{\alpha,0}_0)^\bot$ is a $q_\alpha^{II}$ dimensional subspace of $\fff_\alpha$, it admits an orthonormal basis
 $(\tilde\theta_j^{\alpha,0})$ independent of $h$. Then the function $\tilde \theta_j^\alpha$ defined by 
 $$
 \tilde\theta_j^\alpha:=\tilde\theta_j^{\alpha,0}-\<\tilde\theta_j^{\alpha,0},\theta_0^\alpha\>\theta_0^\alpha
 $$
 form a basis of $(\theta_0^\alpha)^\bot$. Moreover, the $  \tilde\theta_j^\alpha$ admit a classical expansion and since
 $\<\tilde\theta_j^{\alpha,0},\theta_0^\alpha\>=\ooo(h)$ for any $j$, they satisfy 
 $$
 \< \tilde\theta_j^\alpha, \tilde\theta_k^\alpha\>=\delta_{jk}+\ooo(h^2).
 $$
 Defining the $(\theta_j^\alpha)$ as the Graam-Schmidt orthonormalization of the $(\tilde\theta_j^\alpha)$, we get the announced result.
\ep

Observe that since $\uuu^{(0),II}_{\alpha}$ has $q_\alpha^{II}$ elements, the functions $\theta^\alpha_1,\ldots,\theta^\alpha_{q_\alpha^{II}}$ can also be indexed by $\uuu^{(0),II}_{\alpha}$ using any arbitrary bijection.
We end up with a family of functions $(\theta^\alpha_\m)_{\m\in\uuu^{(0),II}_{\alpha}}$ and for convenience we will also denote 
$\theta^\alpha_{\hat \m}=\theta^\alpha_0$.
Then, we define the $q_\alpha^{II}$ quasimodes associated to the $\m\in \uuu^{(0),II}_{\alpha}$ by
\begin{equation}\label{eq:definQMGraam}
f_\m^{(0)}(x)=h^{-\frac d 4}\sum_{\m'\in\widehat\uuu_{\alpha}^{(0),II}}\theta^\alpha_\m(\m')c(\m',h)\hat\chi_{\m'}(x)e^{(\phi(\m)-\phi(x))/h}
\end{equation}
where the normalization factor $c(\m',h)$ is defined above and insures that 
$$\Vert c(\m',h)h^{-\frac d 4}\hat\chi_{\m'}(x)e^{(\phi(\m)-\phi(x))/h}\Vert_{L^2}=1.$$

Before going further and as a preparation for the final analysis we would like to write the quasimode given by \eqref{eq:definQMusual} and 
\eqref{eq:definQMGraam} in the same fashion. 
For this purpose, we define  $\widehat \uuu_{\alpha}^{(0)}$ 
by 
\begin{equation}\label{eq:definhatU0}
\widehat \uuu_{\alpha}^{(0)}=\uuu_{\alpha}^{(0),I}\cup\widehat\uuu_{\alpha}^{(0),II}
\end{equation}
with the convention that $\widehat\uuu_{\alpha}^{(0),II}=\emptyset$ if $\uuu_{\alpha}^{(0),II}=\emptyset$
(observe that  $\widehat \uuu_{\alpha}^{(0)}$ is equal to the set $\puu^{(0)}_\alpha$ defined in \ref{eq:defparenU} if and only if $\uuu^{(0),II}_\alpha\neq\emptyset$).
Then, we define $\theta^\alpha_\m(\m')$ for any $\m\in\uuu_{\alpha}^{(0)}, \m'\in\widehat \uuu_{\alpha}^{(0)}$ in the following way:
\begin{itemize}
\item[-] if $\m\in\uuu_{\alpha}^{(0),II}$ and $\m'\in\widehat\uuu_{\alpha}^{(0),II}$, we keep the above definition.
\item[-] otherwise,  we set 
\end{itemize}
\begin{equation}\label{eq:defin-theta-usual}
\theta_\m^\alpha(\m')=\delta_{\m,\m'}.
\end{equation}
Then formula \eqref{eq:definQMusual} and \eqref{eq:definQMGraam} can be summarized in 
\begin{equation}\label{eq:definQMgeneral} 
f_\m^{(0)}(x)=h^{-\frac d 4}\sum_{\m'\in\widehat\uuu_{\alpha}^{(0)}}\theta^\alpha_\m(\m')c(\m',h)\hat\chi_{\m'}(x)e^{(\phi(\m)-\phi(x))/h}
\end{equation}
with $\widehat\uuu_{\alpha}^{(0)}$ and $\theta^\alpha$ as above.
\begin{defin}\label{definmatT}
For any $\alpha\in\ala$, let us denote by $\Tr^\alpha\in\Mr(\uuu^{(0)}_\alpha,\widehat\uuu^{(0)}_\alpha)$ the matrix given by
$$
\Tr^\alpha=(\theta^\alpha_\m(\m')))_{\m'\in\widehat\uuu^{(0)}_\alpha,\m\in\uuu^{(0)}_\alpha}
$$
\end{defin}
Let us remark that if all points of $\uuu_\alpha^{(0)}$ are of type I, then $\Tr^\alpha$ is just the $q_\alpha\times q_\alpha$ identity matrix, whereas if
$\uuu_\alpha^{(0),II}\neq\emptyset$ it is a $(q_\alpha+1)\times q_\alpha$ matrix.
Observe also that the partitions $\uuu^{(0)}_\alpha=\uuu^{(0),I}_\alpha\sqcup\uuu^{(0),II}_\alpha$ and 
$\widehat\uuu^{(0)}_\alpha=\uuu^{(0),I}_\alpha\sqcup\widehat\uuu^{(0),II}_\alpha$ induce decompositions of the corresponding vector spaces
\be
\Fr(\uuu^{(0)}_\alpha)=\Fr(\uuu^{(0),I}_\alpha)\oplus\Fr(\uuu^{(0),II}_\alpha)
\ee
and 
\be
\Fr(\widehat\uuu^{(0)}_\alpha)=\Fr(\uuu^{(0),I}_\alpha)\oplus\Fr(\widehat\uuu^{(0),II}_\alpha).
\ee
From the above construction, one deduces that in a suitable basis the matrix $\Tr^\alpha$ is block diagonal with $\Id$ on the upper-left corner and a certain orthogonal matrix in the
lower-right corner. More precisely,  there exists an orthogonal matrix
 $\wideparen\Tr^\alpha\in\Mr(\uuu^{(0),II}_\alpha,\widehat\uuu^{(0),II}_\alpha)$ such that for any 
 $f=f^I+f^{II}$ with $f^I\in\Fr(\uuu^{(0),I}_\alpha)$ and $f^{II}\in\Fr(\widehat\uuu^{(0),II}_\alpha)$, one has
\be\label{eq:Trdiagblock}
\Tr f(\m)=f^I(\m)+(\wideparen\Tr^\alpha f^{II})(\m).
\ee
Moreover, the matrix $\wideparen\Tr^\alpha$ is just the matrix $(\theta_\m^\alpha(\m'))_{\m\in \uuu^{(0),II}_\alpha, \m'\in \widehat\uuu^{(0),II}_\alpha}$
whose coefficients are given by Lemma \ref{lem:graamschmidt}. In particular, $\Ran \wideparen\Tr^\alpha=(\R \theta_0^\alpha)^\bot$ where 
$ \theta_0^\alpha$ is defined by \eqref{eq:deftheta0}.\\

For any $\m\in\uuu^{(0)}$, let us introduce the set $F(\m)$ defined as follows. 
If $\m=\underline\m$, let $F(\underline\m)=X$. If $\m\in\ulu^{(0),I}:=\ulu^{(0)}\cap\uuu^{(0),I}$, let 
$F(\m)=\overline{E(\m)}$
and if $\m\in\ulu^{(0),II}:=\ulu^{(0)}\cap\uuu^{(0),II}$ let
\begin{equation}\label{eq:defFmtypeII}
F(\m)=\big(\bigcup_{\m'\in\uuu^{(0),II}_\alpha}\overline{E(\m')}\big)\cup\overline {\widehat E(\m)}
\end{equation}
where $\alpha$ is such that $\m\in \uuu^{(0)}_\alpha$. Observe that we always have $\overline {E(\m)}\subset F(\m)$. 
\begin{proposition}\label{prop:propelemF}
Let   $\m,\m'\in \uuu^{(0)}$ be such that $\m\neq\m'$. The following hold true
\begin{enumerate}
\item[i)] if $\m\rrr\m'$ then
\begin{itemize}
\item[i.a)] if $\m$ or $\m'$ is of type I, then $F(\m)\cap F(\m')\subset \vvv^{(1)}.$
\item[i.b)] if $\m$ and $\m'$ are both of type II, then $F(\m)=F(\m')$.
\end{itemize}
\item[ii)] If $\m'\notin\Cl(\m)$, then 
\begin{itemize}
\item[ii.a)] if $\sigma(\m)=\sigma(\m')$, then $F(\m)\cap F(\m')=\emptyset$.
\item[ii.b)] if $\sigma(\m)>\sigma(\m')$, then either $F(\m)\cap F(\m')=\emptyset$ or 
$F(\m')\subset \mathring{F}(\m)$
\end{itemize}
\end{enumerate}
\end{proposition}
\bp
Let $\m\rrr\m'$ with $\m\neq \m'$. As in the proof of Proposition \ref{prop:propelemE}, one has necessarily $\m,\m'\neq\underline\m$. 
Assume first that $\m$ is of type I. Then $F(\m)=\overline{E(\m)}$. If $\m'$ is also of type I, then
$F(\m')=\overline{E(\m')}$. Moreover since $\m\neq\m'$, it follows from i.a) of Proposition \ref{prop:propelemE} that 
$E(\m)\cap E(\m')=\emptyset$. Therefore, $F(\m)\cap F(\m')$ is either empty or is reduced to a union of saddle points which are separating by definition. 
If $\m'$ is of type II, the same proof works. This completes the proof of i.a).

Suppose now that $\m$ and $\m'$ are both of type II. Since $\m\rrr\m'$, it follows that $\widehat E(\m)=\widehat E(\m')$ and hence $F(\m)=F(\m')$
which shows i.b).

Suppose now that $\m'\notin \Cl(\m)$. Consider first the case where $\sigma(\m)=\sigma(\m')$. Then, one has necessarily 
$F(\m)\cap F(\m')=\emptyset$ otherwise we would have $\m\rrr\m'$.

Suppose now that $\sigma(\m)>\sigma(\m')$ and that $F(\m)\cap F(\m')\neq\emptyset$.
If $\m=\underline\m$, then $F(\m)=X$ and the conclusion is obvious. Suppose now that $\m\in\ulu^{(0)}$ and consider first the case where
$\m$ and $\m'$ are of type I. Then $F(\m)=\overline{E(\m)}$ and $F(\m')=\overline{E(\m')}$ and since $\sigma(\m)>\sigma(\m')$ it follows that  $E(\m)\cap E(\m')\neq\emptyset$. Hence  ii.a) of  Proposition \ref{prop:propelemE} shows that $E_-(\m')\subset E(\m)$ which yields $F(\m')\subset E(\m)= \mathring{F}(\m)$. If $\m$ is of type I and $\m'$ of type II, then
one has $E(\m)\cap \tilde E\neq\emptyset$ with either $\tilde E=E(\m'')$ for some  $\m''\in \Cl(\m')$ or $\tilde E=\widehat E(\m')$.
As before,  $E(\m)$ contains the connected component of $\{\phi<\bsigma(\m)\}$ that contains $\tilde E$ and the same proof works.

Let us now suppose that $\m$ is of type II and $\m'$ is of type I. Then $E(\m')\cap \tilde E\neq\emptyset$ with either $\tilde E=E(\m'')$ 
for some $\m''\in \Cl(\m)$ or $\tilde E=\widehat E(\m)$. In both cases one sees easily that 
$E_-(\m')\subset \tilde E$ which proves the result.

The case where both $\m$ and $\m'$ are of type II is left to the reader.
\ep

Let us now give some informations on the support of the quasimodes. For $\m\in\ulu^{(0)}$, let us introduce the set 
\begin{equation}\label{eq:definFepsilondetla}
F_{\epsilon,\tilde\epsilon,\delta}(\m)=\left(\big(F(\m)\setminus\bigcup_{\s\in\vvv^{(1)}\cap\partial F(\m)}B(\s,\epsilon)\big)+\overline {B(0,\tilde\epsilon)}\right)
\bigcup\left(\bigcup_{\s\in(\uuu^{(1)}\setminus\vvv^{(1)})\cap\partial F(\m)}\overline{B(\s,\delta)}\right)
\end{equation}
If  $\m$ is of type I, it is clear that $F_{\epsilon,\tilde\epsilon,\delta}(\m)=\overline{E_{\epsilon,\tilde\epsilon,\delta}(\m)}$ and if $\m$ is of type II, one has
$$
F_{\epsilon,\tilde\epsilon,\delta}(\m)=\overline{\widehat E_{\epsilon,\tilde\epsilon,\delta}(\alpha)}\cup\big(\cup_{\m'\in\uuu^{(0),II}_\alpha}\overline{E_{\epsilon,\tilde\epsilon,\delta}(\m)}\big)
$$
From the above construction one deduces the following proposition.
\begin{proposition}\label{prop:suppQM0}
There exists $\epsilon_0,\delta_0>0$ such that for all $0<\delta<\delta_0$ and all $0<\tilde\epsilon<\epsilon/4<\epsilon_0/4$,  the following hold true:
\begin{enumerate}
\item[i)] for any $\m,\m'\in \uuu^{(0)}$
$$
F(\m)\cap F(\m')=\emptyset\Longrightarrow F_{\epsilon,\tilde\epsilon,\delta}(\m)\cap F_{\epsilon,\tilde\epsilon,\delta}(\m')=\emptyset
$$
\item[ii)]
 for any $\alpha\in\ala$ and $\m\in\uuu^{(0)}_\alpha$,  one has $\supp(f_\m^{(0)})\subset F_{\epsilon,2\tilde\epsilon, 2\delta}(\m)$ and 
$$
\forall \s\in\uuu^{(1)}\setminus (\vvv^{(1)}\cap\partial F(\m)),\;\; d_\phi f_\m^{(0)}=0 \text{ in } B(\s,\delta).
$$
\end{enumerate}

\end{proposition}
\bp
Observe that 
$$F_{\epsilon,2\tilde\epsilon,2\delta}(\m)\subset F(\m)+\overline{B(0,2\max(\delta,\tilde\epsilon))}$$
Since $F(\m)$ and $F(\m')$ are compact, the first point of the proposition immediately follows.
The second point of the proposition is a direct consequence of c) of Proposition \ref{prop:propcutoff}.
\ep

Recall that the functions $f^{(0)}_\m,\;\m\in\uuu^{(0)}$  depend on  $\epsilon,\tilde\epsilon,\delta$ via the definition of the cut-off function 
$\chi_m$. This family is quasi-orthonormal in the following sense.
\begin{proposition}\label{prop:orthogQM0}
There exists $\epsilon_0,\delta_0,\beta>0$ such that for all $0<\delta<\delta_0$, $0<\tilde\epsilon<\epsilon/4<\epsilon_0/4$
and  all $\m,\m'\in\uuu^{(0)}$, one has 
$$\<f_\m^{(0)},f_{\m'}^{(0)}\>=\delta_{\m,\m'}+\ooo(e^{-\beta/h}).$$
\end{proposition}
\bp Throughout the proof, we assume that $0<\delta<\delta_0$ and $0<\tilde\epsilon<\epsilon/4<\epsilon_0/4$ as in Proposition \ref{prop:propcutoff} and we decrease $\epsilon_0,\delta_0$ if necessary.

Let $\m,\m'$ be two minima.

$\bullet$ Consider first the case where  $\m\rrr\m'$. If $\m=\underline\m$, one has necessarily $\m'=\m$ and hence 
$\<f_\m^{(0)},f_{\m'}^{(0)}\>=\Vert f_\m^{(0)}\Vert^2=1$ by construction.
Consider now the case where $\m,\m'\neq\underline\m$ and 
suppose first that $\m$ or $\m'$  is of type I. If $\m=\m'$, the definition of $c(\m,h)$ shows that $\Vert f_\m\Vert=1$. If $\m\neq\m'$, it follows from
 ii) of Proposition \ref{prop:suppQM0}  that $f^{(0)}_\m$ and $f_{\m'}^{(0)}$ are supported in  $F_{\epsilon,2\tilde\epsilon,2\delta}(\m)$ and $F_{\epsilon,2\tilde\epsilon,2\delta}(\m')$ respectively. 
Moreover, thanks to i) of Proposition \ref{prop:propelemF}, one has $F(\m)\cap F(\m')\subset \vvv^{(1)}\cap \partial F(\m)$. Hence, 
one can chose $\epsilon_0$ sufficiently small, so that 
$F_{\epsilon,2\tilde\epsilon,2\delta}(\m)\cap F_{\epsilon,2\tilde\epsilon,2\delta}(\m')=\emptyset$.
Therefore, $\supp(f^{(0)}_\m)\cap \supp(f^{(0)}_{\m'})=\emptyset$ and hence $f^{(0)}_\m$ and $f_{\m'}^{(0)}$ are orthogonal.

Suppose now that $\m$ and $\m'$ are both of type II. Then, we can write
\begin{equation*}
\begin{split}
f_\m^{(0)}(x)&=h^{-\frac d 4}\sum_{\nu_1\in\widehat\uuu_{\alpha}^{(0)}}\theta_\m(\nu_1)c(\nu_1,h)\hat\chi_{\nu_1}(x)e^{(\phi(\hat\m(\m))-\phi(x))/h}\\
f_{\m'}^{(0)}(x)&=h^{-\frac d 4}\sum_{\nu_2\in\widehat\uuu_{\alpha}^{(0)}}\theta_{\m'}(\nu_2)c(\nu_2,h)\hat\chi_{\nu_2}(x)e^{(\phi(\hat\m(\m))-\phi(x))/h}.
\end{split}
\end{equation*}
Since for $\nu_2\neq\nu_1$, $\hat\chi_{\nu_1}$ and $\hat\chi_{\nu_2}$ have again disjoint support for $\epsilon_0,\delta_0>0$ small enough, we get
\begin{equation*}
\begin{split}
\<f_\m^{(0)},f_{\m'}^{(0)}\>&=h^{-\frac d 2}\sum_{\nu\in\widehat\uuu_{\alpha}^{(0)}}\theta_\m(\nu)\theta_{\m'}(\nu)
\vert c(\nu,h)\vert^2\int_X\vert\hat\chi_\nu(x)\vert^2e^{2(\phi(\hat\m(\m))-\phi(x))/h}dx\\
&=\<\theta_\m,\theta_{\m'}\>_{\fff_\alpha}=\delta_{\m,\m'}
\end{split}
\end{equation*}
This shows in particular that $\Vert f_\m^{(0)}\Vert_{L^2}=1$ for all $\m\in\uuu^{(0)}$.

$\bullet$ Suppose now, that $\m'\notin\Cl(\m)$ (in particular $\m\neq\m'$). If $\sigma(\m')=\sigma(\m)$ then $F(\m)\cap F(\m')=\emptyset$
thanks to ii.a) of Proposition  \ref{prop:propelemF} and i) of Proposition \ref{prop:suppQM0} implies that   $F_{\epsilon,2\tilde\epsilon,2\delta}(\m)\cap F_{\epsilon,2\tilde\epsilon,2\delta}(\m')=\emptyset$. Then, the first part of ii) of Proposition \ref{prop:suppQM0}
 proves that $f_\m^{(0)}$ and $f_{\m'}^{(0)}$ are orthogonal.

Consider now the case where $\sigma(\m)\neq\sigma(\m')$, let say $\sigma(\m)>\sigma(\m')$. 
From ii.b) of Proposition  \ref{prop:propelemF}, we know that either $F(\m')$ is disjoint from $F(\m)$ or $F(\m')\subset \mathring{F}(\m)$. In the first case, we get immediately
$\<f_\m^{(0)},f_{\m'}^{(0)}\>=0$ by the same argument as before. Suppose we are in the second situation, that is $F(\m')\subset \mathring{F}(\m) $. By definition, we have $\phi(\m)\leq \phi(\m')$, and by taking $\epsilon_0, \delta_0>0$ small enough we can insure that $F_{\epsilon,2\tilde\epsilon,2\delta}(\m')\subset \mathring{F}_{\epsilon,2\tilde\epsilon,2\delta}(\m)$.

Suppose first that $\phi(\m)<\phi(\m')$. A priori we don't know if $\m,\m'$ are of type I or II. However, since 
$F_{\epsilon,2\tilde\epsilon,2\delta}(\m')\subset \mathring{F}_{\epsilon,2\tilde\epsilon,2\delta}(\m)$, then 
$$
\<f_\m^{(0)},f_{\m'}^{(0)}\>=\int_{F_{\epsilon,2\tilde\epsilon,2\delta}(\m')}f_\m^{(0)}(x)f_{\m'}^{(0)}(x) dx.
$$
and
\begin{equation}\label{eq:restrictQM}
(f^{(0)}_\m)_{\vert F_{\epsilon,2\tilde\epsilon,2\delta}(\m')}=\tilde c(\m,h)h^{-\frac d 4}e^{(\phi(\m)-\phi(x))/h}
\end{equation}
 where the constant $\tilde c(\m,h)$ is  uniformly bounded with respect to $h$. This is clear if $\m$ is of type I. If $\m$ is of type II and let say $\m\in \uuu^{(0)}_\alpha$, 
 then $F(\m')\subset \mathring{F}(\m)$ implies that there exists $\nu\in \widehat \uuu^{(0)}_\alpha$ such that $F(\m')\subset E(\nu)$ (or $\widehat E(\nu)$). Then the general formula \eqref{eq:definQMGraam} shows \eqref{eq:restrictQM}.
Moreover, by construction, there exists 
a cut-off function $\psi\in\mathscr{C}_c^\infty(\mathring{F}_{\epsilon,2\tilde\epsilon,2\delta}(\m))$ independent of $h$ such that $\inf_{\supp\psi}\phi=\phi(\m')$ and
$$
\vert f_{\m'}^{(0)}(x)\vert\leq h^{-d/4}\psi(x)e^{(\phi(\m')-\phi(x))/h}
$$
and it follows that 
$$
\vert \<f_\m^{(0)},f_{\m'}^{(0)}\>\vert \leq Ch^{-\frac d 2}\int\psi(x)e^{(\phi(\m')+\phi(\m)-2\phi(x))/h} dx\leq C'h^{-d/2}e^{(\phi(\m)-\phi(\m'))/h}.
$$
Since $\phi(\m')>\phi(\m)$, this proves the result.

It remains to study the case where $\phi(\m)=\phi(\m')$. Let $\alpha,\alpha'\in\aaa$ be such that 
$\m\in \uuu^{(0)}_\alpha$ and $\m'\in \uuu^{(0)}_{\alpha'}$. From the above assumption, one has also $\sigma(\m)>\sigma(\m')$ 
and $F(\m')\subset \mathring{F}(\m)$.
Since $\sigma(\m)>\sigma(\m')$ and $\phi(\m)=\phi(\m')$ then $f_{\m'}^{(0)}$ is necessarily of type II. It has the form \eqref{eq:definQMGraam} and since $F_{\epsilon,2\tilde\epsilon,2\delta}(\m')\subset \mathring{F}_{\epsilon,2\tilde\epsilon,2\delta}(\m)$, then
 \eqref{eq:restrictQM} still holds true. Hence, we get
 \begin{equation}\label{eq:orthog1}
 \<f_\m^{(0)},f_{\m'}^{(0)}\>=\tilde c(\m,h)h^{-\frac d 2}\sum_{\nu\in\widehat\uuu_{\alpha'}^{(0),II}}\theta_{\m'}(\nu)c(\nu,h)\int\hat \chi_{\nu}(x)e^{2(\phi(x)-\phi(\m))/h}dx.
 \end{equation}
 On the other hand, by a standard argument based on Laplace method, we know that there exist $r>0$ and $\beta>0$ such that for all $\nu\in \widehat\uuu_{\alpha'}^{(0),II}$, one has
 \begin{equation*}
 \begin{split}
 h^{-\frac d 2}c(\nu,h)\int\hat \chi_{\nu}(x)e^{2(\phi(x)-\phi(\m))/h}dx&
 =h^{-\frac d 2}c(\nu,h)\sum_{\nu'\in H(\nu)}\int_{B(\nu',r)}e^{2(\phi(x)-\phi(\m))/h}dx +\ooo(e^{-\beta/h})\\
& =h^{-\frac d 2} c(\nu,h)\int\vert\hat \chi_{\nu}(x)\vert^2e^{2(\phi(x)-\phi(\m))/h}dx+\ooo(e^{-\beta/h})\\
&=\frac1{c(\nu,h)}+\ooo(e^{-\beta/h})
\end{split}
\end{equation*}
Plugging this in \eqref{eq:orthog1}, we get
 \begin{equation}\label{eq:orthog2}
 \begin{split}
 \<f_\m^{(0)},f_{\m'}^{(0)}\>&=\tilde c(\m,h)\sum_{\nu\in\widehat\uuu_{\alpha'}^{(0),II}}\theta_{\m'}(\nu)\frac 1{c(\nu,h)}+\ooo(e^{-\beta/h})\\
 &=\frac{\tilde c(\m,h)}{c^{\alpha'}_0(h)}\<\theta_{\m'},\theta_0^{\alpha'}\>_{\fff_{\alpha'}}+\ooo(e^{-\beta/h})
 \end{split}
 \end{equation}
Since $\theta_{\m'}$ is orthogonal to $\theta_0^{\alpha'}$ by construction, the first term of the right hand side above vanishes and we get 
$\<f_\m^{(0)},f_{\m'}^{(0)}\>=\ooo(e^{-\beta/h})$. This completes the proof.
\ep

We end up this section by giving some exponential estimate of the action of $\dh$ on the quasimodes.
\begin{lemma}\label{lem:derivQM0}
There exists $C>0$ such that for all $\epsilon>0$ small enough, we have
$$
\dh f_\m^{(0)}=\ooo(e^{-(S(\m)-C\epsilon)/h})
$$
for all $\m\in\uuu^{(0)}$.
\end{lemma}

\bp
The result is classical, but since the quasimodes $f_\m^{(0)}$ are slightly different from the usual ones, we have to check the estimates.
Let $\m\in\uuu^{(0)}$ and let us compute $\dh f_\m^{(0)}$.

If $\m=\underline\m$, then $\dh f_\m^{(0)}=0$ and there is nothing to do.

Suppose now that $\m\neq\underline\m$. From \eqref{eq:definQMgeneral}, one has 
$$
\dh f_\m^{(0)}(x)=h^{1-\frac d 4}\sum_{\m'\in\widehat\uuu_{\alpha}^{(0)}}\theta_\m(\m')c(\m',h)\nabla\hat\chi_{\m'}(x)e^{(\phi(\m)-\phi(x))/h}
$$
All terms of the above sum corresponding to $\m\in \uuu_{\alpha}^{(0)}$ are $\ooo(e^{-(S(\m)-C\epsilon)/h}$ thanks to b) of Proposition \ref{prop:propcutoff}. 
The only new term is the one corresponding to 
$\hat\m(\m)$. Since $\hat\chi_{\hat\m}\in\ccc^\infty_c(\widehat E_{\epsilon, 2\tilde\epsilon,2\delta})$ and is equal to $1$ on $\widehat E_{\epsilon, \tilde\epsilon,\delta}$, we have again 
$$\phi(x)-\phi(\hat\m)=\phi(x)-\phi(\m)\geq S(\m)-C\epsilon$$ on $\supp(\nabla\hat\chi_{\hat\m})$ and the proof is complete.
\ep

\subsection{Quasimodes for $1$-forms}\label{subsec:QMF1F}
This section is devoted to the quasimodes associated to low lying  eigenvalues of $\Delta_\phi^{(1)}$.
The construction of these quasimodes was done in \cite{HeSj85_01} and we refer to that paper for all the proofs.  
Here,  we just describe the main properties of these functions. In this section $\omega_\s$ denotes a small neighbourhood of $\s\in\uuu^{(1)}$ that may be chosen as small as needed independently of $\epsilon_0$ fixed in previous sections.

 Given any saddle point $\s\in\uuu^{(1)}$,  and any appropriate open neighborhood $\omega_\s$ of $\s$, let 
 $P_{\phi,\s}$ denote the operator  $\Delta_\phi^{(1)}$ restricted to $\omega_\s$  with Dirichlet boundary conditions.
 Let $u_{\s}$ denote a normalized fundamental state of $P_{\phi,\s}$. The quasimodes $f_\s^{(1)}$ are then defined by
\begin{equation} \label{d20}
f_{\s}^{(1)}(x) : = \varepsilon_0\Vert \psi_{\s} u_{\s} \Vert^{- 1} \psi_{\s} (x) u_{\s} (x) ,
\end{equation}
where $\psi_\s$ is a well-chosen $C^{\infty}_{0}$ localization function supported in $\omega_\s$ and equal to $1$ near $\s$ and 
$\varepsilon_0=\pm 1$ will be fixed later. By taking $\omega_\s$ sufficiently small, 
we can insure that the $f_{\s}^{(1)}$ have disjoint supports, 
and thanks to c) of Proposition \ref{prop:propcutoff}, we can also shrink $\omega_\s$ in order that 
\be\label{eq:suppQM1}
\forall \s\in \uuu^{(1)}\setminus\vvv^{(1)},\;\forall \m\in\uuu^{(0)},\;(\s\in\supp(\chi_\m)\Longrightarrow \chi_\m=1 \text{ on }\omega_\s).
\ee
Observe that this choice of $\omega_s$ depends on $\delta_0$ but not on $\epsilon_0$.
From this construction, we immediately deduce that 
\begin{equation} \label{d21}
\< f_{\s}^{ (1)} , f_{\s'}^{ (1)} \> = \delta_{\s, \s'},
\end{equation}
and hence the family  $\{ f_{\s}^{ (1)},\; \s\in\uuu^{(1)} \}$ is a free family of $1$-forms. 
From \cite{He88_01} Proposition 5.2.6, one knows that the eigenvalues of $P_{\phi,\s}$ are exponentially small.
Using Agmon estimates, it follows that
there exists $\beta > 0$ independent of $\varepsilon$ such that
\begin{equation} \label{b37}
\Delta_\phi^{(1)} f_{\s}^{(1)} = \ooo ( e^{- \beta / h} ) .
\end{equation}
Combined with the spectral theorem, this proves that the $n_1$ eigenvalues of $\Delta_\phi^{(1)}$ in $[0,\epsilon_1 h]$ are actually 
$ \ooo ( e^{- \beta / h} )$ (see  \cite{He88_01} Proposition 5.2.5, for details).

Furthermore, Theorem 2.5 of \cite{HeSj85_01} implies that these quasimodes have a WKB writing
\begin{equation} \label{d22}
f_{\s}^{(1)} (x) = \varepsilon_0h^{-\frac d 4}\psi_{\s} (x) b^{(1)}_{\s} ( x , h ) e^{- \phi_{+ , \s} (x) / h} ,
\end{equation}
where $b^{(1)}_{\s} ( x , h )$ is a  $1$-form having a semiclassical asymptotic, and $\phi_{+ , \s}$ is the phase generating the outgoing manifold of $\vert \xi\vert^{2} - \vert \nabla_x \phi (x) \vert^{2}$ at $( \s , 0 )$ (see \cite{DiSj99_01} chapter 3 for details on such constructions). In particular, the phase function $\phi_{+ , \s}$ satisfies the eikonal equation $\vert \nabla_{x} \phi_{+ , \s} \vert^{2} = \vert \nabla_{x} \phi \vert^{2}$ and $\phi_{+ , \s} (x) \asymp \vert x - \s\vert^{2}$ near $\s$ (the notation $\asymp$ was defined in the paragraph before section 2.1). For other properties of $\phi_{+ , \s}$ we refer to \cite{HeSj85_01}.

\subsection{Projection onto the eigenspaces} \label{z4}

The next step in our analysis is to project the preceding quasimodes onto the generalized eigenspaces associated to exponentially small eigenvalues. Recall that we have built in the preceding section quasimodes $f_\m^{(0)}$, $\m\in\uuu^{(0)}$ with good orthogonality properties. To each of these  quasimode we will  associate a function in  $E^{(0)}$, the eigenspace associated to $o(h)$ eigenvalues. For this, we first define the spectral projector
\begin{equation} \label{d4}
\Pi^{(0)} = \frac{1}{2\pi i} \int_\gamma (z-\Delta_\phi^{(0)})^{-1} d z ,
\end{equation}
where $\gamma = \partial B ( 0 , \epsilon_0 h  )$ and  $\epsilon_0 > 0$ is such that $\sigma(\Delta_\phi)\cap[0,2\epsilon_0 h]\subset[0,e^{-C/h}]$. 
From the  fact that $\Delta_\phi^{(0)}$ is selfadjoint, we get that
\begin{equation*}
\Vert \Pi^{(0)}\Vert = 1 .
\end{equation*}
We now introduce the projection of the quasimodes constructed above,
$
e_\m^{(0)} =\Pi^{(0)} ( f_\m^{(0)} ) 
$.
We have the following
\begin{lemma}\sl \label{lem:projectQM0}
The system $( e_\m^{(0)} )_{\m\in\uuu^{(0)}}$ is free and spans $E^{(0)}$. Besides, there exists $\beta > 0$ independent of $\epsilon_0$ such that
for all $0<\tilde\epsilon<\epsilon/4<\epsilon_0/4$, one has
\begin{equation*}
e_\m^{(0)} = f_\m^{(0)} + \ooo ( e^{- \beta / h} ) \qquad \text{and} \qquad \< e_{\m}^{(0)} , e_{\m'}^{(0)} \> = \delta_{\m , \m '} + \ooo ( e^{-\beta/h} ) .
\end{equation*}
for all $\m,\m'\in\uuu^{(0)}$.
\end{lemma}
\bp The argument is very classical. We recall it for reader's convenience.
One has 
\be
\begin{split}
e_\m^{(0)}-f_\m^{(0)}&=(\Pi^{(0)}-\Id)f_\m^{(0)}=\frac{1}{2\pi i} \int_\gamma ((z-\Delta_\phi^{(0)})^{-1}-z^{-1})f_\m^{(0)} d z\\
&=\frac{1}{2\pi i} \int_\gamma (z-\Delta_\phi^{(0)})^{-1}z^{-1} \Delta_\phi^{(0)}f_\m^{(0)} d z
\end{split}
\ee
Since $(z-\Delta_\phi^{(0)})^{-1}=\ooo(h^{-1})$ on $\gamma $, it follows from Lemma \ref{lem:derivQM0} that 
$e_\m^{(0)}-f_\m^{(0)}=\ooo(^{-\beta/h})$ for some $\beta >0$. This proves the first point. Combining this information with 
 Proposition \ref{prop:orthogQM0} we get immediately the second point.
\ep

We can do a similar study for $\Delta_\phi^{(1)}$, for which we know that the $n_1$ eigenvalues lying in $[0,\epsilon_1 h]$ are 
actually $\ooo ( e^{-\alpha'/h} )$.
To the family of quasimodes $( f_{\s}^{(1)})_{\s\in\uuu^{(1)}}$, we now associate a family of functions in $E^{(1)}$, the eigenspace associated 
to eigenvalues of
${\Delta_\phi^{(1)}}$ in $[0,\epsilon_1 h]$. Thanks to the spectral properties of the selfadjoint operator $\Delta_\phi^{(1)}$, its spectral projector onto $E^{(1)}$ is given by
\begin{equation} \label{d11}
\Pi^{(1)} = \frac{1}{2 \pi i} \int_\gamma ( z - \Delta_\phi^{(1)} )^{- 1} d z ,
\end{equation}
where $\gamma = \partial B ( 0 , \varepsilon_1 h )$ with  $\epsilon_1$  defined above. In the sequel, we denote
$e_\s^{(1)} = \Pi^{(1)} ( f_{\s}^{(1)} ) $.
The family $(e_\s^{(1)})_{\s}$ satisfies the following estimates
\begin{lemma}\sl \label{lem:projectQM1}
The system $( e_\s^{(1)} )_{\s\in\uuu^{(1)}}$ is free and spans $E^{(1)}$. Besides, we have
\begin{equation*}
e_\s^{(1)} = f_\s^{(1)} + \ooo ( e^{-\beta'/h} ) , \qquad \text{and} \qquad \,\< e_{\s}^{(1)} , e_{\s '}^{(1)} \> = \delta_{\s,\s'} + \ooo ( e^{-\beta'/h} ) .
\end{equation*}
with $\beta'>0$ independent of $\varepsilon$.
\end{lemma}
\bp
Using  the orthonormality of the $f^{(1)}_j $ and \eqref{b37}, the proof is the same as that of Lemma \ref{lem:projectQM0}.
\ep

\section{Preliminary for singular values analysis}\label{sec:PFSVA}

This section is a preparation to the study of  the singular values of the operator $\lll:E^{(0)}\rightarrow E^{(1)}$ defined below  \eqref{eq:actionderivEp}.
We simplify the forthcoming study by  several reductions and changes of basis. 
Let us denote by $\lll^\pi$ the $n_1\times n_0$ matrix given by
\be\label{eq:defLpi}
\lll^\pi_{\s,\m}=\<e^{(1)}_\s,\dh e^{(0)}_\m\>,\;\forall \s\in\uuu^{(1)},\m\in\uuu^{(0)}
\ee
with $e^{(1)}_\s$, $e^{(0)}_\m$ defined in the preceding section.
Since $(e^{(0)}_\m)$ and $(e^{(1)}_\s)$ are almost orthonormal bases (thanks to Lemma \ref{lem:projectQM0} and \ref{lem:projectQM1}), this matrix is close to the matrix of the operator $\lll$ in these bases. We first work on the matrix $\lll^\pi$.

Recall that $\underline\m$ denotes the absolute minimum of $\phi$ associated to the connected component $E(\underline\m)=X$.
Since $\Delta_\phi^{(0)}e_{\underline\m}=0$, the non zero singular values of $\lll^\pi$ are exactly the singular values of  the reduced matrix $\lll^{\pi,'}$ defined by $\lll^{\pi,'}_{\s,\m}=\lll^{\pi}_{\s,\m}$ for all $\s\in\uuu^{(1)},\m\in\ulu^{(0)}$
with $\ulu^{(0)}=\uuu^{(0)}\setminus\{\underline\m\}$.

\begin{lemma}\label{lem:approxSV} 
There exists $\beta''>0$ such that for $\epsilon>0$ sufficiently small, one has 
$$
\lll^{\pi,'}_{\s,\m}=\<f^{(1)}_\s,\dh f^{(0)}_\m\>+\ooo(e^{-(S(\m)+\beta'')/h})
$$
for all $\s\in\uuu^{(1)},\m\in\ulu^{(0)}$.
\end{lemma}

\bp
The trick to get the good error estimate above is now well-known (see for instance proof of Prop. 5.8 in  \cite{HeHiSj11_01}) but we recall the proof for reader's convenience. Let $\s\in\uuu^{(1)},\m\in\ulu^{(0)}$, then thanks to \eqref{eq:intertwin} we have
\begin{equation*}
\begin{split}
\<e^{(1)}_\s,\dh e^{(0)}_\m\>&=\<e^{(1)}_\s,\dh \Pi^{(0)}f^{(0)}_\m\>=\<e^{(1)}_\s,\Pi^{(1)}\dh f^{(0)}_\m\>=\<e^{(1)}_\s,\dh f^{(0)}_\m\>\\
&=\<f^{(1)}_\s,\dh f^{(0)}_\m\>+\<e^{(1)}_\s-f^{(1)}_\s,\dh f^{(0)}_\m\>
\end{split}
\end{equation*}
But from Lemma \ref{lem:derivQM0}, \ref{lem:projectQM1} and Cauchy-Schwarz inequality one gets
$$
\vert \<e^{(1)}_\s-f^{(1)}_\s,\dh f^{(0)}_\m\>\vert\leq C e^{-(\beta'+S(\m)-C\epsilon)/h}.
$$
Since $\beta'$ is independent of $\epsilon$, one can conclude by taking $\epsilon$ small enough and $\beta''=\beta'/2$.
\ep

Let us denote $\lll^{bkw}\in\Mr(\uuu^{(0)},\uuu^{(1)})$ the matrix defined by
\begin{equation}\label{eq:definLbkw}
\lll^{bkw}_{\s,\m}=\<f^{(1)}_\s,\dh f^{(0)}_\m\>,  \forall\s\in\uuu^{(1)},\m\in \uuu^{(0)}.
\end{equation}
Of course, the first column of this matrix is identically zero and it is more interesting to consider the matrix
$\lll^{bkw,'}\in\Mr(\ulu^{(0)},\uuu^{(1)})$ defined by
\begin{equation}\label{eq:definLbkw'}
\lll^{bkw,'}_{\s,\m}=\<f^{(1)}_\s,\dh f^{(0)}_\m\>,  \forall\s\in\uuu^{(1)},\m\in \ulu^{(0)}.
\end{equation}
As we shall see later, the singular values of $\lll^{\pi,'}$ and $\lll^{bkw,'}$ are exponentially close and it is natural to 
study the matrix $\lll^{bkw,'}$. 
For $\s\in\uuu^{(1)}\setminus\vvv^{(1)}$ and $\m\in\ulu^{(0)}$, thanks to ii) of Proposition \ref{prop:suppQM0} 
one has  $\dh f_\m^{(0)}=0$ near $\s$, and hence
\begin{equation}\label{eq:uselesssp}
\<f_\s^{(1)},\dh f_\m^{(0)}\>=0.
\end{equation}
Therefore, the singular values of $\lll^{bkw,'}$ are equal to the singular values of the reduced matrix $\lll^{bkw,''}\in\Mr(\ulu^{(0)},\vvv^{(1)})$
defined by 
\begin{equation}\label{eq:definLr''}
\lll^{bkw,''}_{\s,\m}=\<f_\s^{(1)},\dh f_\m^{(0)}\>,\;\forall \s\in\vvv^{(1)},\m\in\ulu^{(0)}.
\end{equation}
In order to study this matrix, we need to introduce a new enumeration of critical points.
Let us start with few abstract notations.
Assume that $(\iii,\leq)$ and $(\jjj,\leq)$ are two totally ordered sets and let 
$A=(a_{ij})_{i\in\iii,j\in\jjj}$ be the associated matrix (with $i,j$ enumerated in increasing order). Assume that we have partitions $\ppp_\iii,\ppp_\jjj$ of $\iii$ and $\jjj$ respectively
$$
\ppp_\iii=(\iii_1,\ldots,\iii_{N_\iii})\text{ and }\ppp_\jjj=(\jjj_1,\ldots, \jjj_{N_\jjj}).
$$
Assume that each partition admits a total order $\preceq$ (that is we can compare the subsets $\iii_i$). Then we get a total  order  $\preceq$ on $\iii$ (resp. $\jjj$) by using the associated lexicographical order:
$$
i\preceq j\; \text { iff } \;(\exists \iii_\alpha\preceq \iii_\beta,\; i\in \iii_{\alpha}\text{ and }j\in\iii_\beta)\text{ or } (\exists \iii_\alpha,\,i,j\in\iii_\alpha\text{ and }i\leq j).
$$
Hence, there exists a unique  $\alpha:(\iii,\leq)\rightarrow (\iii,\preceq)$ which is strictly increasing (and hence bijective). Similarly, there is a unique 
$\beta:(\jjj,\leq)\rightarrow (\jjj,\preceq)$ which is strictly increasing.
We denote by $A_{\ppp_\iii,\ppp_\jjj}$ the matrix $(a_{\alpha(i),\beta(j)})_{i\in\iii,j\in\jjj}$.
This matrix is obtained from $A$ by intertwining the basis vector, hence it has exactly the same singular values.\\

Let us go back to the matrix $\lll^{bkw,''}$. Consider the partitions of $\ulu^{(0)}$ and $\vvv^{(1)}$ given by 
$$\ppp^{(0)}=\{\uuu^{(0)}_\alpha,\;\alpha\in\ala\}\text{ and }\ppp^{(1)}=\{\vvv^{(1)}_\beta,\;\beta\in\ala\}.$$
At this stage of our analysis, we do not need any specific choice of order on these partitions. We just endow $\ala$ with any total order and  for all $\alpha,\beta\in\ala$ we choose any arbitrary total order on $\uuu^{(0)}_\alpha$ and $\vvv^{(1)}_\beta$.
This gives an order on the above partitions and we 
denote by $\Lr=(\Lr^{\alpha,\beta})_{\alpha,\beta\in\ala}$ the matrix $\lll^{bkw,''}$ associated to these partitions. Observe here that each 
$\Lr^{\alpha,\beta}$ is itself a matrix $\Lr^{\alpha,\beta}=(\Lr^{\alpha,\beta}_{\s,\m})_{\s\in\vvv^{(1)}_\beta,\m\in\uuu^{(0)}_\alpha}$.
\begin{lemma}\label{lem:diagblocLr}
For all $\alpha\neq\beta$, $\Lr^{\alpha,\beta}=0$.
\end{lemma}
\bp
Let $\alpha;\beta\in\ala$ such that $\alpha\neq\beta$ and let $\m\in \uuu^{(0)}_\alpha$ and 
$\s\in \vvv^{(1)}_\beta$. 
If $\bsigma(\alpha)=\bsigma(\beta)$   then $\alpha\neq\beta$ implies that $\s\notin F(\m)$. Shrinking if necessary (by taking $\epsilon_0,\delta_0>0$ small enough) the support of $f_\m^{(0)}$ and $f_\s^{(1)}$, it follows that these functions have disjoint supports so that 
their scalar product vanishes.

If $\bsigma(\alpha)\neq\bsigma(\beta)$, then by construction $\dh f_\m^{(0)}$ is supported near $\{\phi=\bsigma(\alpha)\}$ whereas 
$e_\s^{(1)}$ is supported near $\{\phi=\bsigma(\beta)\}$. Since this two sets are disjoints we get 
$\<f_\s^{(1)},\dh e_\m^{(0)}\>=0$ and the proof is complete.
\ep

From this lemma we deduce that the matrix $\Lr$ admits a block-diagonal structure
\begin{equation}\label{eq:blocdiagD}
\Lr=\diag(\Lr^\alpha,\alpha\in\ala)
\end{equation}
with $\Lr^\alpha:=\Lr^{\alpha,\alpha}$.  Recall from Definition \ref{definmatT}, 
that for any $\alpha\in\ala$, 
 the matrix $\Tr^\alpha\in\Mr(\uuu^{(0)}_\alpha,\widehat\uuu^{(0)}_\alpha)$ is given by 
$\Tr^\alpha=(\theta^\alpha_\m(\m'))_{\m'\in\widehat \uuu^{(0)}_\alpha,\m\in \uuu^{(0)}_{\alpha}}$. We have the following factorization result on $\Lr^\alpha$.
\begin{lemma}\label{lem:factorisD}
We have
$\Lr^{\alpha}=\widehat\Lr^{\alpha} \Tr^{\alpha}$
where the matrix $\widehat\Lr^{\alpha}= (\hat\ell^\alpha_{\s,\m'})_{\s,\m'}\in\Mr(\widehat \uuu^{(0)}_\alpha,\vvv^{(1)}_\alpha)$ is given by
$$
\hat\ell^\alpha_{\s,\m'}=\<f_\s^{(1)},\dh g_{\m'}^{(0)}\>,\;\forall \s\in \vvv^{(1)}_{\alpha},\m'\in\widehat\uuu^{(0)}_{\alpha}
$$
with
$
g_{\m'}^{(0)}(x)=h^{-\frac d4}c(\m',h)\hat\chi_{\m'}(x)e^{\phi(\m')-\phi(x)/h}.
$
\end{lemma}
\bp
Let $\s\in \vvv^{(1)}_{\alpha}, \m\in \uuu^{(0)}_{\alpha}$.
From equation \eqref{eq:definQMgeneral}, one has
\begin{equation*}
\<f^{(1)}_\s,\dh f^{(0)}_\m\>=h^{-\frac d 4}\sum_{\m'\in\widehat\uuu_{\alpha}^{(0)}}\theta^\alpha_\m(\m')c(\m',h)
\<f^{(1)}_\s,hd\hat \chi_{\m'}(x)e^{(\phi(\m)-\phi(x))/h}\>
\end{equation*}
Moreover, the function $\phi$ being constant on $\widehat\uuu_{\alpha}^{(0),II}$, 
we can replace $\phi(\m)$ by $\phi(\m')$ in the above identity and  it follows that 
$$
\<f^{(1)}_\s,\dh f^{(0)}_\m\>=\sum_{\m'\in\widehat\uuu_{\alpha}^{(0)}}\theta^\alpha_\m(\m')\<f^{(1)}_\s,d_{\phi,h}g_{\m'}^{(0)}\>
$$
which is exactly the result to be proved.
\ep

One of the crucial points of our analysis is to compute the coefficient $\hat\ell^\alpha_{\s,\m}$. Given $\m\in \widehat \uuu^{(0)}_\alpha$, we define
\begin{equation}\label{eq:defhphi}
h_\phi(\m)=\Big(\sum_{\m'\in \widehat H_\alpha(\m)}\vert \det\Hess\phi(\m)\vert^{-\frac 12}\Big)^{-\frac 12}
\end{equation}
with $\widehat H_\alpha(\m)$ defined in \eqref{eq:defhatHm}. One has clearly $h_\phi(\m)=\pi^{\frac d 4}\gamma_0(\m)$ with 
$\gamma_0$ given by \eqref{eq:gamma0gen}. Moreover, in the case where 
$H(\m)=\{\m\}$, one has $h_\phi(\m)= \vert\det\Hess\phi(\m)\vert^{\frac 14}$. 
Given $\s\in\vvv^{(1)}$, we denote by $\hat \lambda_1(\s)$ the unique negative eigenvalue of $\Hess \phi(\s)$.
In order to keep uniform notations, we also extend the definition \eqref{eq:defhphi} to saddle points by 
$$
h_\phi(\s)=\vert \det\Hess\phi(\s)\vert^{\frac 14}.
$$
Eventually, we introduce the diagonal matrix $\widehat\Omega^\alpha\in\Mr(\widehat\uuu^{(0)}_\alpha,\widehat\uuu^{(0)}_\alpha)$ defined by
\be\label{eq:defhatOmega}
\widehat\Omega^\alpha f(\m)=e^{-\hat S(\m)/h} f(\m),\;\forall\m\in\widehat\uuu^{(0)}_\alpha
\ee
with $\hat S(\m)=\bsigma(\alpha)-\phi(\m)$.
For $\m\in \uuu^{(0)}_\alpha$, one has of course $\bsigma(\alpha)=\bsigma(\m)$ and hence $\hat S(\m)=S(\m)$ but this fails to be true for $\m=\hat\m(\alpha)$.
We then define the rescaled matrix 
$\widetilde \Lr^{\alpha}=(\tilde\ell^\alpha_{\s,\m})\in\Mr(\widehat\uuu^{(0)}_\alpha,\vvv^{(1)}_\alpha)$ by 
$$
\widehat\Lr^\alpha=\widetilde\Lr^\alpha\widehat\Omega^\alpha
$$
i.e.
\be\label{eq:deftildeL}
\tilde\ell_{\s,\m}^\alpha=e^{\hat S(\m)/h}\hat\ell^\alpha_{\s,\m},\;\forall\s\in\vvv^{(1)}_\alpha,\m\in\widehat\uuu^{(0)}_\alpha.
\ee
Going back to the matrix $\Lr^\alpha$, one has 
$$
\Lr^\alpha=\widetilde\Lr^\alpha\widehat\Omega^\alpha\Tr^\alpha
$$
and using the fact that $\Tr^\alpha f(\m)=f(\m)$ for any $f$ supported on $\uuu^{(0),I}_\alpha$ one gets
\be\label{eq:factorLalpha}
\Lr^\alpha=\widetilde\Lr^\alpha\Tr^\alpha\Omega^\alpha
\ee
with $\Omega^\alpha\in \Mr(\uuu^{(0)}_\alpha,\uuu^{(0)}_\alpha)$ defined by $\Omega^\alpha f(\m)=e^{-S(\m)/h} f(\m)$.
The following lemma  gives 
an asymptotic expansion of the  matrix  $\widetilde\Lr^\alpha$. We recall that $\m_1(\s),\m_2(\s)$ were defined in Lemma \ref{lem:defm1m2}

\begin{lemma}\label{lem:computinter} Let $\alpha\in\ala$ and $\s\in\vvv^{(1)}_\alpha$,  $\m\in\widehat\uuu^{(0)}_\alpha$. The following hold true:
\begin{enumerate}
\item[i)] if $\m\notin\{\m_1(\s),\m_2(\s)\}$, then $\tilde\ell^\alpha_{\s,\m}=0$.
\item[ii)] the coefficients $\tilde\ell^\alpha_{\s,\m}$ admits a classical expansion 
$\tilde\ell^\alpha_{\s,\m}\sim h^{\frac 12} \sum_{k\geq 0}h^{k}\tilde\ell_{\s,\m}^{\alpha,k}$. Moreover, one can chose $\varepsilon_0=\pm 1$ in 
\eqref{d20} in order that  the
leading terms satisfy
\be\label{eq:computinter1}
\tilde\ell_{\s,\m_1(\s)}^{\alpha,0}= \pi^{-\frac 12}\vert\hat\lambda_1(\s)\vert^{\frac 12}
\frac{h_\phi(\m_1(\s))}{h_\phi(\s)}
\ee
and 
in the case where $\m_2(\s)\in\widehat\uuu^{(0)}_\alpha$, 
\be\label{eq:computinter2}
\tilde\ell_{\s,\m_2(\s)}^{\alpha,0}=- \pi^{-\frac 12}\vert\hat\lambda_1(\s)\vert^{\frac 12}
\frac{h_\phi(\m_2(\s))}{h_\phi(\s)}.
\ee
In particular, if $\m_2(\s)\in\widehat\uuu^{(0)}_\alpha$, one has
\be\label{eq:computinter3}
\frac{\tilde\ell_{\s,\m_1(\s)}^{\alpha,0}}{h_\phi(\m_1(\s))}=-\frac{\tilde\ell_{\s,\m_2(\s)}^{\alpha,0}}{h_\phi(\m_2(\s))}
\ee
for all $\s\in\vvv^{(1)}_\alpha$.
\end{enumerate}
\end{lemma}
\bp
Suppose first that $\m\neq\m_1(\s),\m_2(\s)$.
Then, $\supp(\dh g_{\m}^{(0)})=\supp(d\hat\chi_{\m})$ is contained in a small neighborhood 
$\omega$ of $\Gamma(\m)$. Since $\m\neq\m_1(\s),\m_2(\s)$ it follows from 
Lemma \ref{lem:defm1m2} that $\s\notin\omega$ and hence $ \tilde \ell_{\s,\m}^{\alpha}=0$ which proves $i)$.

Let us now compute the coefficients $\tilde \ell_{\s,\m}$ for $\m\in\{\m_1(s),\m_2(\s)\}\cap\widehat\uuu^{(0)}_\alpha$ (observe that this set may be reduced to $\m_1(\s)$). We compute these coefficients in the case where $\m_2(\s)\in \widehat\uuu^{(0)}_\alpha$. If it is not the case, the only non-zero coefficient is $\tilde \ell_{\s,\m_1(\s)}$ that is computed in the same way.
Recall from \eqref{d22}, that the quasimodes on $1$-forms are given by
\begin{equation*} 
f_\s^{(1)}=\varepsilon_0h^{-\frac d 4}\psi_{\s} (x) b^{(1)}_{\s} ( x , h ) e^{- \phi_{+ , \s} (x) / h}.
\end{equation*}
Summing up the construction of \cite{HeKlNi04_01} section 4.2, there exists an open neighborhood $V_\s$ of $\s$ on which one can find a system of local Morse coordinates
$(y,z)\in\R\times\R^{d-1}$ in which $\s$ is the origin and such that the following properties hold true:
\begin{enumerate}
\item in the above coordinate system one has 
$$\phi=\phi(\s)+\frac 12\big(\hat \lambda^1(\s)y^2+\sum_{j=2}^d\hat\lambda_j(\s)z_j^2\big)$$
and
$$
\phi_+=\frac 12\big(-\hat \lambda^1(\s)y^2+\sum_{j=2}^d\hat\lambda_j(\s)z_j^2\big)
$$
where $(\hat\lambda_j(\s))_{j=1,\ldots,d}$ are the eigenvalues of $\operatorname{Hess}(\phi)$ at point $\s$.
\item the amplitude $b^{(1)}_{\s} ( x , h )$ admits a classical expansion
\begin{equation}\label{eq:expQM1}
b^{(1)}_{\s} \sim \sum_{k=0}^\infty h^kw_{\s,k}
\end{equation}
with 
\begin{equation}\label{eq:expQM2}
w_{\s,0}=(-1)^{d-1}\frac{\vert \det\Hess \phi(\s)\vert^{\frac 1 4}}{\pi^\frac d 4} dy \;\text{ on }\, \{z=0\}
\end{equation}\,

\item one can chose the orientation of the $y$ axis so that
$$E(\m_1(\s))\cap V_\s\subset\{y<0\}\cap V_\s \text{ and }E(\m_2(\s))\cap V_\s\subset\{y>0\}\cap V_\s.$$
\end{enumerate}
Moreover, the cut-off function $\chi_\m$ can be constructed so that 
\begin{enumerate}
\item[(4)] in $V_\s$ the functions $\hat\chi_{\m_j}$, $j=1,2$ depend only on the variable $y$,
\end{enumerate}
and one can shrink $\omega_\s$ in order that 
\begin{enumerate}
\item[(5)] $\supp(f_\s^{(1)})$ is contained in $V_\s$.

\end{enumerate}
Observe that the only minor (but important) difference with \cite{HeKlNi04_01} is the property (2), saying that  each $\chi_{m_j}, j=1,2$ is supported in one of the two different half plane $\{y\lessgtr0\}$. Let us now compute the first coefficient in the asymptotic expansion of $\tilde \ell_{\s,\m}^{p,\alpha}$.
Using the above properties, Proposition \ref{prop:propcutoff} and following the computations of \cite{HeKlNi04_01} section 6
we get
\begin{equation*}
\begin{split}
\hat\ell^\alpha_{\s,\m}=\<f_\s^{(1)}&,\dh g_\m^{(0)}\>=\\
&h^{1-\frac d 2}c(\m,h)e(\s,h)\int_{B(\s,\epsilon)}e^{-(\phi_+(x)+\phi(x)-\phi(\m))/h}(\hat \chi_\m'(y)+\ooo(h))dy\wedge dz_2\wedge\ldots\wedge dz_d\\
&+\ooo_\epsilon(e^{-(\phi(\s)-\phi(\m)+c_\epsilon)/h})
\end{split}
\end{equation*}
with 
$$
e(\s,h)=\varepsilon_0(-1)^{d-1}\frac{\vert \det\Hess \phi(\s)\vert^{\frac 1 4}}{\pi^\frac d 4}+\ooo(h)=\varepsilon_0(-1)^{d-1}\pi^{-\frac d 4}h_\phi(\s)+\ooo(h)$$
Using the local form of $\phi$ and $\phi_+$, we get 
\begin{equation*}
\begin{split}
\hat\ell^\alpha_{\s,\m}&=
h^{1-\frac d 2}c(\m,h)e(\s,h)e^{-(\phi(\s)-\phi(\m))/h}\int_{B(\s,\epsilon)}e^{-g_-(z)/h}(\hat\chi'_\m(y)+\ooo(h))dy\wedge dz_2\wedge\ldots\wedge dz_d\\
&+\ooo_\epsilon(e^{-(\phi(\s)-\phi(\m)+c_\epsilon)/h})
\end{split}
\end{equation*}
with $g_-(z)=\sum_{j=2}^d\hat\lambda_j(\s)z_j^2$.
Since $\hat\chi_\m$ depends only on $y$ and $g_-\geq c\nu^2$ on $\vert z\vert_\infty\geq \nu$, the integration domain $B(\s,\epsilon)$ can be replaced by
a smaller one $W_\s=\{\vert y\vert<\epsilon,\;\vert z\vert_\infty\leq\nu_\epsilon\}$ modulo exponentially small error terms. Using also the identity
$\hat S(\m)=\phi(\s)-\phi(\m)$, we get
$$ \hat\ell^\alpha_{\s,\m}=I_\epsilon(h)e^{-\hat S(\m)/h}
+\ooo_\epsilon(e^{-(\hat S(\m)+c_\epsilon)/h})$$
with
\begin{equation*}
\begin{split}
I_\epsilon(h)=h^{1-\frac d 2}c(\m,h)e(\s,h)\int_{W_\s}e^{-g_-(z)/h}(\hat\chi'_\m(y)+\ooo(h))dy\wedge dz_2\wedge\ldots\wedge dz_d
\end{split}
\end{equation*}
The integral in the right hand side can be easily computed by mean of Stoke formula and Laplace method. 
We get
\begin{equation*}
\begin{split}
I_\epsilon(h)&=h^{1-\frac d 2}c(\m,h)e(\s,h)([\hat\chi_\m]_{-\epsilon}^\epsilon+\ooo(h))
\int_{\vert z\vert_\infty\leq\nu_\epsilon}e^{-g_-(z)/h}dz_2\wedge\ldots\wedge dz_d\\
&=h^{\frac 1 2}c(\m,h)e(\s,h)([\hat\chi_\m]_{-\epsilon}^\epsilon+\ooo(h))
\Big(\frac{\pi^{\frac {d-1}2}}{\vert\hat\lambda_2(\s)\ldots\hat\lambda_d(\s)\vert^{\frac 1 2}}\Big)
\end{split}
\end{equation*}
Combining this with the expression of $c(\m,h)$ and $e(\s,h)$, we obtain
\begin{equation*}
\tilde\ell_{\s,\m}^{\alpha,0}=\varepsilon_0(-1)^{d-1}[\hat\chi_\m]_{-\epsilon}^\epsilon \pi^{-\frac 12}\vert\hat\lambda_1(\s)\vert^{\frac 12}
\frac{h_\phi(\m)}{h_\phi(\s)}.
\end{equation*}
We now remark that with our choice of $\hat\chi_m$, one has 
$[\hat\chi_{\m_1}]_{-\epsilon}^\epsilon=-1\text{ and } [\hat\chi_{\m_2}]_{-\epsilon}^\epsilon=1.$
Taking $\varepsilon_0=(-1)^d$,  we get immediately the formula of $ii)$.
\ep

\section{Computation of the approximated singular values}\label{sec:COTASV}

From Lemma \ref{lem:SVdiagmat}, we know that the singular values of a block-diagonal matrix are given by the singular values of each block. Hence, in view of the results of the preceding section, we  study the matrices $\Lr^\alpha$. 
The first step in the analysis is to prove that $\Lr^\alpha$ is injective excepted for $\alpha=\underline\alpha$.

\subsection{Injectivity of the matrix $\Lr^\alpha$}\label{subsec:IOTML}
We first compute the kernel of the matrix $\widetilde\Lr^\alpha$.
\begin{lemma}\label{lem:kernelLalphatilde}
Let $\alpha\in\ala$, then 
\begin{itemize}
\item[-] if $\alpha$ is of type I (that is $\uuu^{(0),II}_\alpha=\emptyset$), then $\widetilde\Lr^{\alpha,0}$ is injective
\item[-] if $\alpha$ is of type II, then  $\Ker( \widetilde\Lr^{\alpha,0})=\R \xi_0$ where 
 $\xi_0\in \R^{\hat q_{\alpha}}\simeq\fff_\alpha$ is defined by
 $$
 \xi_0(\m)=h_\phi(\m)^{-1}
 $$
 for all $\m\in\widehat\uuu^{(0)}_\alpha$.
\end{itemize}
\end{lemma} 
\bp Suppose first that $\alpha$ is of type II.  Let $x\in\Fr_\alpha=\Fr(\widehat\uuu^{(0)}_{\alpha})$ be such that 
$\widetilde\Lr^{\alpha,0}x=0$, then
\be\label{eq:systL}
\sum_{\m\in \widehat\uuu^{(0)}_{\alpha}}\tilde\ell^{\alpha,0}_{\s,\m}x_\m=0,\;\forall \s\in \vvv^{(1)}_{\alpha}.
\ee
From $i)$  of Lemma \ref{lem:computinter} it follows that
$$
\tilde\ell^{\alpha,0}_{\s,\m_1(\s)}x_{\m_1(\s)}=-\tilde\ell^{\alpha,0}_{\s,\m_2(\s)}x_{\m_2(\s)},\;\forall \s\in\vvv^{(1)}_{\alpha}.
$$
Moreover, since $\alpha$ is of type II, then $\m_2(\s)\in \widehat\uuu^{(0)}_\alpha$ for any $\s\in\vvv^{(1)}_\alpha$
and thanks to  \eqref{eq:computinter3}  we get
\begin{equation}\label{eq:kern1}
x_{\m_1(\s)}h_\phi(\m_1(\s))=x_{\m_2(\s)}h_\phi(\m_2(\s)),\;\forall \s\in\vvv^{(1)}_{\alpha}.
\end{equation}
Now, we recall that for any $\s\in\vvv^{(1)}_{\alpha}$, $\m_1(\s)$ and $\m_2(\s)$ are exactly the two minima such that 
$\s=\Gamma_\alpha(\m_1)\cap\Gamma_\alpha(\m_2)$. Therefore, we deduce from \eqref{eq:kern1} that 
$$
\forall \m,\m'\in \widehat \uuu^{(0)}_{\alpha},\;(\Gamma_\alpha(\m)\cap\Gamma_\alpha(\m')\neq\emptyset\Longrightarrow 
h_\phi(\m) x_\m=h_\phi(\m')x_{\m'}).
$$
By definition of the equivalence relation $\rrr$, this implies that  $x_\m h_\phi(\m)$  is constant on $\widehat\uuu^{(0)}_\alpha$, 
which means exactly that $x\in\R\xi_0$.

Suppose now that $\alpha$ is of type I and let  $x\in\Fr(\uuu^{(0)}_\alpha)$ such that $\widetilde \Lr^{\alpha,0}x=0$.   As precedently, one shows that there exists a constant $c$ such that for all $\m\in\uuu^{(0)}_\alpha$, 
$h_\phi(\m)x_\m=c$. 
Recall that the non-empty set $\vvv^{(1),b}_\alpha$ was defined in Lemma \ref{lem:partitionValpha1}. Given $\s_b\in\vvv^{(1),\b}_\alpha$, since $\m_2(\s_b)=\hat\m(\alpha)\notin\widehat\uuu^{(0)}_\alpha$, one has 
$\tilde\ell_{\s_b,\m}^{\alpha,0}=0$ for any $\m\neq \m_1(\s_b)$ and 
$$
\tilde\ell_{\s_b,\m_1(\s_b)}^{\alpha,0}=\pi^{-\frac 12}\vert\hat\lambda_1(\s_b)\vert^{\frac 12}\neq 0
$$
Combined with  \eqref{eq:systL} this shows that  $x_{\m_1(\s_b)}=0$
and   hence $c=0$ which proves that $\Ker( \widetilde \Lr^{\alpha,0})=0$.
\ep

\begin{proposition}\label{prop:injectLalpha} 
Let $\alpha\in\ala$, then the matrix $\wideparen\Lr^\alpha:=\widetilde\Lr^\alpha\Tr^\alpha$ admits a classical expansion 
$\wideparen\Lr^\alpha\sim h^{\frac 12}\sum_j h^j\wideparen\Lr^{\alpha,j}$ and the matrix $\wideparen\Lr^{\alpha,0}$ is injective.
\end{proposition} 
\bp
 Thanks to  Lemma \ref{lem:graamschmidt} and \ref{lem:computinter} the matrices $\widetilde\Lr^\alpha$ and $\Tr^\alpha$ admit some classical expansions
$\widetilde\Lr^\alpha\sim h^{\frac 12}\sum h^j\widetilde\Lr^{\alpha,j}$ and $\Tr^\alpha\sim\sum h^j\Tr^{\alpha,j}$. Therefore,  
$\wideparen\Lr^\alpha$ admits a classical expansions $\wideparen\Lr^\alpha\sim h^{\frac 12}\sum_j h^j\wideparen\Lr^{\alpha,j}$ with $\wideparen\Lr^{\alpha,0}=\widetilde\Lr^{\alpha,0}\Tr^{\alpha,0}$.

Let us now  prove that 
$\wideparen\Lr^{\alpha,0}$ is injective.

Suppose first that $\alpha$ is of type I. Then $\Tr^{\alpha}=\Tr^{\alpha,0}=\Id$ and the result follows immediately from the first part of Lemma 
\ref{lem:kernelLalphatilde}. 

Suppose now that $\alpha$ is of type II and let $x\in\Fr(\uuu^{(0)})$ be such that 
$\widetilde\Lr^{\alpha,0}\Tr^{\alpha,0}x=0$. We decompose 
$x=x^I+x^{II}$ with $x^\bullet$ supported in $\widehat \uuu^{(0),\bullet}$. Thanks to \eqref{eq:Trdiagblock}, one has
$$
\Tr^{\alpha,0} x(\m)=x^I(\m)+(\wideparen \Tr^{\alpha,0} x^{II})(\m)
$$ 
with $\wideparen \Tr^{\alpha,0} :\Fr(\uuu^{(0),II})\rightarrow \Fr(\widehat\uuu^{(0),II})$ 
such that $ \Ran \wideparen \Tr^{\alpha,0}=(\R \theta^{\alpha}_0)^\bot$ where the function $\theta^{\alpha}_0$ is defined by \eqref{eq:deftheta0}.
On the other hand, one has $\ker \widetilde \Lr^{\alpha,0}=\R\xi_0$ and one can decompose
$\xi_0=\xi_0^{I}+\xi_0^{II}$ with $\xi_0^{II}=\theta^{\alpha,0}$. The equation $\widetilde\Lr^{\alpha,0}\Tr^{\alpha,0}x=0$ implies that 
there exists $\lambda\in\R$ such that $\Tr^{\alpha,0} x=\lambda \xi_0$ and hence
$ \wideparen \Tr^{\alpha,0} x^{II}=\lambda\xi_0^{II}$. On the other hand, by construction, $ \Ran \wideparen \Tr^{\alpha,0}=(\xi_0^{II})^\bot$. 
This implies that $\lambda=0$ and proves the result.
\ep

\begin{corollary}
For all $\alpha\in\ala$ the matrix $\Lr^\alpha$ is injective.
\end{corollary}
\bp
This follows directly from the above proposition and the fact that 
\be\label{eq:graded1}
\Lr^\alpha=\widehat \Lr^\alpha\Tr^\alpha=\widetilde\Lr^\alpha\widehat\Omega^\alpha\Tr^\alpha=\wideparen\Lr^\alpha\Omega^\alpha
\ee
 with $\Omega^\alpha$ defined below 
\eqref{eq:factorLalpha} which is invertible.
\ep
\subsection{Graded structure of the matrices $\Lr^\alpha$}\label{subsec:GSOTM} Throughout this section, we assume that $\alpha\in\aaa$ is fixed.
Recall that we defined $\sss_\alpha=S(\uuu^{(0)}_\alpha)$, $p(\alpha)=\sharp \sss_\alpha$ and some integers 
$\nu_1^\alpha<\ldots<\nu^{\alpha}_{p(\alpha)}$ such that 
$$
\sss_\alpha=\{S_{\nu_1^\alpha},\ldots,S_{\nu_{p(\alpha)}^\alpha}\}
$$
with the convention $S_{\nu_1^\alpha}>\ldots>S_{\nu_{p(\alpha)}^\alpha}$. In order to lighten the notation we will drop the indices $\alpha$ and write from now $p=p(\alpha)$, $\nu_j=\nu_j^\alpha$.
To the set of heights $\sss_\alpha$, we can associate a natural partition 
\be\label{eq:gradedpartition}
\widehat \uuu^{(0)}_\alpha=\bigsqcup_{n=1}^{p}\widehat\uuu^{(0)}_{\alpha,n}
\ee
with 
$\widehat\uuu^{(0)}_{\alpha,n}=\{\m\in \widehat\uuu^{(0)}_\alpha,\,\phi(\m)=\bsigma(\alpha)-S_{\nu_n}\}$. We order this partition by deciding that 
$\widehat\uuu^{(0)}_{\alpha,n+1}\prec \widehat\uuu^{(0)}_{\alpha,n}$.
On the other hand, we recall that $\Lr^\alpha=\wideparen\Lr^\alpha\Omega^\alpha$
with $\wideparen\Lr^\alpha=\widetilde\Lr^\alpha\Tr^\alpha$.
Let us 
compute the matrices $\wideparen\Lr^\alpha$ and $\Omega^\alpha$ in the basis given by the above partition of $\widehat\uuu^{(0)}_{\alpha}$. With a slight abuse of notation we still denote by $\wideparen\Lr^\alpha$,  $\Omega^\alpha$ the resulting matrices. 
Since $\hat S(\m)=\bsigma(\alpha)-S_{\nu_k}$ on $\uuu^{(0)}_{\alpha,k}$, it follows from \eqref{eq:defhatOmega} that in the above partition, 
 the matrix $\Omega^\alpha$ writes
\be\label{eq:Omegaalphadiag}
\Omega^\alpha=
\left(
\begin{array}{ccccc}
e^{-S_{\nu_p}/h}I_{r_p}&0&\hdots&\hdots&0\\
0&e^{-S_{\nu_{p-1}}/h}I_{r_{p-1}}&0&\hdots&0\\
\vdots&0&\ddots&\ddots&\vdots\\
\vdots&\ddots&\ddots&\ddots&0\\
0&\hdots&\hdots&0&e^{-S_{\nu_1}/h}I_{r_1}
\end{array}
\right)
\ee
where the $r_j=\sharp\widehat \uuu^{(0)}_{\alpha,j}$ are  such that $r_1+\ldots+r_p=\sharp\uuu^{(0)}_\alpha$. Factorizing by $e^{-S_{\nu_p}/h}$, we get 
$\Omega^\alpha=e^{-S_{\nu_p}/h}\wideparen\Omega^\alpha(\tau)$ with 
\be\label{eq:defomegatau}
\wideparen\Omega^\alpha(\tau)=
\left(
\begin{array}{ccccc}
I_{r_p}&0&\hdots&\hdots&0\\
0&\tau_2I_{r_{p-1}}&0&\hdots&0\\
\vdots&0&\ddots&\ddots&\vdots\\
\vdots&\ddots&\ddots&\ddots&0\\
0&\hdots&\hdots&0&\tau_2\tau_3\ldots\tau_pI_{r_1}
\end{array}
\right)
\ee
where $\tau=(\tau_2,\ldots,\tau_p)\in(\R_+^*)^p$ is defined by $\tau_j=e^{(S_{\nu_{p-(j-2)}}-S_{\nu_{p-(j-1)}})/h}$ for any $j=2,\ldots,p$.
With these new notations, one deduces from \eqref{eq:graded1}, that 
$\Lr^{\alpha,*}\Lr^\alpha= h e^{-2S_{\nu_p}/h}\wideparen\mmm^\alpha(\tau)$
with
\be\label{eq:graded2}
\wideparen\mmm^\alpha(\tau)=\wideparen\Omega^\alpha(\tau)(h^{-1}\wideparen\Lr^{\alpha,*}\wideparen\Lr^\alpha)\wideparen\Omega^\alpha(\tau).
\ee
It turns out that matrix such matrices can be described in a slightly more general setting that is useful to compute their spectrum. We introduce this setting now. Throughout, we denote by $\Sr^+(E)$ the set of symmetric positive definite matrix on a vector space $E$.  We will denote by $\Sr^+_{cl}(E)$ the set of $h$-depending matrices $M(h)\in\Sr^+(E)$ admitting a classical expansion 
$M(h)\sim\sum_j h^j M_j$ with $M_0\in\Sr^+(E)$. We will sometimes forget $E$ and write for short $\Sr^+,\Sr^+_{cl}$.
\begin{defin}\label{defin:gradedmat}
Let $\Er=(E_j)_{j=1,\ldots,p}$ be a sequence of finite dimensional vector spaces $E_j$ of dimension $r_j>0$, let 
$E=\oplus_{j=1,\ldots, p}E_j$  and let $\tau=(\tau_2,\ldots,\tau_p)\in (\R_+^*)^{p-1}$. Suppose that  $\tau \mapsto\mmm(\tau)$ is a smooth map from $(\R_+^*)^{p-1}$ into the set of matrices 
$\Mr(E)$. 
\begin{itemize}
\item[-] We say that $\mmm(\tau)$ is an $(\Er,\tau)$-graded matrix if there exists  
$\mmm'\in\Sr^+(E)$ independent of $\tau$ such that $\mmm(\tau)=\Omega(\tau) \mmm'\Omega(\tau)$ with $\Omega(\tau)\in\Mr(E)$ of the form \eqref{eq:defomegatau}, that is $\Omega=\diag(\varepsilon_j(\tau)I_{r_j},\,j=1,\ldots,p)$ where 
$\varepsilon_1(\tau)=1$ and $\varepsilon_j(\tau)=(\prod_{k=2}^j\tau_k) $ for all $j\geq 2$.

\item[-]We say that a family of $(\Er,\tau)$-graded matrices $\mmm_h(\tau)$, $h\in]0,h_0]$ is classical if one has 
$\mmm_h(\tau)=\Omega(\tau) \mmm'(h)\Omega(\tau)$ for some matrix $\mmm'(h)\in\Sr^+_{cl}(E)$.
\end{itemize}
Throughout, we denote by $\Gr(\Er,\tau)$ the set of $(\Er,\tau)$-graded matrices and by $\Gr_{cl}(\Er,\tau)$ the set of classical $(\Er,\tau)$-graded matrices.
\end{defin}
Let us remark that for $p=1$, a graded matrix is simply a $\tau$-independent symmetric positive definite matrix.

\begin{lemma}\label{lem:Malphataugrad} Suppose that $\mmm_h(\tau)$ is a classical $(\Er,\tau)$-graded family of matrices and that $p\geq 2$. Then
one has 
\be\label{eq:taugradmatrix}
\mmm_h(\tau)=
\left(
\begin{array}{cc}
J(h)&\tau_2B_h(\tau')^*\\
\tau_2B_h(\tau')&\tau_2^2\nnn_h(\tau')
\end{array}
\right)
\ee
with 
\begin{itemize}
\item[-] $J(h)\in\Sr^+_{cl}(E_1)$
\item[-] $\nnn_h(\tau')\in\Gr_{cl}(\Er',\tau')$ with $\tau'=(\tau_3,\ldots,\tau_p)$ and $\Er'=(E_j)_{j=2,\ldots,p}$.
\item[-] $B_h(\tau')\in \Mr(E_1,\oplus_{j=2}^p E_j)$ satisfies
$$
B_h(\tau')^*=(b_{2}(h)^*,\tau_3b_{3}(h)^*,\tau_3\tau_4b_{4}(h)^*,\ldots,\tau_3\ldots\tau_pb_{p}(h)^*)
$$
with $b_{j}(h):E_1\rightarrow E_j$ independent of $\tau$ admitting a classical expansion.
\end{itemize}
Moreover, the matrix $\nnn_h(\tau')-B_h(\tau')J(h)^{-1}B_h(\tau')^*$ belongs to $\Gr_{cl}(\Er',\tau')$.
\end{lemma}
\bp
Assume that $\mmm_h(\tau)=\Omega(\tau) \mmm'(h)\Omega(\tau)$ with $\Omega(\tau)$ of the form  \eqref{eq:defomegatau}.
 First observe that
$$
\Omega(\tau)=
\left(
\begin{array}{cc}
I_{r_p}&0\\
0&\tau_2\Omega'(\tau')
\end{array}
\right)
$$
with 
$$
\Omega'(\tau')=
\left(
\begin{array}{ccccc}
I_{r_{p-1}}&0&\hdots&\hdots&0\\
0&\tau_3I_{r_{p-2}}&0&\hdots&0\\
\vdots&0&\ddots&\ddots&\vdots\\
\vdots&\ddots&\ddots&\ddots&0\\
0&\hdots&\hdots&0&\tau_3\tau_4\ldots\tau_pI_{r_1}
\end{array}
\right).
$$
On the other hand, we can write 
$$
\mmm'(h)=
\left(
\begin{array}{cc}
J(h)&B'(h)^*\\
B'(h)&\nnn'(h)
\end{array}
\right)
$$
with  
$J(h),\nnn'(h)\in\Sr^+_{cl}$ and $B'(h)$  admitting a classical expansion.
Therefore,
$$
\Omega(\tau)\mmm'_h\Omega(\tau)
=\left(
\begin{array}{cc}
J(h)&\tau_2 B'(h)^*\Omega'(\tau')\\
\tau_2 \Omega'(\tau')B'(h)&\tau_2^2\Omega'(\tau')\nnn'(h)\Omega'(\tau')
\end{array}
\right)
$$
which has exactly the  form \eqref{eq:taugradmatrix} with $B_h(\tau')= \Omega'(\tau')B'(h)$ and 
$\nnn_h(\tau')=\Omega'(\tau')\nnn'(h)\Omega'(\tau')$.
By construction, $\nnn_h(\tau')$ belongs to $\Gr_{cl}(\Er',\tau')$ and $B_h(\tau')$ has the required form.

It remains to prove that $\rrr_h:=\nnn_h(\tau')-B_h(\tau')J(h)^{-1}B_h(\tau')^*$ belongs to $\Gr_{cl}(\Er',\tau')$. First observe that since $J(h)$ is symmetric positive definite, this quantity is well-defined. Moreover, one has by construction
\begin{equation*}
\begin{split}
\rrr_h&=\Omega'(\tau')\nnn'(h)\Omega'(\tau')- \Omega'(\tau')B'(h) J(h)^{-1} B'(h)^*\Omega'(\tau')\\
&= \Omega'(\tau')\rrr'(h) \Omega'(\tau')
\end{split}
\end{equation*}
with $\rrr'(h)=\nnn'(h)-B'(h) J(h)^{-1} B'(h)^*$.
Since $J(h)\in\Sr^+_{cl}$, then $J(h)^{-1}\in\Sr^+_{cl}$ and $\rrr'(h)$ admits a classical expansion $\rrr'(h)\sim\sum_j h^j\rrr'_j$ with 
$$
\rrr'_0=J_0-B'_0J_0^{-1}(B_0')^*.
$$
Moreover, since $\mmm'(h)\in \Sr^+_{cl}$ then the matrix 
$$
\mmm'_0=
\left(
\begin{array}{cc}
J_0&(B'_0)^*\\
B'_0&\nnn'_0
\end{array}
\right)
$$
is symmetric definite positive. Hence, it follows directly from Lemma \ref{lem:schurpositive} that $\rrr'_0\in\Sr^+$.
\ep
\subsection{The spectrum of graded matrices}\label{subsec:TSOGM}

Using  Lemma \ref{lem:Malphataugrad}, we define an application $\rrr:\Gr_{cl}(\Er,\tau)\rightarrow \Gr_{cl}(\Er',\tau')$ with $\tau'=(\tau_3,\ldots,\tau_p)$ and
$\Er'=\oplus_{j=2}^pE_j$, by
\be\label{eq:defRtau}
\rrr(\mmm_h(\tau))=\nnn_h(\tau')-B_h(\tau')J(h)^{-1}B_h^*(\tau')
\ee
for any $\mmm_h(\tau)\in \Gr_{cl}(\Er,\tau)$. Of course, the map $\rrr$ depends on $\Er$ and $\tau$, but we ommit this dependance since the set on which  $\rrr$ is acting  will be obvious in the sequel. By a slight abuse of notations we will denote 
$\rrr^k=\rrr\circ\ldots\circ\rrr$ ($k$ times). Obviously, $\rrr^k$ acts from $\Gr(\Er,\tau)$ into 
$\Gr(\Er^{(k)},\tau^{(k)})$ with $\Er^{(k)}=\oplus_{j=k+1}^p E_j$ and $\tau^{(k)}=(\tau_{k+2},\ldots,\tau_p)$.
In the same way, we defined $\rrr$, we can define a map $\jjj:\Gr_{cl}(\Er,\tau)\rightarrow \Sr^+_{cl}(E_1)$ by $\jjj(\mmm_h(\tau))=\mmm_h$ if $p=1$ and $\jjj(\mmm_h(\tau))=J(h)$ for any $\mmm_h(\tau)$ having the form \eqref{eq:taugradmatrix} if $p\geq 2$.

\begin{theorem}\label{th:specmatgrad} Let $\Er=(E_j)_{j=1,\ldots,p}$ be a finite sequence of  vector space $E_j$ of finite dimension $n_j=\dim E_j$
and let $\tau=(\tau_2,\ldots,\tau_p)\in (\R_+^*)^{p-1}$. Suppose that $\mmm_h(\tau)$ is classical $(\Er,\tau)$-graded. There exists $h_0>0$ and $\delta>0$ such that 
uniformly with respect to  $h\in]0,h_0]$ and $\vert\tau\vert_\infty<\delta$, one has
\be\label{eq:locspecgrad0}
\sigma(\mmm_h(\tau))=\bigsqcup_{j=1}^p\varepsilon_j\sigma(\jjj\circ\rrr^{j-1}(\mmm_h(\tau)))(1+\ooo(\vert\tau\vert_\infty^2))
\ee
with  $\varepsilon_j=\varepsilon_j(\tau)$ given in Definition \ref{defin:gradedmat}.
\end{theorem}
\begin{remark}
In the above theorem, the matrix $\jjj\circ\rrr^{j-1}(\mmm_h(\tau))$ is always independent of the parameter $\tau$. Let us denote
$\{\lambda_1^j\leq\ldots\leq \lambda_{n_j}^j\}=\sigma(\jjj\circ\rrr^{j-1}(\mmm_h(\tau))$.
The identity \eqref{eq:locspecgrad0} means that there exists $a,b>0$ independent of $\tau,h$ such that 
$$
\sigma(\mmm_h(\tau))\subset\bigsqcup_{j=1}^p\varepsilon_j[a,b]
$$
and that for all $j=1,\ldots,p$, $\mmm_h(\tau)$ has exactly $n_j$ eigenvalues $\mu_1^j\leq\ldots\leq \mu_{n_j}^j$ 
 in $\varepsilon_j[a,b]$ and 
$$
\mu_n^j=\varepsilon_j(\lambda_n^j+\ooo(\vert\tau\vert_\infty^2)).
$$
\end{remark}
\bp
We prove the theorem by induction on $p$. Throughout the proof the notation $\ooo(.)$ is uniform with respect to the parameters $h$ and $\tau$.
For $p=1$, $\mmm_h(\tau)=\mmm_h\in\Sr^+_{cl}(E_1)$ is independent of $\tau$ and $\jjj \rrr^{0}(\mmm_h(\tau))=\jjj\mmm_h(\tau)=\mmm_h$
which proves the statement.

Suppose now that $p\geq 2$ and 
let $\mmm_h(\tau)\in\Gr_{cl}(\Er,\tau)$. We have
$$
\mmm_h(\tau)=\left(
\begin{array}{cc}
J(h)&\tau_2B_h(\tau')^*\\
\tau_2B_h(\tau')&\tau_2^2\nnn_h(\tau')
\end{array}
\right)
$$
with $J(h), B_h(\tau')$ and $\nnn_h(\tau')$ as in Lemma \ref{lem:Malphataugrad}.
In order to lighten the notation we will drop the variable $\tau, \tau'$ in the proof below.
For $\lambda\in\C$, let
\begin{equation}\label{eq:schur1}
\ppp(\lambda):=\mmm_h-\lambda
=
\left(
\begin{array}{cc}
J(h)-\lambda &\tau_2 B_h^*\\
\tau_2 B_h&\tau_2^2\nnn_h-\lambda
\end{array}
\right).
\end{equation}  
This is  an holomorphic function, and since it is non trivial, its inverse is well defined 
excepted for a finite number of values of $\lambda$ which are exactly the  eigenvalues of $\mmm_h$. Moreover $\lambda\in\C\mapsto \ppp(\lambda)^{-1}$ is meromorphic with poles in $\sigma(\mmm_h)$ and 
for any $\mu$ in $\sigma(\mmm_h)$, the rank of the residue of $\ppp(\lambda)^{-1}$ at $\mu$ is exactly the multiplicity of $\mu$ as an eigenvalue.

Let us first prove that $\mmm_h$ admits at least $n_1$ eigenvalues of size $1$.
Let $\lambda_n^1=\lambda_n^1(h)$, $n=1,\ldots,n_1$ denote the increasing sequence of eigenvalues of the positive definite matrix $J(h)$. 
Since $J(h)=J_0+\ooo(h)$ with $J_0\in \Sr^+$, then the $\lambda_n^1(h)$  satisfy
$\lambda_n^1(h)=\lambda_{n,0}^1+\ooo(h)$ with $\lambda^1_{n,0}$ eigenvalue of $J_0$. In particular $\lambda^1_{n,0}>0$ for all $n=1,\ldots,n_1$ and hence there exists $c_1,d_1>0$ and $h_0>0$ such that for $h\in]0,h_0]$ and all $n=1,\ldots,n_1$, one has
$\lambda_n^1(h)\in [c_1,d_1]$.
Let $n\in\{1,\ldots,n_1\}$ be fixed and consider 
$D_n=D_n(h,\tau_2)=\{z\in\C,\;\vert z-\lambda^1_n\vert\leq M\tau_2^2\}$ for some $M>0$ that will be chosen large enough later and 
$\tilde D_n=\{z\in\C,\;\vert z-\lambda^1_n\vert\leq 2M\tau_2^2\}$.
Observe that for $h,\tau_2>0$ small enough, the disks $\tilde D_n$ are disjoint. By definition, one has $\nnn_h(\tau')=\ooo(1)$ and since
 $\lambda_n^\ell\geq c_1>0$, this implies that for $\tau_2>0$ small enough with respect to $c$ and  $\lambda\in \tilde D_n$, the matrix  $\tau_2^2\nnn_h(\tau')-\lambda$ is  invertible, and 
$(\tau_2^2\nnn_h(\tau')-\lambda)^{-1}=\ooo(1)$. 
Moreover, for $\lambda\in\tilde D_n\setminus D_n$, $J(h)-\lambda$ is invertible and $(J(h)-\lambda)^{-1}=\ooo(\tau_2^{-2}M^{-1})$. This implies that for 
$M>0$ large enough, $J(h)-\lambda-\tau_2^2 B_h^*(\tau_2^2\nnn_h-\lambda)^{-1}B_h$ is invertible with 
\be\label{eq:invertshur1}
\begin{split}
\big(J(h)-\lambda-\tau_2^2 B_h^*(\tau_2^2\nnn_h-\lambda)^{-1}B_h\big)^{-1}&=(J(h)-\lambda)^{-1}\Big(I-\tau_2^2 B_h^*(\tau_2^2\nnn_h-\lambda)^{-1}B_h(J(h)-\lambda)^{-1}\Big)^{-1}\\
&=(J(h)-\lambda)^{-1}(1+\ooo(M^{-1})).
\end{split}
\ee
Hence, the standard Schur complement procedure shows that for $\lambda\in\tilde D_n\setminus D_n$, 
$\ppp(\lambda)$ is invertible with  inverse $\eee(\lambda)$  given by
\begin{equation}\label{eq:schur2}
\eee(\lambda)=
\left(
\begin{array}{cc}
E(\lambda)&-\tau_2 E(\lambda)B_h^*(\tau_2^2\nnn_h-\lambda)^{-1}\\
-\tau_2(\tau_2^2\nnn_h-\lambda)^{-1}B_hE(\lambda)&E_0(\lambda)
\end{array}
\right)
\end{equation}
with 
\begin{equation*}
E(\lambda)=\Big(J(h)-\lambda-\tau_2^2 B_h^*(\tau_2^2\nnn_h-\lambda)^{-1}B_h\Big)^{-1}
\end{equation*}
and 
\begin{equation*}
E_0(\lambda)=(\tau_2^2\nnn_h-\lambda)^{-1}+\tau_2^2(\tau_2^2 \nnn_h-\lambda)^{-1}B_hE(\lambda)B_h^*(\tau_2^2\nnn_h-\lambda)^{-1}
\end{equation*}
 By the functional calculus and Cauchy formula, the number of eigenvalues of $\mmm_h$ (counted with multiplicity) in 
$D_n$ is equal to  the rank of the  projector
$$\Pi_n=\frac 1{2i\pi}\int_{\partial D_n}\eee(\lambda)d\lambda.$$
Let us denote $R_{\tau_2}(\lambda)=-\tau_2 (\tau_2^2\nnn_h-\lambda)^{-1}B_h$ and  $R^\dagger_{\tau_2}(\lambda)=-\tau_2 B_h^*(\tau_2^2\nnn_h-\lambda)^{-1}$ , then
\begin{equation}\label{eq:schur40}
\Pi_n=\frac 1{2i\pi}\int_{\partial D_n}
\left(
\begin{array}{cc}
E(\lambda)&E(\lambda)R^\dagger_{\tau_2}(\lambda)\\
R_{\tau_2}(\lambda) E(\lambda)&(\tau_2^2 \nnn_h-\lambda)^{-1}+R_{\tau_2}(\lambda)E(\lambda)R^\dagger_{\tau_2}(\lambda)
\end{array}
\right)d\lambda
\end{equation}
Since $\nnn_h=\ooo(1)$, then $\tau_2^2\nnn_h-\lambda$ is invertible in  $\tilde D_n$  and it follows that
$$
\Pi_n=\frac 1{2i\pi}\int_{\partial D_n}
\left(
\begin{array}{cc}
E(\lambda)&E(\lambda)R^\dagger_{\tau_2}(\lambda)\\
R_{\tau_2}(\lambda) E(\lambda)&R_{\tau_2}(\lambda)E(\lambda)R^\dagger_{\tau_2}(\lambda)
\end{array}
\right)d\lambda
$$
that can be written
$$
\Pi_n=\frac 1{2i\pi}\int_{\partial D_n}
\rrr_{\tau_2}(\lambda)\left(
\begin{array}{cc}
E(\lambda)&E(\lambda)\\
E(\lambda)&E(\lambda)
\end{array}
\right)\rrr^\dagger_{\tau_2}(\lambda) d\lambda
$$
with 
 $\rrr_{\tau_2}^\bullet=\left(
\begin{array}{cc}
I&0\\
0&R^\bullet_{\tau_2}
\end{array}
\right)
$. Moreover,  since $R_{\tau_2}^\bullet$ is invertible and holomorphic in   $\tilde D_n$ then so is $\rrr_{\tau_2}^\bullet$. Therefore, 
$rk(\Pi_n)=rk(\tilde \Pi_n)$ where
$$
\tilde \Pi_n=\frac 1{2i\pi}\int_{\partial D_n}
\left(
\begin{array}{cc}
E(\lambda)&E(\lambda)\\
E(\lambda)&E(\lambda)
\end{array}
\right) d\lambda=\left(\begin{array}{cc}
E_n&E_n\\
E_n&E_n
\end{array}
\right) 
$$
with 
$$
E_n=\frac 1{2i\pi}\int_{\partial D_n}\Big(J(h)-\lambda-\tau_2^2 B_h^*(\tau_2^2\nnn_h-\lambda)^{-1}B_h\Big)^{-1}d\lambda.
$$
Since 
for $M$ large enough independent of $(h,\tau)$, the matrix   $(I-\tau_2^2 B_h^*(\tau_2^2\nnn_h-\lambda)^{-1}B_h(J(h)-\lambda)^{-1})^{-1}$ is holomorphic in $ \tilde D_n$, it follows from \eqref{eq:invertshur1} that 
the rank of $E_n$ is exactly the multiplicity of $\lambda_n^1$ and hence the rank of $\Pi_n$ is exactly the multiplicity of $\lambda_n^1$. 
This proves that $\mmm_h$ admits at least $n_1$ eigenvalues $\mu_1^1\leq\ldots\leq \mu_{n_1}^1$ in the interval
$[c_1-M\tau_2^2,d_1+M\tau_2^2]$ and that these eigenvalues satisfy
\be\label{eq:locspec11}
\mu_n^1=\lambda_n^1+\ooo(\tau_2^2),\;\;\forall n=1,\ldots,n_1.
\ee

Let us now study the eigenvalues below $\tau_2^2$. Throughout the proof, we denote $t=\vert\tau'\vert_\infty$.
Thanks to the last part of Lemma \ref{lem:Malphataugrad}, the matrix $\zzz_h(\tau'):=\rrr(\mmm_h(\tau))=\nnn_h-B_hJ(h)^{-1}B_h^*$ is classical $(\Er',\tau')$-graded. 
Hence, it follows from the induction hypothesis that uniformly with respect to $h$, one has
\be\label{eq:locspec2}
\sigma(\zzz_h(\tau'))=\bigsqcup_{j=2}^p\tilde \varepsilon_j\sigma(\jjj\circ\rrr^{j-2}(\zzz_h(\tau')))(1+\ooo(\vert\tau'\vert_\infty^2))
\ee
with $\tilde \varepsilon_j=(\prod_{l=3}^{j}\tau_l)^2$ for $j\geq 3$ and $\tilde \varepsilon_2=1$.
Moreover, by definition, one has $\zzz_h=\rrr(\mmm_h(\tau))$, hence \eqref{eq:locspec2} rewrites
\be\label{eq:locspec3}
\sigma(\zzz_h(\tau'))=\bigsqcup_{j=2}^p\tilde \varepsilon_j\sigma(\jjj\circ\rrr^{j-1}(\mmm_h(\tau)))(1+\ooo(\vert\tau'\vert_\infty^2)).
\ee
Since $\mmm_h(\tau')\in\Gr_{cl}(\Er,\tau)$, then for all $j=2,\ldots,p$ the matrix $\jjj\circ\rrr^{j-1}(\mmm_h(\tau))$ belongs to $\Sr^+_{cl}(E_j)$.
For $j=2,\ldots,p$, let  $\lambda_1^j(h)\leq \ldots\leq \lambda_{n_j}^j(h)$ denote the eigenvalues of the symmetric matrix 
$\jjj\circ\rrr^{j-1}(\mmm_h(\tau))$. 
As above, this implies that there exists $c_j,d_j>0$ and $h_0>0$ such that for all $h\in]0,h_0]$ the 
eigenvalues $\lambda_{n_j}^j(h)$ satisfy  $\lambda_{n_j}^j(h)\in [c_j,d_j]$ for all $n=1,\ldots,n_j$.
Suppose now that $j\in\{2,\ldots,p\}$ and $n\in\{1,\ldots,n_j\}$ are fixed and consider 
$D'_{j,n}=\{z\in\C,\;\vert z-\varepsilon_j\lambda^j_n\vert\leq Mt^2\varepsilon_j\}$ for some $M>0$ to be chosen large enough and 
$\tilde D'_{j,n}=\{z\in\C,\;\vert z-\varepsilon_j\lambda^j_n\vert\leq 2Mt^2\varepsilon_j\}$. As above, we introduce also 
the corresponding projector
$$
\Pi'_{j,n}=\frac 1{2i\pi}\int_{\partial D'_{j,n}}\eee(\lambda)d\lambda.
$$
Since $J_0$ is invertible, we know that for $\lambda$ in   $\tilde D'_{j,n}$ and $h,t$ small enough, $J(h)-\lambda$ is invertible and once again the Schur complement formula permits to 
write the invert of $\ppp(\lambda)$
\begin{equation}\label{eq:schur5}
\eee(\lambda)=
\left(
\begin{array}{cc}
E_0(\lambda)&-\tau_2 (J(h)-\lambda)^{-1}B_h^*E(\lambda)\\
-\tau_2 E(\lambda)B_h(J(h)-\lambda)^{-1}&E(\lambda)
\end{array}
\right)
\end{equation}
with 
$$
E(\lambda)=\Big(\tau_2^2 \nnn_h-\lambda-\tau_2^2 B_h(J(h)-\lambda)^{-1}B_h^*\Big)^{-1}
$$
and 
$$
E_0(\lambda)=(J(h)-\lambda)^{-1}+\tau_2^2(J(h)-\lambda)^{-1}B_h^*E(\lambda)B_h(J(h)-\lambda)^{-1}.
$$
Setting $\lambda=\tau_2^2 z$, we get (using the relation $\varepsilon_j=\tau_2^2\tilde \varepsilon_j$)
$$
\Pi'_{j,n}=\frac {\tau_2^2}{2i\pi}\int_{\partial \hat D'_{j,n}}\eee(\tau_2^2 z)dz
$$
with $\hat D'_n=\{z\in\C,\;\vert z-\tilde\varepsilon_j\lambda^j_n\vert\leq Mt^2\tilde \varepsilon_j\}$.
Moreover, for  $\vert z-\tilde\varepsilon_j\lambda^j_n\vert= Mt^2\tilde \varepsilon_j$, the matrix $J(h)$ is invertible with $J(h)^{-1}=\ooo(1)$, hence we  have 
\begin{equation*}
\begin{split}
E(\tau_2^2z)&=\tau_2^{-2}\big(\nnn_h-z- B_h(J(h)-\tau_2^2 z)^{-1}B_h^*\big)^{-1}\\
&=\tau_2^{-2}\big(\zzz_h-z+\ooo(\tau_2^2\vert z\vert))^{-1}\\
&=\tau_2^{-2}(\zzz_h-z)^{-1}(I+\ooo(\tau_2^2\tilde \varepsilon_j\Vert(\zzz_h-z)^{-1}\Vert)).
\end{split}
\end{equation*}
Moreover, by definition of $\hat D'_{j,n}$ and thanks to \eqref{eq:locspec2}, one has $\dist(z,\sigma(\zzz_h))\geq  \frac 12Mt^2\tilde \varepsilon_j$ 
for any $z\in \partial \hat D'_{j,n}$.
Hence $\Vert (\zzz_h-z)^{-1}\Vert\leq 2(Mt^2\tilde \varepsilon_j)^{-1}$ and since $t\geq \tau_2$, it follows that 
$$
E(\tau_2^2z)
=\tau_2^{-2}(\zzz_h-z)^{-1}(I+\ooo(M^{-1})).
$$
Integrating along $\partial \tilde D'_{j,n}$ and working as above, we get 
$$
\Pi'_{j,n}=\frac 1{2i\pi}\int_{\partial D'_{j,n}}
\rrr_{\tau_2}^\dagger(\lambda)\left(
\begin{array}{cc}
E(\lambda)&E(\lambda)\\
E(\lambda)&E(\lambda)
\end{array}
\right)\rrr_{\tau_2}(\lambda) d\lambda
$$
with 
 $\rrr_{\tau_2}^\bullet=\left(
\begin{array}{cc}
R^\bullet_{\tau_2}&0\\
0&I
\end{array}
\right)
$.
The same argument as above show that $rk(\Pi_n)=rk(E'_n)$ 
with
\begin{equation*}
\begin{split}
E'_n=\frac {\tau_2^2}{2i\pi}\int_{\partial \hat D'_{j,n}}E(\tau_2^2 z)dz
=\frac 1{2i\pi}\int_{\partial \hat D'_{j,n}}\Big(  \zzz_h-z )^{-1}(I+\ooo(M^{-1}))^{-1}dz
\end{split}
\end{equation*}
This shows again that the rank of $E'_n$ (and hence $\Pi'_{j,n}$) is exactly the multiplicity of $\lambda_n^j$. 
Therefore, for any $j=2,\ldots,p$,  $\mmm_h$ admits 
$n_j$ eigenvalues $\mu_1^1\leq\ldots\mu_{n_1}^1$ in the interval $\varepsilon_j[c_j-Mt^2,d_j+Mt^2]$ and that these eigenvalues satisfy
\be
\mu_n^j=\varepsilon_j(\lambda_n^j+\ooo(\vert\tau\vert_\infty^2)),\;\;\forall n=1,\ldots,n_j.
\ee
Combining this estimate with \eqref{eq:locspec11} and using the fact the $\dim(E)=\sum_{j=1}^p r_j$, we obtain the $\mu_n^j$ are the only eigenvalues 
of $\mmm_h$. This completes the proof. 
 \ep
 
 \subsection{The singular values of $\Lr^\alpha$}\label{subsec:TSVOL}
 Given $\m,\m_1,\m_2\in\uuu^{(0)}$ and $\s\in\vvv^{(1)}$, we denote
\begin{equation}\label{eq:definupsilon2}
\begin{split}
\upsilon_2(\s,\m,\m_1,\m_2)=
  \pi^{-\frac 12}\vert\hat\lambda_1(\s)\vert^{\frac 12}\Big(
\frac{h_\phi(\m_1)}{h_\phi(\s)}\delta_{\m,\m_1}
-\frac{h_\phi(\m_2)}{h_\phi(\s)}\delta_{\m,\m_2}\Big)
\end{split}
\end{equation}
and 
\begin{equation}\label{eq:definupsilon1}
\upsilon_1(\s,\m,\m_1)=
 \pi^{-\frac 12}\vert\hat\lambda_1(\s)\vert^{\frac 12}
\frac{h_\phi(\m_1)}{h_\phi(\s)}\delta_{\m,\m_1}.
\end{equation}
 Let us define the matrix $\Upsilon^\alpha\in\Mr(\widehat\uuu^{(0)}_\alpha,\vvv^{(1)}_\alpha)$ by 
\begin{equation}\label{eq:definmatrixUpsilon}
\Upsilon^\alpha(\s,\m)=
\left\{
\begin{array}{c}
  \upsilon_2(\s,\m,\m_1(\s),\m_2(\s))\;\;\text{ if } \m_2(\s)\in \widehat\uuu^{(0)}_\alpha\\
\upsilon_1(\s,\m,\m_1(\s))\;\; \text{ if } \m_2(\s)\notin \widehat\uuu^{(0)}_\alpha\phantom{****}
\end{array}
\right.
\end{equation}
where the indexes $\m,\s$ are enumerated according to the partitions of Section \ref{subsec:GSOTM}.
Observe that with this notation, the conclusion of Lemma \ref{lem:computinter} rephrases as $\widetilde\Lr^{\alpha,0}=\Upsilon$. Moreover, the above expression can be simplified according to the type of $\alpha$. More precisely, 
\begin{itemize}
\item[-] if $\alpha$ is of type I, then $\m_2(\s)\in \widehat\uuu^{(0)}_\alpha$ if and only if $\s\in\vvv_\alpha^{(1),i}$ 
\item[-] if $\alpha$ is of type II, then $\m_2(\s)$ is always in $\widehat\uuu^{(0)}_\alpha$.
\end{itemize}

\begin{theorem}\label{thm:SVLalpha} Let $\mmm^{\alpha}=\Lr^{\alpha,*}\Lr^\alpha$.
There exist $c>0$ such that counted with multiplicity, one has
$$
\sigma(\mmm^\alpha)=\bigsqcup_{j=1}^{p(\alpha)}he^{-2h^{-1}S_{\nu_j^\alpha}}\sigma(M^{\alpha,j})(1+\ooo(e^{-c/h}))
$$
where the matrices $M^{\alpha,j}$ have a classical expansion $M^{\alpha,j}\sim\sum h^kM^{\alpha,j}_k$ whose leading term is given by
$$
M^{\alpha,j}_0=\jjj\rrr^{j-1}(\zzz^\alpha)
$$
where $\zzz^\alpha= \Omega^\alpha\Tr^{\alpha,0}\Upsilon^{\alpha,*}\Upsilon^\alpha\Tr^{\alpha,0}\Omega^\alpha$ belongs to $\Gr(\Er,\tau)$ with 
$\Er=(\Fr( \widehat\uuu^{(0)}_{\alpha,j}))_{j=1,\ldots,p}$ and $\tau=(\tau_j)_{j=1,\ldots,p}$ with $\tau_j=e^{(S_{\nu_{p-(j-2)}}-S_{\nu_{p-(j-1)}})/h}$
\end{theorem}
\bp
One has 
$$
\mmm^\alpha=\Lr^{\alpha,*}\Lr^\alpha=he^{-2S_{p_1}/h}\wideparen\mmm^\alpha
$$
with $\wideparen\mmm^\alpha$ given by  \eqref{eq:graded2}:
$$
\wideparen\mmm^\alpha(\tau)=\wideparen\Omega^\alpha(\tau)^*\wideparen\mmm^{\alpha,'}\wideparen\Omega^\alpha(\tau)
$$
with
$\wideparen\mmm^{\alpha,'}=(h^{-1}\wideparen\Lr^{\alpha,*}\wideparen\Lr^\alpha)$. Of course, this matrix is symmetric positive and thanks to Proposition \ref{prop:injectLalpha}, it admits a classical expansion
$$\wideparen\mmm^{\alpha,'}\sim\sum_k h^k \wideparen\mmm^{\alpha,'}_k$$
with $\wideparen\mmm^{\alpha,'}_0=(\wideparen\Lr^{\alpha,0})^*\wideparen\Lr^{\alpha,0}=\Tr^{\alpha,0}\Upsilon^{\alpha,*}\Upsilon^\alpha\Tr^{\alpha,0}\in\Sr^+$.
This shows  that $\wideparen\mmm^{\alpha,'}$ belongs to $\Sr^+_{cl}$. Hence
$\wideparen\mmm^\alpha$ is classical $(\Er,\tau)$-graded with $\Er=(\Fr( \widehat\uuu^{(0)}_{\alpha,j}))_{j=1,\ldots,p}$ and 
$\tau=(\tau_2,\ldots,\tau_p)$, $\tau_j=e^{(S_{\nu_{p-(j-2)}}-S_{\nu_{p-(j-1)}})/h}$ and the conclusion follows directly from Theorem \ref{th:specmatgrad}.
\ep
 
 \section{ Proof of main theorem}\label{sec:POMT}
In this section we explain how one can deduce Theorem \ref{th:main} from Theorem \ref{thm:SVLalpha}.
As in \cite{HeKlNi04_01}, the general idea is to compare the singular values of the successive reduced matrix by mean of Fan inequalities.
As a preparation, we shall compare the matrices  $\lll^{\pi,'}$ and $\lll^{bkw,'}$ defined in Section 3. 
First, observe that thanks to \eqref{eq:uselesssp}, \eqref{eq:definLr''}, \eqref{eq:factorLalpha}, one has
\be\label{eq:factorLbkw}
\lll^{bkw,'}=\Jr\lll^{bkw,''}=\Jr\Lr=\Jr\widetilde\Lr\Tr\Omega
\ee
with 
 $\Jr:\Fr(\vvv^{(1)})\rightarrow \Fr(\uuu^{(1)})$ defined by $\Jr_{\s,\s'}=\delta_{\s,\s'}$,  $\widetilde \Lr=\diag(\widetilde\Lr^\alpha,\alpha\in\ala)$, 
 $\Tr=\diag(\Tr^\alpha,\alpha\in\ala)$ and $\Omega=\diag(\Omega^\alpha,\alpha\in\ala)$.
\begin{lemma}\label{lem:comparLpi-Lbkw}
There exists $\gamma>0$ such that 
$$
\lll^{\pi,'}=(\Jr+\ooo(e^{-\gamma/h}))\Lr.
$$
\end{lemma}
\bp
First observe that thanks to  Lemma \ref{lem:approxSV}, one has
\be\label{eq:comparLpi-Lbkw1}
\lll^{\pi,'}=\lll^{bkw,'}+\rrr
\ee
with $\rrr:\Fr(\ulu^{(0)})\rightarrow\Fr(\uuu^{(1)})$
satisfying 
\be\label{eq:comparLpi-Lbkw2}
\rrr_{\s,\m}=\ooo(e^{-(S(\m)+\gamma)/h}),\;\forall\m\in\ulu^{(0)}
\ee
 for some $\gamma>0$.
 Using \eqref{eq:factorLbkw}, we get 
 $$
 \lll^{\pi,'}=\Jr\widetilde\Lr\Tr\Omega+\widetilde\rrr\Omega
 $$
 with $\widetilde\rrr=\ooo(e^{-\gamma/h})$.
Hence, we have to prove that there exists 
$\overline\rrr:\Fr(\vvv^{(1)})\rightarrow\Fr(\uuu^{(1)})$ such that 
$
\widetilde\rrr=\overline\rrr \widetilde\Lr\Tr
$
and 
$\overline\rrr=\ooo(e^{-\gamma/h})$.
From Proposition \ref{prop:injectLalpha}, we know that the matrix $\Wr:=(\widetilde\Lr\Tr)^*\widetilde\Lr\Tr$ 
is invertible with inverse uniformly bounded with respect to $h$. This allows to define
$\overline\rrr:=\widetilde\rrr\Wr^{-1}(\widetilde\Lr\Tr)^*$. Thanks to the above remarks, one has $\overline\rrr=\ooo(e^{-\gamma/h})$ and by construction
$$
\overline\rrr \widetilde\Lr\Tr=\widetilde\rrr\Wr^{-1}(\widetilde\Lr\Tr)^*\widetilde\Lr\Tr=\widetilde\rrr
$$
which completes the proof.
\ep

We are now ready to prove Theorem \ref{th:main}. Until the end of this section, $\gamma,C>0$ denote some 
constants independent on $h$ that may change from line to line. We shall also denote by $SV(M)$ the singular values of any matrix $M$.

From Section \ref{subsec:GSOP}, we know that the $n_0$ exponentially small eigenvalues of $\Delta_\phi^{(0)}$ are the square of the singular values of the matrix $\lll$. Thanks to Lemmas \ref{lem:projectQM0} and \ref{lem:projectQM1}, one has 
$$
\lll=(\Id+\ooo(e^{-\gamma/h}))\lll^\pi(\Id+\ooo(e^{-\gamma/h}))
$$
and it follows from Fan inequality (Lemma \ref{lem:Fan}) that 
$$
SV(\lll)= SV(\lll^\pi)(1+\ooo(e^{-\gamma/h})).
$$
Hence, we are reduced to compute the singular values of $\lll^\pi$. 
Since the first column of $\lll^{\pi}$ is the null vector, it follows that the non zero singular values of 
$\lll^\pi$ are the singular values of $\lll^{\pi,'}$. 
From Lemma \ref{lem:comparLpi-Lbkw}, one knows that 
\be\label{eq:finalproof1}
\lll^{\pi,'}=(\Jr+\ooo(e^{-\gamma/h}))\Lr.
\ee
and since $\Jr^*\Jr=\Id$ this implies for $h$ small enough
\be\label{eq:finalproof2}
\Lr=(\Jr^*+\ooo(e^{-\gamma/h}))\lll^{\pi,'}.
\ee
Using the fact that $\Vert \Jr\Vert=\Vert\Jr^*\Vert=1$,  \eqref{eq:finalproof1} and  \eqref{eq:finalproof2} combined with Lemma \ref{lem:Fan} show that 
$$
SV(\lll^{\pi,'})=(1+\ooo(e^{-\gamma/h}))SV(\Lr).
$$
Combined with Theorem \ref{thm:SVLalpha} this proves Theorem \ref{th:main}.
 \section{Some particular cases and examples}
 In this section, we rephrase Theorem \ref{thm:SVLalpha} in the particular situations $p(\alpha)=1$ and  $p(\alpha)=2$.
 \subsection{The case $p(\alpha)=1$}\label{sec:COSVF1}
In this section we assume that  $p(\alpha)=1$. Then, the set $\sss_\alpha$ is reduced to a singleton $\sss_\alpha=\{S_{\nu_1^\alpha}\}$.
Moreover, the points of $\uuu^{(0)}_\alpha$ are either all of type I, or all of type II.  

\subsubsection{The case where $\alpha$ is of type II}
We first assume that $\alpha$ is of type II. Then all the points $\m\in\uuu^{(0)}_\alpha$ are of type II and Theorem \ref{thm:SVLalpha} takes the following form
\begin{theorem}\label{thm:spectreDtypeII}
Let $\alpha\in\aaa$ be such that $p(\alpha)=1$ and all the points of $\uuu^{(0)}_\alpha$ are of type II. Then the matrix $\Lr^{\alpha}$ has exactly
$q_{\alpha}=\sharp\uuu_{\alpha}^{(0)}$ singular values counted with multiplicity, 
$\rho_{\alpha,\mu}(h)$, $\mu=1,\ldots,q_{\alpha}$. They have the following form
$$
\rho_{\alpha,\mu}(h)=h^{\frac 12} \zeta_{\alpha,\mu}(h)e^{-S_{\nu_1^\alpha}/h}
$$
where  $\zeta_{\alpha,\mu}\sim \sum_{r=0}^\infty h^{r} \zeta_{\alpha,\mu,r}$ is a classical symbol 
such that the $\zeta_{\alpha,\mu,0}$, $\mu=1,\ldots,q_{\alpha}$ are the non zero singular values of the matrix 
$\Upsilon^\alpha\in\Mr(\widehat\uuu^{(0)}_\alpha,\vvv^{(1)}_\alpha)$ given by
$$
\Upsilon^\alpha_{\s,\m}= \pi^{-\frac 12}\vert\hat\lambda_1(\s)\vert^{\frac 12}\Big(
\frac{h_\phi(\m_1(\s))}{h_\phi(\s)}\delta_{\m,\m_1(\s)}
-\frac{h_\phi(\m_2(\s))}{h_\phi(\s)}\delta_{\m,\m_2(\s)}\Big),\;\;\forall \s\in\vvv^{(1)}_\alpha,\;\m\in\widehat\uuu^{(0)}_\alpha
$$
with $\m_1(\s),\m_2(\s)$ defined in Lemma \ref{lem:defm1m2}.
\end{theorem}

Observe that the description of the approximated small eigenvalues of $\Delta_\phi$ in the above theorem is very close in spirit to that obtained in non degenerate situations. Though, the different eigenvalues $\rho_{\alpha,\mu}$ are linked one each other, the only minima  involved
in the computation of the prefactor $\zeta_{\alpha,\mu}$ are associated to the typical height $S_{\nu_1^\alpha}$. In that sense, we can say that the above formula is a generalized Eyring-Kramers formula.\\

As already mentioned in the introduction, the matrix $\Upsilon^\alpha$ enjoys a nice interpretation in terms of graphs theory. In order to simplify, suppose that the function $\phi$ is such that the coefficients of $\Upsilon^\alpha$ are either $1$ or $-1$. Define a graph $\ggg_\alpha$ associated to the 
equivalence class $\alpha$ in the following way. The vertices of the graph are the minima $\m\in\widehat\uuu^{(0)}_\alpha$ and the edges are the 
saddle points $\s\in\vvv^{(1)}_\alpha$. The two vertices associated to the edge  $\s\in\vvv^{(1)}_\alpha$ are just $\m_1(\s)$ and $\m_2(\s)$.
With this definition it turns out that the matrix $\Upsilon^\alpha$ is the transpose of the incidence matrix of  a certain oriented version of the graph 
$\ggg$. As a consequence, the $\vert \zeta_{\alpha,\mu,0}\vert^2$ are the eigenvalues of the corresponding graph Laplacian 
 $\Delta_\ggg=(\delta_{\m,\m'})_{\m,\m'\in\widehat\uuu^{(0)}}$ defined by
\begin{equation}\label{eq:deflaplacegraphe}
\delta_{\m,\m'}=
\left\{
\begin{array}{c}
\operatorname{d}(\m)\text{ if }\m=\m'\phantom{++++++++++++++++++++}\\ 
-1\text{ if } \m\neq \m'\text{ and there is an edge between }\m \text{ and }\m'\\
0\text{ otherwise}\phantom{++++++++++++++++++++++}
\end{array}
\right.
\end{equation}
where the degree $\operatorname{d}(\m)$  is the number of edges incident to the vertex $\m$.

Figure \ref{figsublevel} in the introduction presents an example of potential $\phi$ having one unique saddle value $\sigma$ and such that  all local minima are absolute minima. We represent also in Figure \ref{figsublevel} the graph associated to the non trivial equivalence class (that is the one which is not reduced to one element).

In the case where the coefficients of $\Upsilon^\alpha$ are not necessarily equal to $\pm1$, the same interpretation is available with weighted graphs.
We refer to  \cite{CvDoSa95} for definitions and standard  results on graphs theory.

%

\subsubsection{ The case where $\alpha$ is of type I}\label{subsec:SVOI}
In this section, we compute explicitly the singular values of $\Lr^{\alpha}$, when $\alpha$ is of type I.
\begin{theorem}\label{thm:spectreDtypeI}
Let $\alpha\in\aaa$ be such that $p(\alpha)=1$ and all the points of $\uuu^{(0)}_\alpha$ are of type I. Then, the matrix $\Lr^{\alpha}$ has exactly
$q_{\alpha}:=\sharp\uuu_{\alpha}^{(0)}$ singular values counted with multiplicity. These  singular values 
$\rho_{\alpha,\mu}(h)$, $\mu=1,\ldots,q_{\alpha}$ have the following form
$$
\rho_{\alpha,\mu}(h)=\zeta_{\alpha,\mu}(h)e^{-S_{\nu_1^\alpha}/h}
$$
where  $\zeta_{\alpha,\mu}\sim h^{\frac 12}\sum_{r=0}^\infty h^{r} \zeta_{\alpha,\mu,r}$ has a classical expansion
such that $\zeta_{\alpha,\mu,0}$ are the $q_{\alpha}$ singular values of the matrix $\Upsilon^{\alpha}$
 given by
$$
\Upsilon_{\s,\m}= \pi^{-\frac 12}\vert\hat\lambda_1(\s)\vert^{\frac 12}\Big(
\frac{h_\phi(\m_1(\s))}{h_\phi(\s)}\delta_{\m,\m_1(\s)}
-\frac{h_\phi(\m_2(\s))}{h_\phi(\s)}\delta_{\m,\m_2(\s)}\Big)$$
if $\s\in\vvv^{(1),\i}_{\alpha}$ and 
$$
\Upsilon_{\s,\m}=
\pi^{-\frac 12}\vert\hat\lambda_1(\s)\vert^{\frac 12}
\frac{h_\phi(\m_1(\s))}{h_\phi(\s)}\delta_{\m,\m_1(\s)}
$$
if $\s\in\vvv^{(1),\b}_{\alpha}$. Moreover, these singular values are non zero.
\end{theorem}

As in the case of points of type II  we can interpret the matrix $\widetilde \Lr^{\alpha,0}$ in terms of graphs.
However, some saddle points are now associated to only one minimum. In terms of graph, this 
leads to some  edges having only one vertex which means that
we are  dealing with hypergraphs.

\subsection{The case $p(\alpha)=2$}\label{subsec:COSVF2}
In all this section we assume that $p(\alpha)=2$. Then $\phi$ takes two different values $\varphi_-<\varphi_+$ on $\uuu^{(0)}_\alpha$.
One has $\sss_\alpha=\{S_{\nu^\alpha_+}<S_{\nu^\alpha_-}\}$ with $S_{\nu^\alpha_\pm}=\sigma(\alpha)-\varphi_\pm$.
\subsubsection{The case where  $\alpha$ is of type II}\label{subsec:TCWAIOTII}
The partition \eqref{eq:gradedpartition} takes the form
 $\widehat\uuu^{(0)}_\alpha=\widehat\uuu^{(0)}_{\alpha,+}\sqcup\widehat\uuu^{(0)}_{\alpha,-}$ 
 with $\widehat\uuu^{(0)}_{\alpha,\pm}=\{\m\in\uuu^{(0)}_\alpha,\;\phi(\m)=\varphi_\pm\}$.
 Since $\alpha$ is of type II, then $\m_2(\s)\in \widehat\uuu^{(0)}_\alpha$ for all $\s$. It is then convenient to introduce the partition of $\vvv^{(1)}_\alpha$ given by
\be\label{eq:defV1alphapm}
\vvv^{(1)}_\alpha=\vvv^{(1)}_{\alpha,+}\sqcup\vvv^{(1)}_{\alpha,+-}\cup\vvv^{(1)}_{\alpha,-}
\ee
with 
$
\vvv^{(1)}_{\alpha,+}=\{\s\in \vvv^{(1)}_\alpha,\;\m_1(\s),\m_2(\s)\in \widehat\uuu^{(0)}_{\alpha,+}\}
$
and 
$
\vvv^{(1)}_{\alpha,-}=\{\s\in \vvv^{(1)}_\alpha,\;\m_1(\s),\m_2(\s)\in \widehat\uuu^{(0)}_{\alpha,-}\},
$
where the functions $\m_1,\m_2$ are defined by Lemma \ref{lem:defm1m2}.
 In the case $\s\in \vvv^{(1)}_{\alpha,+-}$, it follows from the choice of Lemma \ref{lem:defm1m2} that 
$\m_1(\s)\in\widehat\uuu^{(0)}_{\alpha,+}$ and $\m_2(\s)\in\widehat \uuu^{(0)}_{\alpha,-}$.
We order the above partitions by deciding $\widehat\uuu^{(0)}_{\alpha,+}\prec\widehat\uuu^{(0)}_{\alpha,-}$ and
$\vvv^{(1)}_{\alpha,+}\prec\vvv^{(1)}_{\alpha,+-}\prec\vvv^{(1)}_{\alpha,-}$. 
Then, the  matrix $\Yr^\alpha:=h^{-\frac 12}e^{-S_{\nu_+^\alpha}/h}\widehat\Lr^\alpha$ has the form
$$
\Yr^\alpha=\left(
\begin{array}{cc}
\iota&0\\
b_{+-}&\tau b_{-+}\\
0&\tau a
\end{array}
\right)
$$
where $\tau=e^{(S_{\nu_+^\alpha}-S_{\nu_-^\alpha})/h}$ and the matrices $\iota,b_{+-},b_{-+}$  admit a classical expansion whose
principal terms are given by the following formula
\begin{itemize}
\item[-] for all $\s \in\vvv^{(1)}_{\alpha,+}$ and $\m\in\widehat\uuu^{(0)}_{\alpha,+}$ one has
$
\iota^0_{\s,\m}=\upsilon_2(\s,\m,\m_1(\s),\m_2(\s))
$
\item[-] 
for all $\s \in\vvv^{(1)}_{\alpha,-}$ and $\m\in\widehat\uuu^{(0)}_{\alpha,-}$ one has
$
a^0_{\s,\m}=\upsilon_2(\s,\m,\m_1(\s),\m_2(\s))
$
\item[-] 
for all $\s \in\vvv^{(1)}_{\alpha,+-}$, $\m\in\widehat\uuu^{(0)}_{\alpha,+}$ and $\m'\in\widehat\uuu^{(0)}_{\alpha,-}$ one has
$
(b^0_{+-})_{\s,\m}=\upsilon_1(\s,\m,\m_1(\s))
$
and 
$
(b^0_{-+})_{\s,\m'}=-\upsilon_1(\s,\m',\m_2(\s))
$
\end{itemize}
 with  $\upsilon_2,\upsilon_1$ given by  \eqref{eq:definupsilon2}, \eqref{eq:definupsilon1}.
By a standard block-matrix computation one has
\begin{equation}\label{eq:blocDalpha}
(\Yr^{\alpha})^*\Yr^{\alpha}=
\left(
\begin{array}{cc}
J&\tau \widehat B\\
\tau \widehat B^*&\tau^2 \widehat  A
\end{array}
\right)
\end{equation}
with $ J=\iota^*\iota+b_{+-}^*b_{+-}$, $\widehat B=b_{+-}^*b_{-+}$ and
$\widehat A=a^*a+b_{-+}^*b_{-+}$. All these matrices admit a classical expansion, $\widehat  A\simeq\sum_{k\geq 0}h^k\widehat  A^k$, 
$\widehat B\simeq\sum_{k\geq 0}h^k\widehat B^k$, $ J=\sum_{k\geq 0}h^k J^k$ and one has 
$ J^0=\iota^{0,*}\iota^0+b_{+-}^{0,*}b_{+-}^0$, $\widehat B^0=b_{+-}^{0,*}b_{-+}^0$ and
$\widehat A^0=a^{0,*}a^0+b_{-+}^{0,*}b^0_{-+}$, where we use the notation $(c^j)^*=c^{j,*}$.

\begin{theorem}\label{thm:spectreDH2ntypeII}
 The matrix $\Lr^\alpha$ has exactly $q_{\alpha,\pm}=\sharp \uuu^{(0)}_{\alpha,\pm}$
singular values $\lambda^\pm_{\alpha,\mu}(h)$, $\mu=1,\ldots , q_{\alpha,\pm}$ counted with multiplicity which are of order
$h^{\frac 12}e^{-S_{\nu^\alpha_\pm}/h}$. 
These singular values  have the following form
$$
\lambda_{\alpha,\mu}^\pm(h)=\zeta^\pm_{\alpha,\mu}(h)e^{-S_{\nu_\pm^\alpha}/h}
$$
where $\zeta^\pm_{\alpha,\mu}\sim h^{\frac 12}\sum_kh^{k}\zeta^\pm_{\alpha,\mu,k}$ is a classical symbol such that $(\zeta^\pm_{\alpha,\mu,0})^2$ are the $q_{\alpha,\pm}$ non-zero eigenvalues 
of the matrices $G^\pm$ given by $G^+=J^0$ and 
$G^-=\widehat A^0-(\widehat B^0)^*(J^0)^{-1}\widehat B^0$, where $\widehat A^0,J^0$ and $\widehat B^0$ are defined below \eqref{eq:blocDalpha}.
\end{theorem}

Let us make a few comments on this theorem. First, observe that the prefactor 
$\zeta^\pm=\zeta^\pm_{\alpha,\mu}$ obeys two different laws wether we are in the $+$ or $-$ case. In the $+$ case, $\zeta^+$ is determined by the matrix 
$J^0$ which depends only on the minima $\m\in \uuu^{(0)}_\alpha$ such that $S(\m)=S_{\nu_+}$. In that sense, the behavior of $\zeta^+$ obeys 
law similar to the generalized Eyring-Kramers law of Theorem \ref{thm:spectreDtypeII}.
In the $-$ case, the situation is different since the matrix $G^-$ involves values of $\phi$ on all minima and not only those for which
 $S(\m)=S_{\nu_-}$. Hence the term $(\widehat B^0)^*(J^0)^{-1}\widehat B^0$ in the definition of $G^-$ can be understood as a tunneling term between minima associate to both  heights.
 
\begin{figure}[!h]
 \center
  \scalebox{0.3}{ 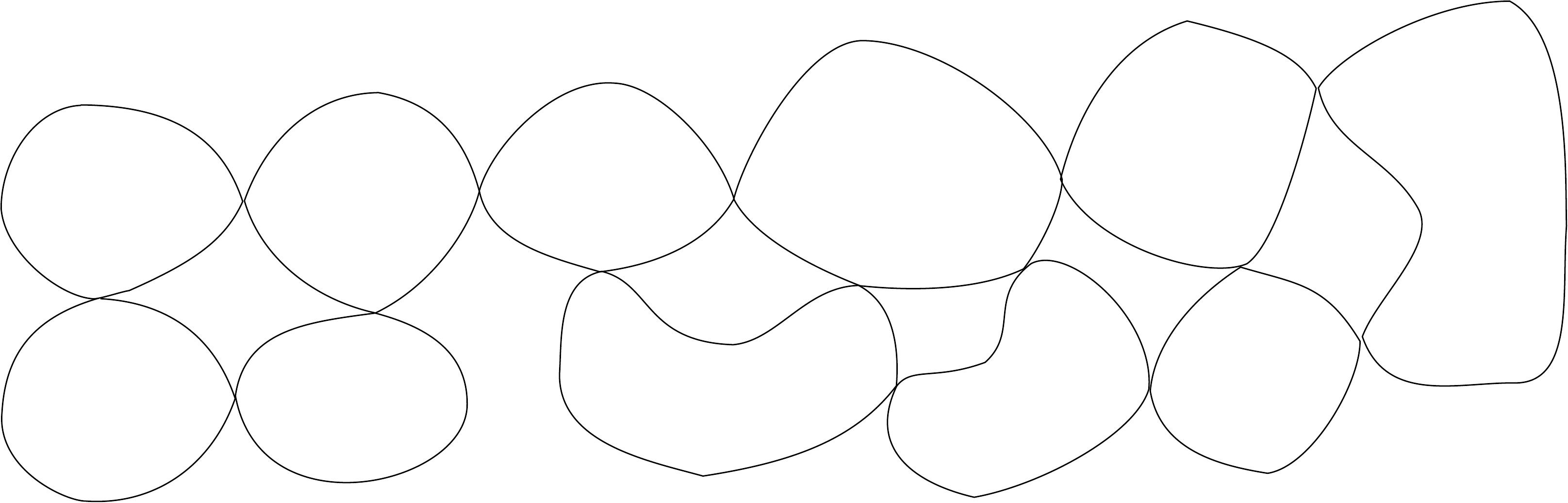 } \\ \vspace{0.2in}
\scalebox{0.5}{ 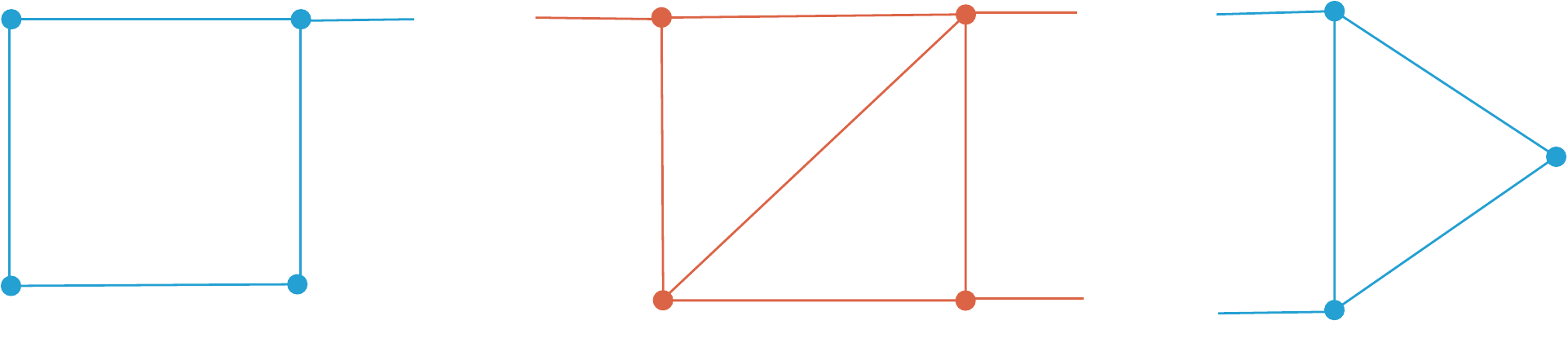}
  \caption{{\em Top:} The sublevel set $\{\phi<\sigma\}$  associated to a potential $\phi$ having a unique saddle value and two minimal values. 
   {\em Bottom: }The associated hypergraphs. }
    \label{figpuits2niveaux}
\end{figure}

 This interpretation is confirmed by the following example. Suppose that $\phi$ has two distinct minimal values and one saddle value.
 Figure \ref{figpuits2niveaux} below represents such a potential. The blue wells correspond to the absolute minimal value and the red one to the other minimal value. All the saddle points are supposed to be at the same level. Then, the matrices $\widehat A^0$ and $J^0$ can be viewed as the Laplacians of the hypergraphs built as follows. First we consider the graph $\ggg$ associated to all the minima whose vertex are the minima and edges are the saddle points between two minima (without distinction on the level of the minima). The blue and red hypergraphs $\ggg_b$ and $\ggg_r$ are obtained by cutting the graph $\ggg$ on edges between a blue and a red minimum.
Eventually, the matrix $B$ links blue  and red minima.
 
\subsubsection{The case where  $\alpha$ is of type I}\label{subsec:COTIP}
In this section we assume that $\alpha$ is of type I. The partition \eqref{eq:gradedpartition} takes the form
$
\uuu^{(0)}_\alpha=\uuu^{(0)}_{\alpha,-}\sqcup\uuu^{(0)}_{\alpha,+}
$
with 
$$
\uuu^{(0)}_{\alpha,\pm}=\{\m\in \uuu^{(0)}_\alpha,\,\phi(\m)=\varphi_\pm\}.
$$
We order the two elements of $\Pr^{(0)}_{\alpha}$ by deciding $\uuu^{(0)}_{\alpha,+}\prec \uuu^{(0)}_{\alpha,-}$.
In order to deal with the saddle points, we introduce the partition 
$\Pr^{(1)}_{\alpha}$ which is a mix of partitions used in Lemma \ref{lem:partitionValpha1} and Section  \ref{subsec:TCWAIOTII}:
$$
\vvv^{(1)}_\alpha=\vvv^{(1)}_{\alpha,+,b}\sqcup\vvv^{(1)}_{\alpha,+,i}\sqcup\vvv^{(1)}_{\alpha,+-}\sqcup\vvv^{(1)}_{\alpha,-,b}\sqcup\vvv^{(1)}_{\alpha,-,i}
$$
with 
\be
\begin{split}
\vvv^{(1)}_{\alpha,+,-}&=\{\s\in \vvv^{(1),i}_\alpha,\;\m_1(\s)\in \uuu^{(0)}_{\alpha,+},\m_2(\s)\in \uuu^{(0)}_{\alpha,-}\}\\
\vvv^{(1)}_{\alpha,+,i}&=\{\s\in \vvv^{(1),i}_\alpha,\;\m_1(\s),\m_2(\s)\in \uuu^{(0)}_{\alpha,+}\}\\
\vvv^{(1)}_{\alpha,+,b}&=\{s\in\vvv^{(1),b}_\alpha,\m_1(\s)\in \uuu^{(0)}_{\alpha,+}\}\\
\vvv^{(1)}_{\alpha,-,i}&=\{\s\in \vvv^{(1),i}_\alpha,\;\m_1(\s),\m_2(\s)\in \uuu^{(0)}_{\alpha,-}\}\\
\vvv^{(1)}_{\alpha,-,b}&=\{s\in\vvv^{(1),b}_\alpha,\m_1(\s)\in \uuu^{(0)}_{\alpha,-}\}.
\end{split}
\ee
Here the function $\m_1,\m_2$ are defined by Lemma \ref{lem:defm1m2}. One has the following
\begin{theorem}\label{thm:spectreDH2ntypeI}
Assume  that $p(\alpha)=2$ and $\alpha$ is of type I. The matrix $\Lr^\alpha$ has exactly $q_{\alpha,\pm}=\sharp \uuu^{(0)}_{\alpha,\pm}$
singular values $\lambda^\pm_{\alpha,\mu}(h)$, $\mu=1,\ldots , q_{\alpha,\pm}$ counted with multiplicity which are of order
$h^{\frac 12}e^{-S_{\nu^\alpha_\pm}/h}$. 
These singular values have the following form
$$
\lambda_{\alpha,\mu}^\pm(h)=\zeta^\pm_{\alpha,\mu}(h)e^{-S_{\nu^\alpha_\pm}/h}
$$
where $\zeta^\pm_{\alpha,\mu}\sim h^{\frac 12} \sum_kh^{k}\zeta^\pm_{\alpha,\mu,k}$ is a classical symbol such that $(\zeta^\pm_{\alpha,\mu,0})^2$ are the $q_{\alpha,\pm}$  eigenvalues  (which  are non zero)
of the matrices $G^\pm$ given by $G^+=J^0$ and 
$G^-=A^0-(B^0)^*(J^0)^{-1}B^0$, where $A^0, B^0$ and $J^0$ are defined by
$$ J^0=\iota^{0,*}\iota^0+b_{+-}^{0,*}b_{+-}^0,\; B^0=b_{+-}^{0,*}b_{-+}^0,\;A^0=a^{0,*}a^0+b_{-+}^{0,*}b^0_{-+}$$
with the matrix, $a^0,b^0_{+-},b^0_{-+}$ and $\iota^0$ defined by
\begin{itemize}
\item[-] for all $\s \in\vvv^{(1)}_{\alpha,+,i}$ and $\m\in\uuu^{(0)}_{\alpha,+}$ one has
$
\iota^0_{\s,\m}=\Upsilon_2(\s,\m,\m_1(\s),\m_2(\s))
$
\item[-]for all $\s \in\vvv^{(1)}_{\alpha,+,b}$ and $\m\in\uuu^{(0)}_{\alpha,+}$ one has
$
\iota^0_{\s,\m}=\Upsilon_1(\s,\m,\m_1(\s))
$

\item[-] 
for all $\s \in\vvv^{(1)}_{\alpha,-,i}$ and $\m\in\uuu^{(0)}_{\alpha,-}$ one has
$
a^0_{\s,\m}=\Upsilon_2(\s,\m,\m_1(\s),\m_2(\s))
$
\item[-] 
for all $\s \in\vvv^{(1)}_{\alpha,-,b}$ and $\m\in\uuu^{(0)}_{\alpha,-}$ one has
$
a^0_{\s,\m}=\Upsilon_1(\s,\m,\m_1(\s))
$
\item[-] 
for all $\s \in\vvv^{(1)}_{\alpha,+-}$, $\m\in\uuu^{(0)}_{\alpha,+}$ and $\m'\in\uuu^{(0)}_{\alpha,-}$ one has
$
(b^0_{+-})_{\s,\m}=\Upsilon_1(\s,\m,\m_1(\s))
$
and 
$
(b^0_{-+})_{\s,\m'}=-\Upsilon_1(\s,\m',\m_2(\s))
$
\end{itemize}
\end{theorem}


 \subsection{Some examples}\label{subsec:SE}

\subsubsection{Computations in dimension one with $p(\alpha)=1$}
\begin{figure}
 \center
  \scalebox{0.5}{ 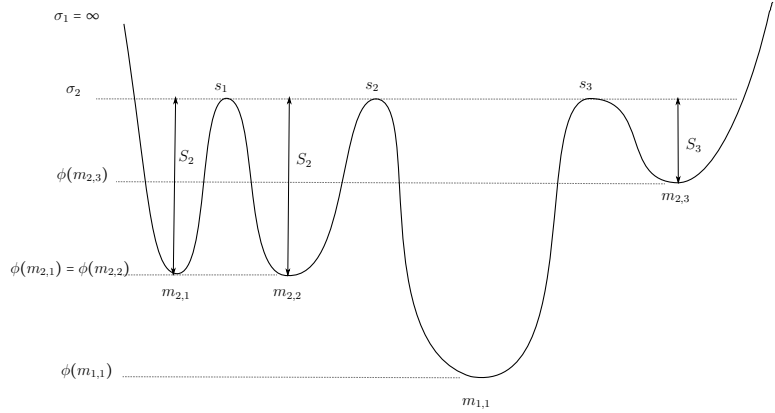}
  \caption{A potential with $p(\alpha)=1$ for all $\alpha$}
    \label{fig6}
\end{figure}

Let us compute the small eigenvalues of the potential $\phi$ represented in Figure \ref{fig6}.

As already noticed in the discussion below Theorem \ref{th:main}, there are exactly three equivalence classes for $\rrr$ in that case:
$\uuu^{(0)}_1=\{\m_{1,1}\}$, $\uuu^{(0)}_2=\{\m_{2,1},\m_{2,2}\}$ and $\uuu^{(0)}_3=\{\m_{2,3}\}$.
Let us denote by $\s_1$ the saddle point between $\m_{2,1}$ and $\m_{2,2}$, $\s_2$ the saddle point between $\m_{2,2}$ and $\m_{1,1}$
and $\s_3$ the saddle point between $\m_{1,1}$ and $\m_{2,3}$.
Denote also $S_2=\phi(\s_1)-\phi(\m_{2,1})=\phi(\s_1)-\phi(\m_{2,2})$ and $S_3=\phi(\s_3)-\phi(\m_{2,3})$. Observe also that for all 
$\m\in\uuu^{(0)}$, one has $H(\m)=\{\m\}$.
Then the matrix $\lll^{bkw}$ defined by \eqref{eq:definLbkw}, admits the following form
$$
\lll^{bkw}=\Big(\frac h \pi\Big)^{\frac 12}
\left(
\begin{array}{cccc}
0&d^2_{1,1}&d^2_{1,2}&0\\
0&d^2_{2,1}&d^2_{2,2}&0\\
0&0&0&d^3
\end{array}
\right)
$$
with the coefficients given by the following formula:
$$
d^2_{1,1}=\big(\vert\phi''(\s_1)\phi''(\m_{2,1})\vert^{\frac 14}+\ooo(h)\big)e^{-S_2/h},\;\;\;
d^2_{1,2}=-\big(\vert\phi''(\s_1)\phi''(\m_{2,2})\vert^{\frac 14}+\ooo(h)\big)e^{-S_2/h}
$$
$$
d^2_{2,1}=0,\;\;\;d^2_{2,2}=\big(\vert\phi''(\s_2)\phi''(\m_{2,2})\vert^{\frac 14}+\ooo(h)\big)e^{-S_2/h},\;\;\;\;
d^3=\big(\vert\phi''(\s_3)\phi''(\m_{2,3})\vert^{\frac 14}+\ooo(h)\big)e^{-S_3/h}
$$
The corresponding squares of singular values are then
$$\lambda_0=0,\;\;\;\lambda_3=\frac h \pi\big(\vert\phi''(\s_3)\phi''(\m_{2,3})\vert^{\frac 12}+\ooo(h)\big)e^{-2S_3/h}\;\;
\text{ and }\;\;\lambda_2^\pm=\frac h\pi\big(\mu_2^\pm +\ooo(h)\big)e^{-2S_2/h}$$
where $\mu_2^\pm$ are the square of the singular values of the matrix 
$$
\widetilde \ddd^2=
\left(\begin{array}{cc}
a&-b\\
0&c
\end{array}
\right)
$$
with $a=\vert\phi''(\s_1)\phi''(\m_{2,1})\vert^{\frac 14}$, $b=\vert\phi''(\s_1)\phi''(\m_{2,2})\vert^{\frac 14}$ and 
$c=\vert\phi''(\s_2)\phi''(\m_{2,2})\vert^{\frac 14}$.
It follows that 
$$
(\widetilde \ddd^2)^*\widetilde \ddd^2=\left(\begin{array}{cc}
a^2&-ab\\
-ab&b^2+c^2
\end{array}
\right)
$$
whose eigenvalues could be computed handily. For instance, if $\vert\phi''(\s)\vert=\vert\phi''(\m)\vert=1$ for all $\s\in\uuu^{(1)}$ and 
$\m\in\uuu^{(0)}$, one has 
$$
(\widetilde \ddd^2)^*\widetilde \ddd^2=\left(\begin{array}{cc}
1&-1\\
-1&2
\end{array}
\right)
$$
whose eigenvalues are
$
\mu_2^\pm=\frac 3 2\pm\frac {\sqrt 5} 2
$

We would like to conclude this example by noticing that 
one has necessarily $\mu_2^+\neq\mu_2^-$. Indeed, if one computes the characteristic polynomial of the above matrix, one finds
$P(x)=x^2-(a^2+b^2+c^2)x+a^2c^2$ whose discriminant is given by 
$$\Delta=(a^2+b^2+c^2)^2-4a^2c^2=((a-c)^2+b^2)((a+c)^2+b^2).$$
Since $\phi$ is  a Morse function, one has $b\neq 0$ and hence $\Delta>0$.

\subsubsection{Computations in dimension one with $p(\alpha)=2$}
Suppose now that the potential $\phi$ is as represented in Figure \ref{fig4}.
As already noticed there are exactly two equivalence classes for $\rrr$ in that case:
$\uuu^{(0)}_1=\{\m_{1,1}\}$, $\uuu^{(0)}_2=\{\m_{2,1},\m_{2,2},\m_{2,3}\}$, and again, one has $H(\m)=\{\m\}$ for all $\m\in\uuu^{(0)}$.
Let us denote by $\s_1$ the saddle point between $\m_{2,1}$ and $\m_{2,2}$, $\s_2$ the saddle point between $\m_{2,2}$ and $\m_{2,3}$
and $\s_3$ the saddle point between $\m_{2,3}$ and $\m_{1,1}$.
Denote also $S_2=\phi(\s_1)-\phi(\m_{2,1})=\phi(\s_1)-\phi(\m_{2,2})$ and $S_3=\phi(\s_2)-\phi(\m_{2,3})$.
Then the matrix $\lll^{bkw,''}$ admits the following form in the basis $(f^{(0)}_{\m_{2,3}},f^{(0)}_{\m_{2,1}},f^{(0)}_{\m_{2,2}})$ and 
$(f^{(1)}_{\s_{3}},f^{(1)}_{\s_{2}},f^{(1)}_{\s_{1}} )$
$$
\lll^{bkw,''}=\Big(\frac h \pi\Big)^{\frac 12} e^{-S_3/h}
\left(
\begin{array}{cccc}
\iota&0&0\\
b_1&0&b_2 e^{-(S_2-S_3)/h}\\
0&a_1e^{-(S_2-S_3)/h}&a_2e^{-(S_2-S_3)/h}
\end{array}
\right)
$$
with 
the leading terms of the coefficients given by the following formula:
$$
\iota^0=-\vert\phi''(\s_3)\phi''(\m_{2,3})\vert^{\frac 14},\;\;\;
b^0_1=\vert\phi''(\s_2)\phi''(\m_{2,3})\vert^{\frac 14},\;\;\;b^0_2= \vert\phi''(\s_2)\phi''(\m_{2,2})\vert^{\frac 14}$$
and
$$
a^0_1=\vert\phi''(\s_1)\phi''(\m_{2,1})\vert^{\frac 14},\;\;\;a^0_2=-\vert\phi''(\s_1)\phi''(\m_{2,2})\vert^{\frac 14}
$$
In order to simplify the computation, assume that $\phi''(\m)=1$ for all $\m\in\uuu^{(0)}$ and
$\phi''(\s_1)=\phi''(\s_2)=-1$. Denote $\theta=\vert \phi''(\s_3)\vert$ and $\tau=e^{-(S_2-S_3)/h}$. Then
$$\lll^{bkw,''}=\Big(\frac h \pi\Big)^{\frac 12} e^{-S_3/h}
\left(\left(
\begin{array}{ccc}
-\theta&0&0\\
1&0&\tau\\
0&\tau&-\tau
\end{array}
\right)+\ooo(h)\right).
$$
Hence, we can apply Theorem \ref{thm:spectreDH2ntypeI} with 
$$
a^0=(1\;-1),\;\iota_0=-\theta,\,b_{+-}^0=1,\;b_{-+}^0=(0\;\,1).
$$
It follows that the singular values of order $e^{-S_2/h}$ are 
$\mu_\pm(h)=\big(\frac h \pi\big)^{\frac 12} e^{-S_2/h}(\sqrt{\lambda_\pm}+\ooo(h))$ with
$\lambda_\pm$ eigenvalues of $M^0:=A^0-(B^0)^*(J^0)^{-1}B^0$ with 
$$
A^0=\left(
\begin{array}{cc}
1&-1\\
-1&2
\end{array}
\right)
,\;
B^0=(0\;-1),\; J^0=1+\theta^2.
$$
Hence 
$$
M^0=\left(
\begin{array}{cc}
1&-1\\
-1&2-\nu
\end{array}
\right)
$$
with $\nu=\frac 1{1+\theta^2}\in]0,1[$. The eigenvalues of this matrix are 
$$
\lambda_\pm=\frac{3-\nu}2\pm\frac{\sqrt{(3-\nu)^2-4(1-\nu)}} 2
$$
This can be seen as perturbations by the well of height $S_3$ of the eigenvalues $\lambda_\pm$ computed in the previous example (obtained by taking $\nu=0$ in the above formula).

\subsubsection{Computations  in higher dimension}
\begin{figure}
 \center
  \scalebox{0.7}{ 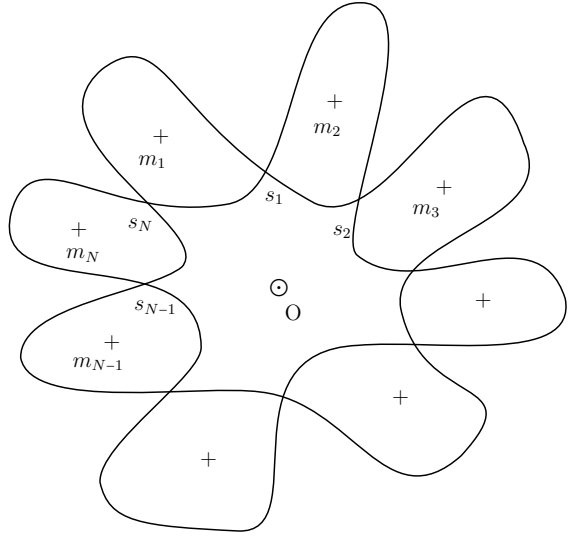}
  \caption{$N$ wells in dimension two}
    \label{fig7}
\end{figure}

Consider the case of potential $\phi$ having $N\geq 3$ minima $\m_1,\ldots,\m_N$ and one local maximum at the origin as presented in Figure \ref{fig7}. Assume also that there are exactly $N$ saddle points $\s_1,\ldots,\s_N$, all at the same height $\phi(\s_j)=\sigma_2$ and that 
the set $\{\phi<\sigma_2\}$ has exactly $N$ connected components $E_1,\ldots, E_N$, each $E_j$ containing  the minimum $m_j$ and that for all 
$j=1,\ldots,N$, $\{\s_j\}=\overline E_j\cap\overline E_{j+1}$ with the convention $E_{N+1}=E_1$.
Assume in addition that all the $\phi(\m_j)$ are equal and denote $S=\sigma_2-\phi(\m_1)$.
Then, Assumption H(1,2) is satisfied. Let us choose $\m_1$ as  the global minimum associated to $\sigma_1=\infty$.
Then all the other minima are associated to the saddle value $\sigma_2$. It is clear that they all belong to the same equivalence class and that they are all of type II. Moreover, for all $\m\in\uuu^{(0)}\setminus\{\m_1\}$, one has $H(\m)=\{\m\}$. Then, we can apply Theorem \ref{thm:spectreDtypeII} to get the spectrum of the Witten Laplacian associated to $\phi$. It follows that the eigenvalues are given by $\lambda_1=0$ and for all $n=2,\ldots,N$
\begin{equation}\label{eq:symexample1}
\lambda_n(h)=b_n(h)e^{-2S/h}(1+\ooo(e^{-\alpha/h}))
\end{equation}
where $b_n$  admit a classical expansion $b_n(h)\simeq\frac h \pi\sum_{k\geq 0}b_{n,k}h^k$. Moreover, one has
$b_{n,0}=\mu_n^2$ where the $\mu_n$, $n=2,\ldots,N$  are the  non zero  singular values of the matrix
$$
\lll:=\left(
\begin{array}{cccccc}
\alpha_1\beta_1&-\alpha_2\beta_1&0&\dots&\ldots&0\\
0&\alpha_2\beta_2&-\alpha_3\beta_2& 0&\ldots&0\\
0&0&\alpha_3\beta_3&\ddots&\ldots&0\\
\vdots&&&&&\vdots\\
0&\ldots&\ldots&0&\alpha_{N-1}\beta_{N-1}&-\alpha_{N}\beta_{N-1}\\
-\alpha_1\beta_N&0&\ldots&\ldots&0&\alpha_N\beta_N
\end{array}
\right)
$$
where we denote $\alpha_j=\varphi''(\m_j)^{\frac 14}$ and $\beta_j=(-\varphi''(\s_j))^{\frac 14}$.

If one assumes additionally that $\alpha_j$ and $\beta_j$ are independent of $j$, let say $\alpha_j=\alpha$ and $\beta_j=\beta$, then 
$\lll=\alpha\beta \aaa$ with 
$$
\Ar=\left(
\begin{array}{ccccccc}
1&-1&0&\ldots&\dots&0&0\\
0&1&-1&0&\ldots&\ldots&0\\
0&0&1&-1&0&\ldots&0\\
\vdots&\vdots&\ddots&\ddots&\ddots&\ddots&\vdots\\
\vdots&\vdots&&\ddots&\ddots&\ddots&0\\
0&0&\ldots&\ldots&0&1&-1\\
-1&0&\ldots&\ldots&0&0&1
\end{array}
\right)
$$
The singular values of $\Ar$ are the square roots of the eigenvalues of 
$$
\Ar^*\Ar=\left(
\begin{array}{cccccc}
2&-1&0&\dots&0&-1\\
-1&2&-1&\ldots&0&0\\
0&-1&2&-1&\ldots&0\\
\vdots&&&&&\vdots\\
0&0&\ldots&-1&2&-1\\
-1&0&0&0&-1&2
\end{array}
\right)
$$
which are known to be $\nu_k=2(1-\cos(\frac{2k\pi}N))$, $k=0,\ldots,N-1$.
In particular, for all $2\leq k<\frac N2$, $\nu_k$ has multiplicity $2$ since 
$\nu_k=\nu_{N-k}$.

Suppose now that the potential $\phi$ is invariant by rotation of angle $\frac {2\pi} N$, then \eqref{eq:symexample1} still holds true with $b_n(h)$ being the singular values of a matrix of  the form
$$
\Ar=\theta(h)
\left(
\begin{array}{ccccccc}
1&-1&0&\ldots&\dots&0&0\\
0&1&-1&0&\ldots&\ldots&0\\
0&0&1&-1&0&\ldots&0\\
\vdots&\vdots&\ddots&\ddots&\ddots&\ddots&\vdots\\
\vdots&\vdots&&\ddots&\ddots&\ddots&0\\
0&0&\ldots&\ldots&0&1&-1\\
-1&0&\ldots&\ldots&0&0&1
\end{array}
\right)
$$
with $\theta(h)\simeq\sum_{k\geq 0}h^{k}\theta_k$. Hence, the above computation is still available and it follows that for  $2\leq k<\frac N2$, 
$b_{k}(h)=b_{N-k}(h)$. This permits to recover the results of \cite{HeHiSj11_01}, section 7.4.

\appendix
\section{Some results in linear algebra}
We collect here some helpful results of linear algebra.
\begin{lemma}\label{lem:Fan}{\bf (Fan inequalities)}
Let $A,B$ be two matrices and denote by $\mu_n$ the singular values of any matrix. Then
$$\mu_n(AB)\leq \Vert B\Vert\mu_n(A)$$
$$\mu_n(AB)\leq \Vert A\Vert\mu_n(B)$$
where $\Vert C\Vert$ denotes the norm of $C:\R^p\rightarrow\R^q$ with $\R^\bullet$ endowed with $\ell^2$ norms.
\end{lemma}
\bp
See \cite{Si79}.
\ep

\begin{lemma}\label{lem:SVdiagmat}
Let $A=\diag(A_1,\ldots,A_N)$ be a block diagonal matrix. Then the singular values of $A$ are the singular values of the $A_n$ counted with multiplicities.
\end{lemma}
\bp
It is straightforward, since 
$A^*A=\diag(A_1^*A_1,\ldots,A_N^*A_N).$
\ep

\begin{lemma}\label{lem:devt-asympt-sv}
Let $E$, $F$ be two finite dimensional vector spaces and $A(h):E\rightarrow F$ be a family of linear operators depending on  a parameter $h\in]0,1]$.
Assume that $A(h)$ admits a classical expansion $A(h)\sim\sum_{k\geq 0}h^k A_k$ and that the matrix $A_0$ has non zero singular values.
Then, for $h>0$ small enough the singular values $\mu_n(h)$ of $A(h)$ admit a classical expansion 
$$\mu_n(h)\sim\sum_{k\geq 0}h^k\mu_n^k$$
where the $\mu_n^0$ are the singular values of $A_0$.
\end{lemma}
\bp
Since the singular values of $A(h)$ are the eigenvalues of $A^*A$ which is selfadjoint, the result follows easily from Kato's perturbation theory
of analytic families of selfadjoint operators
(\cite{Ka66}, chap 2, section 1) applied to the expansion of $A^*A$ in $h$ powers cut at finite rank.
\ep

\begin{lemma}\label{lem:vs-proj}
Let $A$ be a $p\times (q+1)$ matrix and $T$ a $(q+1)\times q$ matrix. Assume that $T^*T=\Id$ and that $\ker A= \Ran(T)^\bot$. Then the singular values of 
$A$ are $\{0,z_1,\ldots,z_q\}$ where $z_1,\ldots,z_q$ are the singular values of $AT$.
\end{lemma}
\bp
First observe that since $\ker A= \Ran(T)^\bot$, then $0$ is a singular value of multiplicity one of $A$. Let us denote $\tilde\xi_0$ a unit vector such that
$\ker A=\R\tilde\xi_0$. 
By definition, there exists an orthonormal basis  $\xi_1,\ldots,\xi_q$ of $\R^q$ such that 
\begin{equation}\label{eq:projSV1}
T^*A^*AT\xi_k=z_k^2\xi_k
\end{equation}
for all $k=1,\ldots,q$. Let us denote $\tilde\xi_k=T\xi_k$. Since $T^*T=Id$, then $\tilde\xi_k$ is an orthonormal family of $\R^{q+1}$. Moreover, since
 $\ker A= \Ran(T)^\bot$, then $\Xi=\{\tilde\xi_0,\ldots,\tilde\xi_q\}$ is an orthonormal basis of $\R^{q+1}$. Moreover, for all $k=1,\ldots,q$, it follows from 
 \eqref{eq:projSV1}
that
 $$
 \vert A\tilde\xi_k\vert^2=\vert AT\xi_k\vert^2=z_k^2
 $$
 This shows that the matrix $A^*A$ in the basis $\Xi$ is exactly $\diag(0,z_1^2,\ldots,z_q^2)$ and proves the result.
\ep
%
%
\begin{lemma}\label{lem:schurpositive}
Let $\mmm$ be a real matrix. Assume that $\mmm$ is symmetric definite positive and that it admits a block decomposition
$$
\mmm=
\left(
\begin{array}{cc}
J&B^*\\
B&N
\end{array}
\right).
$$
Then $J$ and $N-B^*J^{-1}B$ are symmetric definite positive.
\end{lemma}
\bp
This is quite standard, but we recall the proof for  reader's convenience.
Of course $J$ and $N-B^*J^{-1}B$ are symmetric. 
Moreover, since $\mmm$ is positive definite, then 
$$
\<J x,x\>=\<\mmm\binom{x}{0},\binom{x}{0}\>\geq c\vert x\vert^2
$$
for some $c>0$. This shows that $J$ is definite positive.
On the other hand, 
denoting
$$
\Omega=
\left(
\begin{array}{cc}
I&-J^{-1}B^*\\
0&I
\end{array}
\right)
$$
one has
$$
\Omega^*\mmm\Omega=
\left(
\begin{array}{cc}
J&0\\
0&N-BJ^{-1}B^*
\end{array}
\right).
$$
Since $\mmm$ is positive definite, this implies that $N-BJ^{-1}B^*$ is positive definite.
\ep


\section{Link between $\rrr$ and the Generic Assumption}
\begin{proposition} Suppose that the Generic Assumption is satisfied, that is for all $\m\in \uuu^{(0)}$ one has the following:
\begin{itemize}
\item[i)] $\phi_{\vert E(\m)}$ has a unique point of minimum
\item[ii)] if  $E$ is a connected component of $\{\phi<\sigma(\m)\}$ such that $E\cap\vvv^{(1)}\neq\emptyset$, 
there exists a unique $\s\in\vvv^{(1)}$ such that $\phi(\s)=\sup E\cap\vvv^{(1)}$.
 In particular, $E \cap\phi^{-1}(]-\infty,\phi(\s)[)$ is the union of exactly two different connected components.
\end{itemize}
Then for all $\m\in \uuu^{(0)}$,  $\Cl(\m)$ is reduced to $\{\m\}$.
\end{proposition}
\bp If $\m=\underline \m$ there is nothing to prove£. Suppose that $\m\in\ulu^{0)}$ and 
apply Assumption ii) to $E_-(\m)$. One has evidently $\vvv^{(1)}\cap E_-(\m)\neq\emptyset$ since it contains $\overline {E(\m)}\subset E_-(\m)$ and  
$E(\m)$ is a critical component. Hence, $E_-(\m)\cap\{\phi<\sigma(\m)\}$ has exactly two connected components which are necessarily  $\widehat E(\m)$ and $E(\m)$.
Suppose now that $\m'\rrr\m$. Then $\sigma(\m')=\sigma(\m)$ and hence $\m'\notin \widehat E(\m)$. Therefore $\m'\in E(\m)$ which implies 
$\m=\m'$.
\ep
\begin{remark}
There exists some  functions $\phi$ such that $\Cl(\m)=\{\m\}$ for all $\m\in\uuu^{(0)}$ and that do not satisfy the Generic Assumption.
Take for instance $\phi:\R\rightarrow \R$ with $2N+1$ minima and $2N$ saddle points such that 
\begin{itemize}
\item[-] $\m_1<\s_1<\m_2<\s_2<\ldots<\m_{2N}<\s_{2N}<\m_{2N+1}$
\item[-] there exists $\sigma$ such that $\phi(\s_j)=\sigma$ for all $j=1,\ldots,2N$.
\item[-] there exists $\alpha<\beta<\sigma$ such that 
$\forall j=1,\ldots,N$, $\phi(\m_{2j+1})=\beta$, and $\forall j=1,\ldots,N$, $\phi(\m_{2j})=\alpha$.
\end{itemize}
Then, $\Cl(\m_j)=\{\m_j\}$ for all $j$ and since 
$\max\{\phi(\s),\;\s\in \R\}=\sigma$ is equal to $\phi(\s_j)$ for all $j=1,\ldots,2N$, the ii) of (GA) is not satisfied.
\end{remark}
\newpage
\section{List of symbols}
We enumerate different notations used in the paper and give the first place they appear. \\

\begin{tabular}{ccc}
$\uuu^{(0)},\uuu^{(1)}$&page 6\\
$n_0,n_1$& page 6\\
$\Fr(.)$&page 7\\ 
$\vvv^{(1)}$&  Definition \ref{def:SSP} \\
$\Cr,\underline\Sigma,\Sigma$&  Definition \ref{def:SSP} \\
$S,\bsigma$& above \eqref{eq:defSrond}\\
$\sss$& \eqref{eq:defSrond}\\
$\ulu^{(0)}$&\eqref{eq:defulu0}\\
$E$&\eqref{eq:definE}\\
$\Gamma(\m)$& below \eqref{eq:definE}\\
$H(\m)$&\eqref{eq:defHm}\\
$E_-$&\eqref{eq:definEmoins}\\
$\widehat E$&\eqref{eq:defEhat}\\
$\hat\m$&\eqref{eq:defmhat}\\
$\uuu^{(0),I},\uuu^{(0),II}$&  Definition \ref{def:type}\\
$\rrr$& Definition \ref{defin:releq}\\
$\uuu^{(0)}_\alpha$&\eqref{decompU0-1}\\
$\aaa,\ala$& below \eqref{decompU0-1}\\
$q_\alpha$&below \eqref{decompU0-1}\\
$\uuu^{(0),I}_\alpha, \uuu^{(0),II}_\alpha$&below \eqref{decompU0-1}\\
$\sss_\alpha$& \eqref{eq:defSalpha}\\
$p(\alpha)$& \eqref{eq:defSalpha}\\
$\nu_j^\alpha$&  below \eqref{eq:defSalpha}\\
\end{tabular}
\hspace{2cm}
\begin{tabular}{ccc}
$\vvv^{(1)}_\alpha$&\eqref{eq:definV1alpha}\\
$\puu^{(0)}_\alpha$&\eqref {eq:defparenU}\\
$\Gamma_\alpha$&\eqref{eq:defGamalpha}\\
$\vvv^{(1),b}_\alpha,\;\vvv^{(1),i}_\alpha$& Lemma \ref{lem:partitionValpha1}\\
$\widehat\uuu^{(0),II}_\alpha$&\eqref{eq:defhatU0II}\\
$\theta_0^\alpha(\m)$&\eqref{eq:deftheta0}\\
$\widehat H_\alpha(\m)$&\eqref{eq:defhatHm}\\
$\widehat\uuu^{(0)}_\alpha$&\eqref{eq:definhatU0}\\
$\Tr^\alpha$& Definition \ref{definmatT}\\
$\lll^\pi$&\eqref{eq:defLpi}\\
$\lll^{\pi,'}$&below \eqref{eq:defLpi}\\
$\lll^{bkw}$&\eqref{eq:definLbkw}\\
$\lll^{bkw,'}$&\eqref{eq:definLbkw'}\\
$\lll^{bkw,''}$&\eqref{eq:definLr''}\\
$\Lr$&above Lemma \ref{lem:diagblocLr}\\
$\Lr^\alpha$&\eqref{eq:blocdiagD}\\
$\widehat\Lr^\alpha$& Lemma \ref{lem:factorisD}\\
$h_\phi(\m)$& \eqref{eq:defhphi}\\
$\widetilde\Lr^\alpha$&\eqref{eq:deftildeL}\\
$\wideparen\Lr^\alpha$&\eqref{eq:graded1}\\
$\Sr^+,\Sr^+_{cl}$&above Definition \ref{defin:gradedmat}\\
$\Gr(\Er,\tau),\Gr_{cl}(\Er,\tau)$& Definition \ref{defin:gradedmat}\\
$\upsilon_2$& \eqref{eq:definupsilon2}\\
$\upsilon_1$& \eqref{eq:definupsilon1}\\
\end{tabular}
\newpage
\bibliographystyle{amsplain}
\bibliography{ref_witten}

\end{document}

%% file: sublevel2.pdf_tex
\begingroup%
  \makeatletter%
  \providecommand\color[2][]{%
    \errmessage{(Inkscape) Color is used for the text in Inkscape, but the package 'color.sty' is not loaded}%
    \renewcommand\color[2][]{}%
  }%
  \providecommand\transparent[1]{%
    \errmessage{(Inkscape) Transparency is used (non-zero) for the text in Inkscape, but the package 'transparent.sty' is not loaded}%
    \renewcommand\transparent[1]{}%
  }%
  \providecommand\rotatebox[2]{#2}%
  \ifx\svgwidth\undefined%
    \setlength{\unitlength}{482.83986883bp}%
    \ifx\svgscale\undefined%
      \relax%
    \else%
      \setlength{\unitlength}{\unitlength * \real{\svgscale}}%
    \fi%
  \else%
    \setlength{\unitlength}{\svgwidth}%
  \fi%
  \global\let\svgwidth\undefined%
  \global\let\svgscale\undefined%
  \makeatother%
  \begin{picture}(1,0.66465309)%
    \put(0,0){\includegraphics[width=\unitlength,page=1]{sublevel2.pdf}}%
    \put(0.68721414,0.36358853){\color[rgb]{0,0,0}\makebox(0,0)[lb]{\smash{}}}%
    \put(0.82007458,0.12204357){\color[rgb]{0,0,0}\makebox(0,0)[lb]{\smash{}}}%
    \put(0.15598199,0.23321976){\color[rgb]{0,0,0}\makebox(0,0)[lb]{\smash{O}}}%
    \put(0.55473403,0.26434911){\color[rgb]{0,0,0}\makebox(0,0)[lb]{\smash{O}}}%
    \put(0.69555721,0.21246684){\color[rgb]{0,0,0}\makebox(0,0)[lb]{\smash{O}}}%
    \put(0,0){\includegraphics[width=\unitlength,page=2]{sublevel2.pdf}}%
    \put(0.05518223,0.2020904){\color[rgb]{0,0,0}\makebox(0,0)[lb]{\smash{X}}}%
    \put(0.24195828,0.09684364){\color[rgb]{0,0,0}\makebox(0,0)[lb]{\smash{X}}}%
    \put(0.21379364,0.3325372){\color[rgb]{0,0,0}\makebox(0,0)[lb]{\smash{X}}}%
    \put(0.41242849,0.3369842){\color[rgb]{0,0,0}\makebox(0,0)[lb]{\smash{X}}}%
    \put(0.57696927,0.12945526){\color[rgb]{0,0,0}\makebox(0,0)[lb]{\smash{X}}}%
    \put(0.73113362,0.34884301){\color[rgb]{0,0,0}\makebox(0,0)[lb]{\smash{X}}}%
    \put(0.80228634,0.10870239){\color[rgb]{0,0,0}\makebox(0,0)[lb]{\smash{X}}}%
  \end{picture}%
\endgroup%

%% file: graph2.pdf_tex
\begingroup%
  \makeatletter%
  \providecommand\color[2][]{%
    \errmessage{(Inkscape) Color is used for the text in Inkscape, but the package 'color.sty' is not loaded}%
    \renewcommand\color[2][]{}%
  }%
  \providecommand\transparent[1]{%
    \errmessage{(Inkscape) Transparency is used (non-zero) for the text in Inkscape, but the package 'transparent.sty' is not loaded}%
    \renewcommand\transparent[1]{}%
  }%
  \providecommand\rotatebox[2]{#2}%
  \ifx\svgwidth\undefined%
    \setlength{\unitlength}{315.49859545bp}%
    \ifx\svgscale\undefined%
      \relax%
    \else%
      \setlength{\unitlength}{\unitlength * \real{\svgscale}}%
    \fi%
  \else%
    \setlength{\unitlength}{\svgwidth}%
  \fi%
  \global\let\svgwidth\undefined%
  \global\let\svgscale\undefined%
  \makeatother%
  \begin{picture}(1,0.37620459)%
    \put(0,0){\includegraphics[width=\unitlength,page=1]{graph2.pdf}}%
  \end{picture}%
\endgroup%

%% file: sublevel-disjoint.pdf_tex
\begingroup%
  \makeatletter%
  \providecommand\color[2][]{%
    \errmessage{(Inkscape) Color is used for the text in Inkscape, but the package 'color.sty' is not loaded}%
    \renewcommand\color[2][]{}%
  }%
  \providecommand\transparent[1]{%
    \errmessage{(Inkscape) Transparency is used (non-zero) for the text in Inkscape, but the package 'transparent.sty' is not loaded}%
    \renewcommand\transparent[1]{}%
  }%
  \providecommand\rotatebox[2]{#2}%
  \ifx\svgwidth\undefined%
    \setlength{\unitlength}{482.83986883bp}%
    \ifx\svgscale\undefined%
      \relax%
    \else%
      \setlength{\unitlength}{\unitlength * \real{\svgscale}}%
    \fi%
  \else%
    \setlength{\unitlength}{\svgwidth}%
  \fi%
  \global\let\svgwidth\undefined%
  \global\let\svgscale\undefined%
  \makeatother%
  \begin{picture}(1,0.66465309)%
    \put(0,0){\includegraphics[width=\unitlength,page=1]{sublevel-disjoint.pdf}}%
    \put(0.68721414,0.36358853){\color[rgb]{0,0,0}\makebox(0,0)[lb]{\smash{}}}%
    \put(0.82007458,0.12204357){\color[rgb]{0,0,0}\makebox(0,0)[lb]{\smash{}}}%
    \put(0.15598199,0.23321976){\color[rgb]{0,0,0}\makebox(0,0)[lb]{\smash{O}}}%
    \put(0.55473403,0.26434911){\color[rgb]{0,0,0}\makebox(0,0)[lb]{\smash{O}}}%
    \put(0.69555721,0.21246684){\color[rgb]{0,0,0}\makebox(0,0)[lb]{\smash{O}}}%
    \put(0,0){\includegraphics[width=\unitlength,page=2]{sublevel-disjoint.pdf}}%
    \put(0.05518223,0.2020904){\color[rgb]{0,0,0}\makebox(0,0)[lb]{\smash{X}}}%
    \put(0.24195828,0.09684364){\color[rgb]{0,0,0}\makebox(0,0)[lb]{\smash{X}}}%
    \put(0.21379364,0.3325372){\color[rgb]{0,0,0}\makebox(0,0)[lb]{\smash{X}}}%
    \put(0.41242849,0.3369842){\color[rgb]{0,0,0}\makebox(0,0)[lb]{\smash{}}}%
    \put(0.57696927,0.12945526){\color[rgb]{0,0,0}\makebox(0,0)[lb]{\smash{X}}}%
    \put(0.73113362,0.34884301){\color[rgb]{0,0,0}\makebox(0,0)[lb]{\smash{X}}}%
    \put(0.80228634,0.10870239){\color[rgb]{0,0,0}\makebox(0,0)[lb]{\smash{X}}}%
    \put(0.42799314,0.35032538){\color[rgb]{0,0,0}\makebox(0,0)[lb]{\smash{A}}}%
    \put(0.5058165,0.48373681){\color[rgb]{0,0,0}\makebox(0,0)[lb]{\smash{}}}%
  \end{picture}%
\endgroup%

%% file: graph1-disjoint.pdf_tex
\begingroup%
  \makeatletter%
  \providecommand\color[2][]{%
    \errmessage{(Inkscape) Color is used for the text in Inkscape, but the package 'color.sty' is not loaded}%
    \renewcommand\color[2][]{}%
  }%
  \providecommand\transparent[1]{%
    \errmessage{(Inkscape) Transparency is used (non-zero) for the text in Inkscape, but the package 'transparent.sty' is not loaded}%
    \renewcommand\transparent[1]{}%
  }%
  \providecommand\rotatebox[2]{#2}%
  \ifx\svgwidth\undefined%
    \setlength{\unitlength}{313.89620384bp}%
    \ifx\svgscale\undefined%
      \relax%
    \else%
      \setlength{\unitlength}{\unitlength * \real{\svgscale}}%
    \fi%
  \else%
    \setlength{\unitlength}{\svgwidth}%
  \fi%
  \global\let\svgwidth\undefined%
  \global\let\svgscale\undefined%
  \makeatother%
  \begin{picture}(1,0.37812506)%
    \put(0,0){\includegraphics[width=\unitlength,page=1]{graph1-disjoint.pdf}}%
  \end{picture}%
\endgroup%

%% file: figure1.pdf_tex
\begingroup%
  \makeatletter%
  \providecommand\color[2][]{%
    \errmessage{(Inkscape) Color is used for the text in Inkscape, but the package 'color.sty' is not loaded}%
    \renewcommand\color[2][]{}%
  }%
  \providecommand\transparent[1]{%
    \errmessage{(Inkscape) Transparency is used (non-zero) for the text in Inkscape, but the package 'transparent.sty' is not loaded}%
    \renewcommand\transparent[1]{}%
  }%
  \providecommand\rotatebox[2]{#2}%
  \ifx\svgwidth\undefined%
    \setlength{\unitlength}{537.22076635bp}%
    \ifx\svgscale\undefined%
      \relax%
    \else%
      \setlength{\unitlength}{\unitlength * \real{\svgscale}}%
    \fi%
  \else%
    \setlength{\unitlength}{\svgwidth}%
  \fi%
  \global\let\svgwidth\undefined%
  \global\let\svgscale\undefined%
  \makeatother%
  \begin{picture}(1,0.46120302)%
    \put(0,0){\includegraphics[width=\unitlength,page=1]{figure1.pdf}}%
    \put(0.22726922,0.25163307){\color[rgb]{0,0,0}\makebox(0,0)[lb]{\smash{}}}%
    \put(0.22411027,0.24952715){\color[rgb]{0,0,0}\makebox(0,0)[lb]{\smash{}}}%
    \put(0.24411698,0.25637155){\color[rgb]{0,0,0}\makebox(0,0)[lb]{\smash{}}}%
    \put(0.25096138,0.25637155){\color[rgb]{0,0,0}\makebox(0,0)[lb]{\smash{}}}%
    \put(0.25412034,0.26005696){\color[rgb]{0,0,0}\makebox(0,0)[lb]{\smash{}}}%
    \put(0.23095466,0.25373909){\color[rgb]{0,0,0}\makebox(0,0)[lb]{\smash{}}}%
    \put(0.24253751,0.24952715){\color[rgb]{0,0,0}\makebox(0,0)[lb]{\smash{}}}%
    \put(0.55962438,0.45333257){\color[rgb]{0,0,0}\makebox(0,0)[lb]{\smash{$E_{1,1}$}}}%
    \put(0.26886212,0.03524468){\color[rgb]{0,0,0}\makebox(0,0)[lb]{\smash{}}}%
    \put(0.30072607,0.00411604){\color[rgb]{0,0,0}\makebox(0,0)[lb]{\smash{$m_{1,1}$}}}%
    \put(0,0){\includegraphics[width=\unitlength,page=2]{figure1.pdf}}%
    \put(0.58066961,0.32272701){\color[rgb]{0,0,0}\makebox(0,0)[lb]{\smash{$E_{2,1}$}}}%
    \put(0.7015492,0.39919495){\color[rgb]{0,0,0}\makebox(0,0)[lb]{\smash{}}}%
    \put(0.83941058,0.32587486){\color[rgb]{0,0,0}\makebox(0,0)[lb]{\smash{$E_{2,3}$}}}%
    \put(0.52482315,0.59225222){\color[rgb]{0,0,0}\makebox(0,0)[lt]{\begin{minipage}{0.0685007\unitlength}\raggedright \end{minipage}}}%
    \put(0.71603307,0.3252375){\color[rgb]{0,0,0}\makebox(0,0)[lb]{\smash{$E_{2,2}$}}}%
    \put(0.60895988,0.0107408){\color[rgb]{0,0,0}\makebox(0,0)[lb]{\smash{$m_{2,1}$}}}%
    \put(0.71766752,0.24081382){\color[rgb]{0,0,0}\makebox(0,0)[lb]{\smash{$m_{2,2}$}}}%
    \put(0.862205,0.10294009){\color[rgb]{0,0,0}\makebox(0,0)[lb]{\smash{$m_{2,3}$}}}%
    \put(0,0){\includegraphics[width=\unitlength,page=3]{figure1.pdf}}%
    \put(0.40334221,0.23048305){\color[rgb]{0,0,0}\makebox(0,0)[lb]{\smash{$E_{3,1}$}}}%
    \put(0.44901292,0.22930234){\color[rgb]{0,0,0}\makebox(0,0)[lb]{\smash{$E_{3,2}$}}}%
    \put(0.78765737,0.22592236){\color[rgb]{0,0,0}\makebox(0,0)[lb]{\smash{$E_{3,3}$}}}%
    \put(0.20167853,0.11944844){\color[rgb]{0,0,0}\makebox(0,0)[lb]{\smash{$E_{4,1}$}}}%
    \put(0.20837969,0.00403965){\color[rgb]{0,0,0}\makebox(0,0)[lb]{\smash{$m_{4,1}$}}}%
    \put(0.39675662,0.15295422){\color[rgb]{0,0,0}\makebox(0,0)[lb]{\smash{}}}%
    \put(0.45557788,0.15220965){\color[rgb]{0,0,0}\makebox(0,0)[lb]{\smash{$m_{3,2}$}}}%
    \put(0.40494692,0.1507205){\color[rgb]{0,0,0}\makebox(0,0)[lb]{\smash{$m_{3,1}$}}}%
    \put(0.78344542,0.1004637){\color[rgb]{0,0,0}\makebox(0,0)[lb]{\smash{$m_{3,3}$}}}%
    \put(0,0){\includegraphics[width=\unitlength,page=4]{figure1.pdf}}%
    \put(0.0809043,0.30480886){\color[rgb]{0,0,0}\makebox(0,0)[lb]{\smash{$\sigma_2$}}}%
    \put(0,0){\includegraphics[width=\unitlength,page=5]{figure1.pdf}}%
    \put(0.08722221,0.21635812){\color[rgb]{0,0,0}\makebox(0,0)[lb]{\smash{$\sigma_3$}}}%
    \put(0,0){\includegraphics[width=\unitlength,page=6]{figure1.pdf}}%
    \put(0.10365139,0.10728385){\color[rgb]{0,0,0}\makebox(0,0)[lb]{\smash{$\sigma_4$}}}%
    \put(0,0){\includegraphics[width=\unitlength,page=7]{figure1.pdf}}%
    \put(-0.00070846,0.43997417){\color[rgb]{0,0,0}\makebox(0,0)[lb]{\smash{$\sigma_1=\infty$}}}%
    \put(0,0){\includegraphics[width=\unitlength,page=8]{figure1.pdf}}%
    \put(0.91423475,0.22741151){\color[rgb]{0,0,0}\makebox(0,0)[lb]{\smash{$E_{3,4}$}}}%
    \put(0.91721305,0.15667708){\color[rgb]{0,0,0}\makebox(0,0)[lb]{\smash{$m_{3,4}$}}}%
  \end{picture}%
\endgroup%

%% file: figure2.pdf_tex
\begingroup%
  \makeatletter%
  \providecommand\color[2][]{%
    \errmessage{(Inkscape) Color is used for the text in Inkscape, but the package 'color.sty' is not loaded}%
    \renewcommand\color[2][]{}%
  }%
  \providecommand\transparent[1]{%
    \errmessage{(Inkscape) Transparency is used (non-zero) for the text in Inkscape, but the package 'transparent.sty' is not loaded}%
    \renewcommand\transparent[1]{}%
  }%
  \providecommand\rotatebox[2]{#2}%
  \ifx\svgwidth\undefined%
    \setlength{\unitlength}{537.22076635bp}%
    \ifx\svgscale\undefined%
      \relax%
    \else%
      \setlength{\unitlength}{\unitlength * \real{\svgscale}}%
    \fi%
  \else%
    \setlength{\unitlength}{\svgwidth}%
  \fi%
  \global\let\svgwidth\undefined%
  \global\let\svgscale\undefined%
  \makeatother%
  \begin{picture}(1,0.46120302)%
    \put(0,0){\includegraphics[width=\unitlength,page=1]{figure2.pdf}}%
    \put(0.22726922,0.25163307){\color[rgb]{0,0,0}\makebox(0,0)[lb]{\smash{}}}%
    \put(0.22411027,0.24952715){\color[rgb]{0,0,0}\makebox(0,0)[lb]{\smash{}}}%
    \put(0.24411698,0.25637155){\color[rgb]{0,0,0}\makebox(0,0)[lb]{\smash{}}}%
    \put(0.25096138,0.25637155){\color[rgb]{0,0,0}\makebox(0,0)[lb]{\smash{}}}%
    \put(0.25412034,0.26005696){\color[rgb]{0,0,0}\makebox(0,0)[lb]{\smash{}}}%
    \put(0.23095466,0.25373909){\color[rgb]{0,0,0}\makebox(0,0)[lb]{\smash{}}}%
    \put(0.24253751,0.24952715){\color[rgb]{0,0,0}\makebox(0,0)[lb]{\smash{}}}%
    \put(0.55962438,0.45333257){\color[rgb]{0,0,0}\makebox(0,0)[lb]{\smash{$E_{1,1}$}}}%
    \put(0.26886212,0.03524468){\color[rgb]{0,0,0}\makebox(0,0)[lb]{\smash{}}}%
    \put(0.30072607,0.00411604){\color[rgb]{0,0,0}\makebox(0,0)[lb]{\smash{$m_{1,1}$}}}%
    \put(0,0){\includegraphics[width=\unitlength,page=2]{figure2.pdf}}%
    \put(0.58066961,0.32272701){\color[rgb]{0,0,0}\makebox(0,0)[lb]{\smash{$E_{2,1}$}}}%
    \put(0.7015492,0.39919495){\color[rgb]{0,0,0}\makebox(0,0)[lb]{\smash{}}}%
    \put(0.83941058,0.32587486){\color[rgb]{0,0,0}\makebox(0,0)[lb]{\smash{$E_{2,3}$}}}%
    \put(0.52482315,0.59225222){\color[rgb]{0,0,0}\makebox(0,0)[lt]{\begin{minipage}{0.0685007\unitlength}\raggedright \end{minipage}}}%
    \put(0.71603307,0.3252375){\color[rgb]{0,0,0}\makebox(0,0)[lb]{\smash{$E_{2,2}$}}}%
    \put(0.60895988,0.0107408){\color[rgb]{0,0,0}\makebox(0,0)[lb]{\smash{$m_{2,1}$}}}%
    \put(0.71766752,0.24081382){\color[rgb]{0,0,0}\makebox(0,0)[lb]{\smash{$m_{2,2}$}}}%
    \put(0.862205,0.10294009){\color[rgb]{0,0,0}\makebox(0,0)[lb]{\smash{$m_{2,3}$}}}%
    \put(0,0){\includegraphics[width=\unitlength,page=3]{figure2.pdf}}%
    \put(0.40334221,0.23048305){\color[rgb]{0,0,0}\makebox(0,0)[lb]{\smash{$E_{3,1}$}}}%
    \put(0.44901292,0.22930234){\color[rgb]{0,0,0}\makebox(0,0)[lb]{\smash{$E_{3,2}$}}}%
    \put(0.79138023,0.22592236){\color[rgb]{0,0,0}\makebox(0,0)[lb]{\smash{$E_{3,3}$}}}%
    \put(0.20167853,0.11944844){\color[rgb]{0,0,0}\makebox(0,0)[lb]{\smash{$E_{4,1}$}}}%
    \put(0.20837969,0.00403965){\color[rgb]{0,0,0}\makebox(0,0)[lb]{\smash{$m_{4,1}$}}}%
    \put(0.39675662,0.15295422){\color[rgb]{0,0,0}\makebox(0,0)[lb]{\smash{}}}%
    \put(0.45557788,0.15220965){\color[rgb]{0,0,0}\makebox(0,0)[lb]{\smash{$m_{3,2}$}}}%
    \put(0.40494692,0.1507205){\color[rgb]{0,0,0}\makebox(0,0)[lb]{\smash{$m_{3,1}$}}}%
    \put(0.78344542,0.1004637){\color[rgb]{0,0,0}\makebox(0,0)[lb]{\smash{$m_{3,3}$}}}%
    \put(0,0){\includegraphics[width=\unitlength,page=4]{figure2.pdf}}%
    \put(0.0809043,0.30480886){\color[rgb]{0,0,0}\makebox(0,0)[lb]{\smash{$\sigma_2$}}}%
    \put(0,0){\includegraphics[width=\unitlength,page=5]{figure2.pdf}}%
    \put(0.08722221,0.21635812){\color[rgb]{0,0,0}\makebox(0,0)[lb]{\smash{$\sigma_3$}}}%
    \put(0,0){\includegraphics[width=\unitlength,page=6]{figure2.pdf}}%
    \put(0.10365139,0.10728385){\color[rgb]{0,0,0}\makebox(0,0)[lb]{\smash{$\sigma_4$}}}%
    \put(0,0){\includegraphics[width=\unitlength,page=7]{figure2.pdf}}%
    \put(-0.00070846,0.43997417){\color[rgb]{0,0,0}\makebox(0,0)[lb]{\smash{$\sigma_1=\infty$}}}%
    \put(0,0){\includegraphics[width=\unitlength,page=8]{figure2.pdf}}%
    \put(0.91423475,0.22741151){\color[rgb]{0,0,0}\makebox(0,0)[lb]{\smash{$E_{3,4}$}}}%
    \put(0.91721305,0.15667708){\color[rgb]{0,0,0}\makebox(0,0)[lb]{\smash{$m_{3,4}$}}}%
    \put(0,0){\includegraphics[width=\unitlength,page=9]{figure2.pdf}}%
    \put(0.23592888,0.22964522){\color[rgb]{0,0,0}\makebox(0,0)[lb]{\smash{$\widetilde E_{3}$}}}%
    \put(0.30815243,0.32271679){\color[rgb]{0,0,0}\makebox(0,0)[lb]{\smash{$\widetilde E_2$}}}%
    \put(0.85988094,0.22517779){\color[rgb]{0,0,0}\makebox(0,0)[lb]{\smash{$\widetilde E_3'$}}}%
    \put(0,0){\includegraphics[width=\unitlength,page=10]{figure2.pdf}}%
    \put(0.30368499,0.11870387){\color[rgb]{0,0,0}\makebox(0,0)[lb]{\smash{$\widetilde E_4$}}}%
  \end{picture}%
\endgroup%

%% file: figure3.pdf_tex
\begingroup%
  \makeatletter%
  \providecommand\color[2][]{%
    \errmessage{(Inkscape) Color is used for the text in Inkscape, but the package 'color.sty' is not loaded}%
    \renewcommand\color[2][]{}%
  }%
  \providecommand\transparent[1]{%
    \errmessage{(Inkscape) Transparency is used (non-zero) for the text in Inkscape, but the package 'transparent.sty' is not loaded}%
    \renewcommand\transparent[1]{}%
  }%
  \providecommand\rotatebox[2]{#2}%
  \ifx\svgwidth\undefined%
    \setlength{\unitlength}{574.43597208bp}%
    \ifx\svgscale\undefined%
      \relax%
    \else%
      \setlength{\unitlength}{\unitlength * \real{\svgscale}}%
    \fi%
  \else%
    \setlength{\unitlength}{\svgwidth}%
  \fi%
  \global\let\svgwidth\undefined%
  \global\let\svgscale\undefined%
  \makeatother%
  \begin{picture}(1,0.5360791)%
    \put(0,0){\includegraphics[width=\unitlength,page=1]{figure3.pdf}}%
    \put(0.56992409,0.00964064){\color[rgb]{0,0,0}\makebox(0,0)[lb]{\smash{}}}%
    \put(0.56148328,0.00542017){\color[rgb]{0,0,0}\makebox(0,0)[lb]{\smash{$m_{1,1}=\hat m$}}}%
    \put(0.19711956,0.14891479){\color[rgb]{0,0,0}\makebox(0,0)[lb]{\smash{$m_{2,1}$}}}%
    \put(0.34342778,0.14891479){\color[rgb]{0,0,0}\makebox(0,0)[lb]{\smash{$m_{2,2}$}}}%
    \put(0.8526928,0.27552765){\color[rgb]{0,0,0}\makebox(0,0)[lb]{\smash{$m_{2,3}$}}}%
    \put(0.10877837,0.41062951){\color[rgb]{0,0,0}\makebox(0,0)[lb]{\smash{$\sigma_2$}}}%
    \put(0.05605585,0.51209551){\color[rgb]{0,0,0}\makebox(0,0)[lb]{\smash{$\sigma_1=\infty$}}}%
    \put(0,0){\includegraphics[width=\unitlength,page=2]{figure3.pdf}}%
    \put(0.0595643,0.30325177){\color[rgb]{0,0,0}\makebox(0,0)[lb]{\smash{$\phi(m_{2,3})$}}}%
    \put(-0.00124057,0.17986456){\color[rgb]{0,0,0}\makebox(0,0)[lb]{\smash{$\phi(m_{2,1})=\phi(m_{2,2})$}}}%
    \put(0,0){\includegraphics[width=\unitlength,page=3]{figure3.pdf}}%
    \put(0.08981928,0.04604698){\color[rgb]{0,0,0}\makebox(0,0)[lb]{\smash{$\phi(\hat m)$}}}%
    \put(0,0){\includegraphics[width=\unitlength,page=4]{figure3.pdf}}%
    \put(0.19289915,0.41902221){\color[rgb]{0,0,0}\makebox(0,0)[lb]{\smash{$E_{2,1}$}}}%
    \put(0.34905501,0.42464944){\color[rgb]{0,0,0}\makebox(0,0)[lb]{\smash{$E_{2,2}$}}}%
    \put(0.55726281,0.42605627){\color[rgb]{0,0,0}\makebox(0,0)[lb]{\smash{$\widetilde E_2$}}}%
    \put(0.81892938,0.42605627){\color[rgb]{0,0,0}\makebox(0,0)[lb]{\smash{$E_{2,3}$}}}%
    \put(0,0){\includegraphics[width=\unitlength,page=5]{figure3.pdf}}%
  \end{picture}%
\endgroup%

%% file: figure4.pdf_tex
\begingroup%
  \makeatletter%
  \providecommand\color[2][]{%
    \errmessage{(Inkscape) Color is used for the text in Inkscape, but the package 'color.sty' is not loaded}%
    \renewcommand\color[2][]{}%
  }%
  \providecommand\transparent[1]{%
    \errmessage{(Inkscape) Transparency is used (non-zero) for the text in Inkscape, but the package 'transparent.sty' is not loaded}%
    \renewcommand\transparent[1]{}%
  }%
  \providecommand\rotatebox[2]{#2}%
  \ifx\svgwidth\undefined%
    \setlength{\unitlength}{574.43597208bp}%
    \ifx\svgscale\undefined%
      \relax%
    \else%
      \setlength{\unitlength}{\unitlength * \real{\svgscale}}%
    \fi%
  \else%
    \setlength{\unitlength}{\svgwidth}%
  \fi%
  \global\let\svgwidth\undefined%
  \global\let\svgscale\undefined%
  \makeatother%
  \begin{picture}(1,0.5360791)%
    \put(0,0){\includegraphics[width=\unitlength,page=1]{figure4.pdf}}%
    \put(0.56992409,0.00964064){\color[rgb]{0,0,0}\makebox(0,0)[lb]{\smash{}}}%
    \put(0.79501369,0.00542017){\color[rgb]{0,0,0}\makebox(0,0)[lb]{\smash{$m_{1,1}=\hat m$}}}%
    \put(0.19711956,0.14891479){\color[rgb]{0,0,0}\makebox(0,0)[lb]{\smash{$m_{2,1}$}}}%
    \put(0.34342778,0.14891479){\color[rgb]{0,0,0}\makebox(0,0)[lb]{\smash{$m_{2,2}$}}}%
    \put(0.56570366,0.26708679){\color[rgb]{0,0,0}\makebox(0,0)[lb]{\smash{$m_{2,3}$}}}%
    \put(0.10877837,0.41062951){\color[rgb]{0,0,0}\makebox(0,0)[lb]{\smash{$\sigma_2$}}}%
    \put(0.05605585,0.51209551){\color[rgb]{0,0,0}\makebox(0,0)[lb]{\smash{$\sigma_1=\infty$}}}%
    \put(0,0){\includegraphics[width=\unitlength,page=2]{figure4.pdf}}%
    \put(0.0595643,0.30325177){\color[rgb]{0,0,0}\makebox(0,0)[lb]{\smash{$\phi(m_{2,3})$}}}%
    \put(-0.00124057,0.17986456){\color[rgb]{0,0,0}\makebox(0,0)[lb]{\smash{$\phi(m_{2,1})=\phi(m_{2,2})$}}}%
    \put(0,0){\includegraphics[width=\unitlength,page=3]{figure4.pdf}}%
    \put(0.08981928,0.04604698){\color[rgb]{0,0,0}\makebox(0,0)[lb]{\smash{$\phi(\hat m)$}}}%
    \put(0,0){\includegraphics[width=\unitlength,page=4]{figure4.pdf}}%
    \put(0.19289915,0.41902221){\color[rgb]{0,0,0}\makebox(0,0)[lb]{\smash{$E_{2,1}$}}}%
    \put(0.34905501,0.42464944){\color[rgb]{0,0,0}\makebox(0,0)[lb]{\smash{$E_{2,2}$}}}%
    \put(0.81892938,0.42605627){\color[rgb]{0,0,0}\makebox(0,0)[lb]{\smash{$\widetilde E_2$}}}%
    \put(0.55585603,0.42746307){\color[rgb]{0,0,0}\makebox(0,0)[lb]{\smash{$E_{2,3}$}}}%
    \put(0,0){\includegraphics[width=\unitlength,page=5]{figure4.pdf}}%
  \end{picture}%
\endgroup%

%% file: figure5.pdf_tex
\begingroup%
  \makeatletter%
  \providecommand\color[2][]{%
    \errmessage{(Inkscape) Color is used for the text in Inkscape, but the package 'color.sty' is not loaded}%
    \renewcommand\color[2][]{}%
  }%
  \providecommand\transparent[1]{%
    \errmessage{(Inkscape) Transparency is used (non-zero) for the text in Inkscape, but the package 'transparent.sty' is not loaded}%
    \renewcommand\transparent[1]{}%
  }%
  \providecommand\rotatebox[2]{#2}%
  \ifx\svgwidth\undefined%
    \setlength{\unitlength}{614.74972919bp}%
    \ifx\svgscale\undefined%
      \relax%
    \else%
      \setlength{\unitlength}{\unitlength * \real{\svgscale}}%
    \fi%
  \else%
    \setlength{\unitlength}{\svgwidth}%
  \fi%
  \global\let\svgwidth\undefined%
  \global\let\svgscale\undefined%
  \makeatother%
  \begin{picture}(1,0.36367114)%
    \put(0,0){\includegraphics[width=\unitlength,page=1]{figure5.pdf}}%
    \put(0.6217848,-0.12112582){\color[rgb]{0,0,0}\makebox(0,0)[lb]{\smash{}}}%
    \put(0.6250519,0.00506473){\color[rgb]{0,0,0}\makebox(0,0)[lb]{\smash{$m_{1,1}=\hat m$}}}%
    \put(0.27342787,0.00901507){\color[rgb]{0,0,0}\makebox(0,0)[lb]{\smash{$m_{2,1}$}}}%
    \put(0.41014156,0.00901507){\color[rgb]{0,0,0}\makebox(0,0)[lb]{\smash{$m_{2,2}$}}}%
    \put(0.88601025,0.12732497){\color[rgb]{0,0,0}\makebox(0,0)[lb]{\smash{$m_{2,3}$}}}%
    \put(0.19087988,0.25356719){\color[rgb]{0,0,0}\makebox(0,0)[lb]{\smash{$\sigma_2$}}}%
    \put(0.14161476,0.3483793){\color[rgb]{0,0,0}\makebox(0,0)[lb]{\smash{$\sigma_1=\infty$}}}%
    \put(0,0){\includegraphics[width=\unitlength,page=2]{figure5.pdf}}%
    \put(0.14489314,0.15323101){\color[rgb]{0,0,0}\makebox(0,0)[lb]{\smash{$\phi(m_{2,3})$}}}%
    \put(-0.00115922,0.03421714){\color[rgb]{0,0,0}\makebox(0,0)[lb]{\smash{$\phi(m_{1,1}=\phi(m_{2,1})=\phi(m_{2,2})$}}}%
    \put(0,0){\includegraphics[width=\unitlength,page=3]{figure5.pdf}}%
    \put(0.26948422,0.26140952){\color[rgb]{0,0,0}\makebox(0,0)[lb]{\smash{$E_{2,1}$}}}%
    \put(0.41539977,0.26666773){\color[rgb]{0,0,0}\makebox(0,0)[lb]{\smash{$E_{2,2}$}}}%
    \put(0.60995382,0.2679823){\color[rgb]{0,0,0}\makebox(0,0)[lb]{\smash{$\widetilde E_2$}}}%
    \put(0.85446095,0.2679823){\color[rgb]{0,0,0}\makebox(0,0)[lb]{\smash{$E_{2,3}$}}}%
    \put(0,0){\includegraphics[width=\unitlength,page=4]{figure5.pdf}}%
  \end{picture}%
\endgroup%

%% file: puits2niveaux.pdf_tex
\begingroup%
  \makeatletter%
  \providecommand\color[2][]{%
    \errmessage{(Inkscape) Color is used for the text in Inkscape, but the package 'color.sty' is not loaded}%
    \renewcommand\color[2][]{}%
  }%
  \providecommand\transparent[1]{%
    \errmessage{(Inkscape) Transparency is used (non-zero) for the text in Inkscape, but the package 'transparent.sty' is not loaded}%
    \renewcommand\transparent[1]{}%
  }%
  \providecommand\rotatebox[2]{#2}%
  \ifx\svgwidth\undefined%
    \setlength{\unitlength}{930.31989812bp}%
    \ifx\svgscale\undefined%
      \relax%
    \else%
      \setlength{\unitlength}{\unitlength * \real{\svgscale}}%
    \fi%
  \else%
    \setlength{\unitlength}{\svgwidth}%
  \fi%
  \global\let\svgwidth\undefined%
  \global\let\svgscale\undefined%
  \makeatother%
  \begin{picture}(1,0.32047547)%
    \put(0,0){\includegraphics[width=\unitlength,page=1]{puits2niveaux.pdf}}%
    \put(0.06314747,0.06747868){\color[rgb]{0,0,0}\makebox(0,0)[lb]{\smash{}}}%
    \put(0.06314747,0.06226679){\color[rgb]{0,0,0}\makebox(0,0)[lb]{\smash{}}}%
    \put(0.06227873,0.19603882){\color[rgb]{0,0,0}\makebox(0,0)[lb]{\smash{}}}%
    \put(0.05923857,0.19690745){\color[rgb]{0,0,0}\makebox(0,0)[lb]{\smash{}}}%
    \put(0.14561837,0.12750204){\color[rgb]{0,0,0}\makebox(0,0)[lb]{\smash{o}}}%
    \put(0.46624533,0.13978659){\color[rgb]{0,0,0}\makebox(0,0)[lb]{\smash{o}}}%
    \put(0.6001471,0.11030365){\color[rgb]{0,0,0}\makebox(0,0)[lb]{\smash{o}}}%
    \put(0.72913499,0.12750204){\color[rgb]{0,0,0}\makebox(0,0)[lb]{\smash{o}}}%
    \put(0.8409245,0.17786875){\color[rgb]{0,0,0}\makebox(0,0)[lb]{\smash{o}}}%
    \put(0,0){\includegraphics[width=\unitlength,page=2]{puits2niveaux.pdf}}%
    \put(0.06331185,0.18401102){\color[rgb]{0,0,0}\makebox(0,0)[lb]{\smash{x}}}%
    \put(0.06822563,0.05993697){\color[rgb]{0,0,0}\makebox(0,0)[lb]{\smash{x}}}%
    \put(0.22055416,0.06239385){\color[rgb]{0,0,0}\makebox(0,0)[lb]{\smash{x}}}%
    \put(0.22423947,0.19752405){\color[rgb]{0,0,0}\makebox(0,0)[lb]{\smash{x}}}%
    \put(0.37165421,0.20489477){\color[rgb]{0,0,0}\makebox(0,0)[lb]{\smash{x}}}%
    \put(0.56206502,0.20735168){\color[rgb]{0,0,0}\makebox(0,0)[lb]{\smash{x}}}%
    \put(0.74756195,0.2294639){\color[rgb]{0,0,0}\makebox(0,0)[lb]{\smash{x}}}%
    \put(0.92937336,0.24420537){\color[rgb]{0,0,0}\makebox(0,0)[lb]{\smash{x}}}%
    \put(0.78564403,0.07344999){\color[rgb]{0,0,0}\makebox(0,0)[lb]{\smash{x}}}%
    \put(0.45764614,0.06116538){\color[rgb]{0,0,0}\makebox(0,0)[lb]{\smash{xs}}}%
  \end{picture}%
\endgroup%

%% file: hypergraph1.pdf_tex
\begingroup%
  \makeatletter%
  \providecommand\color[2][]{%
    \errmessage{(Inkscape) Color is used for the text in Inkscape, but the package 'color.sty' is not loaded}%
    \renewcommand\color[2][]{}%
  }%
  \providecommand\transparent[1]{%
    \errmessage{(Inkscape) Transparency is used (non-zero) for the text in Inkscape, but the package 'transparent.sty' is not loaded}%
    \renewcommand\transparent[1]{}%
  }%
  \providecommand\rotatebox[2]{#2}%
  \ifx\svgwidth\undefined%
    \setlength{\unitlength}{553.54601766bp}%
    \ifx\svgscale\undefined%
      \relax%
    \else%
      \setlength{\unitlength}{\unitlength * \real{\svgscale}}%
    \fi%
  \else%
    \setlength{\unitlength}{\svgwidth}%
  \fi%
  \global\let\svgwidth\undefined%
  \global\let\svgscale\undefined%
  \makeatother%
  \begin{picture}(1,0.21585932)%
    \put(0,0){\includegraphics[width=\unitlength,page=1]{hypergraph1.pdf}}%
  \end{picture}%
\endgroup%

%% file: figure6.pdf_tex
\begingroup%
  \makeatletter%
  \providecommand\color[2][]{%
    \errmessage{(Inkscape) Color is used for the text in Inkscape, but the package 'color.sty' is not loaded}%
    \renewcommand\color[2][]{}%
  }%
  \providecommand\transparent[1]{%
    \errmessage{(Inkscape) Transparency is used (non-zero) for the text in Inkscape, but the package 'transparent.sty' is not loaded}%
    \renewcommand\transparent[1]{}%
  }%
  \providecommand\rotatebox[2]{#2}%
  \ifx\svgwidth\undefined%
    \setlength{\unitlength}{574.43597208bp}%
    \ifx\svgscale\undefined%
      \relax%
    \else%
      \setlength{\unitlength}{\unitlength * \real{\svgscale}}%
    \fi%
  \else%
    \setlength{\unitlength}{\svgwidth}%
  \fi%
  \global\let\svgwidth\undefined%
  \global\let\svgscale\undefined%
  \makeatother%
  \begin{picture}(1,0.53508432)%
    \put(0,0){\includegraphics[width=\unitlength,page=1]{figure6.pdf}}%
    \put(0.56992409,0.00864587){\color[rgb]{0,0,0}\makebox(0,0)[lb]{\smash{}}}%
    \put(0.5913262,0.00542017){\color[rgb]{0,0,0}\makebox(0,0)[lb]{\smash{$m_{1,1}$}}}%
    \put(0.19711956,0.14792001){\color[rgb]{0,0,0}\makebox(0,0)[lb]{\smash{$m_{2,1}$}}}%
    \put(0.34342778,0.14792001){\color[rgb]{0,0,0}\makebox(0,0)[lb]{\smash{$m_{2,2}$}}}%
    \put(0.8526928,0.27453287){\color[rgb]{0,0,0}\makebox(0,0)[lb]{\smash{$m_{2,3}$}}}%
    \put(0.0714747,0.41460856){\color[rgb]{0,0,0}\makebox(0,0)[lb]{\smash{$\sigma_2$}}}%
    \put(0.05605585,0.51110073){\color[rgb]{0,0,0}\makebox(0,0)[lb]{\smash{$\sigma_1=\infty$}}}%
    \put(0,0){\includegraphics[width=\unitlength,page=2]{figure6.pdf}}%
    \put(0.0595643,0.30225699){\color[rgb]{0,0,0}\makebox(0,0)[lb]{\smash{$\phi(m_{2,3})$}}}%
    \put(-0.00124057,0.17886978){\color[rgb]{0,0,0}\makebox(0,0)[lb]{\smash{$\phi(m_{2,1})=\phi(m_{2,2})$}}}%
    \put(0,0){\includegraphics[width=\unitlength,page=3]{figure6.pdf}}%
    \put(0.06495016,0.0450522){\color[rgb]{0,0,0}\makebox(0,0)[lb]{\smash{$\phi(m_{1,1})$}}}%
    \put(0,0){\includegraphics[width=\unitlength,page=4]{figure6.pdf}}%
    \put(0.26644857,0.41908501){\color[rgb]{0,0,0}\makebox(0,0)[lb]{\smash{$s_1$}}}%
    \put(0.46440674,0.42057715){\color[rgb]{0,0,0}\makebox(0,0)[lb]{\smash{$s_2$}}}%
    \put(0.74493034,0.4220693){\color[rgb]{0,0,0}\makebox(0,0)[lb]{\smash{$s_3$}}}%
    \put(0,0){\includegraphics[width=\unitlength,page=5]{figure6.pdf}}%
    \put(0.88469473,0.34099597){\color[rgb]{0,0,0}\makebox(0,0)[lb]{\smash{$S_3$}}}%
    \put(0,0){\includegraphics[width=\unitlength,page=6]{figure6.pdf}}%
    \put(0.37338577,0.32110069){\color[rgb]{0,0,0}\makebox(0,0)[lb]{\smash{$S_2$}}}%
    \put(0.21919724,0.32507976){\color[rgb]{0,0,0}\makebox(0,0)[lb]{\smash{$S_2$}}}%
  \end{picture}%
\endgroup%

%% file: figure7.pdf_tex
\begingroup%
  \makeatletter%
  \providecommand\color[2][]{%
    \errmessage{(Inkscape) Color is used for the text in Inkscape, but the package 'color.sty' is not loaded}%
    \renewcommand\color[2][]{}%
  }%
  \providecommand\transparent[1]{%
    \errmessage{(Inkscape) Transparency is used (non-zero) for the text in Inkscape, but the package 'transparent.sty' is not loaded}%
    \renewcommand\transparent[1]{}%
  }%
  \providecommand\rotatebox[2]{#2}%
  \ifx\svgwidth\undefined%
    \setlength{\unitlength}{301.01303234bp}%
    \ifx\svgscale\undefined%
      \relax%
    \else%
      \setlength{\unitlength}{\unitlength * \real{\svgscale}}%
    \fi%
  \else%
    \setlength{\unitlength}{\svgwidth}%
  \fi%
  \global\let\svgwidth\undefined%
  \global\let\svgscale\undefined%
  \makeatother%
  \begin{picture}(1,0.95007499)%
    \put(0,0){\includegraphics[width=\unitlength,page=1]{figure7.pdf}}%
    \put(0.2559671,0.69333164){\color[rgb]{0,0,0}\makebox(0,0)[lb]{\smash{}}}%
    \put(0.25217039,0.68194155){\color[rgb]{0,0,0}\makebox(0,0)[lb]{\smash{}}}%
    \put(0.56919508,0.7616724){\color[rgb]{0,0,0}\makebox(0,0)[lb]{\smash{+}}}%
    \put(0.25806698,0.69905896){\color[rgb]{0,0,0}\makebox(0,0)[lb]{\smash{+}}}%
    \put(0.1116924,0.53197181){\color[rgb]{0,0,0}\makebox(0,0)[lb]{\smash{+}}}%
    \put(0.17054129,0.33074644){\color[rgb]{0,0,0}\makebox(0,0)[lb]{\smash{+}}}%
    \put(0.34329128,0.12002943){\color[rgb]{0,0,0}\makebox(0,0)[lb]{\smash{+}}}%
    \put(0.68689294,0.23203225){\color[rgb]{0,0,0}\makebox(0,0)[lb]{\smash{+}}}%
    \put(0.83496431,0.40478224){\color[rgb]{0,0,0}\makebox(0,0)[lb]{\smash{+}}}%
    \put(0.76472533,0.60790583){\color[rgb]{0,0,0}\makebox(0,0)[lb]{\smash{+}}}%
    \put(0.23318689,0.65726292){\color[rgb]{0,0,0}\makebox(0,0)[lb]{\smash{$m_1$}}}%
    \put(0.54641491,0.71459334){\color[rgb]{0,0,0}\makebox(0,0)[lb]{\smash{$m_2$}}}%
    \put(0.72296149,0.56993882){\color[rgb]{0,0,0}\makebox(0,0)[lb]{\smash{$m_3$}}}%
    \put(0.10409898,0.49020805){\color[rgb]{0,0,0}\makebox(0,0)[lb]{\smash{$m_N$}}}%
    \put(0.11548911,0.29657618){\color[rgb]{0,0,0}\makebox(0,0)[lb]{\smash{$m_{N-1}$}}}%
    \put(0.45909073,0.59271907){\color[rgb]{0,0,0}\makebox(0,0)[lb]{\smash{$s_1$}}}%
    \put(0.58058525,0.53197181){\color[rgb]{0,0,0}\makebox(0,0)[lb]{\smash{$s_2$}}}%
    \put(0.22559347,0.39718882){\color[rgb]{0,0,0}\makebox(0,0)[lb]{\smash{$s_{N-1}$}}}%
    \put(0.21420338,0.54526027){\color[rgb]{0,0,0}\makebox(0,0)[lb]{\smash{$s_N$}}}%
    \put(0.46668414,0.42756241){\color[rgb]{0,0,0}\makebox(0,0)[lb]{\smash{$\bigodot$}}}%
    \put(0.49515936,0.38010361){\color[rgb]{0,0,0}\makebox(0,0)[lb]{\smash{O}}}%
  \end{picture}%
\endgroup%